\documentstyle[amscd,amssymb,verbatim,12pt]{amsart}
\theoremstyle{plain}
\newtheorem{Thm}{Theorem}
\newtheorem{Cor}{Corollary}
\newtheorem{Lem}{Lemma}

\theoremstyle{definition}
\newtheorem{Def}{Definition}

\theoremstyle{remark}
\newtheorem{Rem}{Remark}

\newtheorem{Not}{Notation}

\numberwithin{equation}{section}
\newcommand{\SW}{Seiberg-Witten}

\newcommand{\RR}{\Bbb{R}}

\newcommand{\ZZ}{\Bbb{Z}}
\newcommand{\NN}{\Bbb{N}}
\newcommand{\CC}{\Bbb{C}}
\newcommand{\HH}{\Bbb{H}}

\newcommand{\1}{\text{\normalshape{\bf{1}}}}
\newcommand{\ii}{\text{\normalshape{\bf{i}}}}
\newcommand{\jj}{\text{\normalshape{\bf{j}}}}
\newcommand{\kk}{\text{\normalshape{\bf{k}}}}

\newcommand{\SS}{\cal{S}}

\newcommand{\LV}{\Lambda V}
\newcommand{\ilv}{1_{_{ \Lambda V}} }  

\newcommand{\End}{ \operatorname{End} }
\newcommand{\Hom}{ \operatorname{Hom} }

\renewcommand{\exp}{ \operatorname{exp} }
\renewcommand{\ker}{ \operatorname{Ker} }

\renewcommand{\Im}{\operatorname{Im} }
\renewcommand{\Re}{\operatorname{Re} }
\newcommand{\sspan}{ \operatorname{span} }

\newcommand{\cm}{\underset{^{_{c}}}{ \centerdot }}  
\newcommand{\q}{  \operatorname{\bf q} }  
\newcommand{\gr}{ \operatorname{gr} }   
\newcommand{\gri}{ \operatorname{gr}_{i} }   
\newcommand{\tg}{ \underline{\tau} }  
\newcommand{\hh}{ \operatorname{\frak{h}} } 
\newcommand{\wwp}{\text{\normalshape{\bf{w}}}^{+}}  
\newcommand{\wwm}{\text{\normalshape{\bf{w}}}^{-}}  

\newcommand{\lso}{ \operatorname{\frak{so}} }   %
\renewcommand{\o}{ \operatorname{\frak{o}} }   %
\newcommand{\lgl}{ \operatorname{\frak{gl}} }   %
\newcommand{\su}{ \operatorname{SU} }   %
\newcommand{\so}{ \operatorname{SO} }   %
\renewcommand{\sp}{ \operatorname{Sp} }   %
\newcommand{\pin}{ \operatorname{Pin} }   %
\newcommand{\spin}{ \operatorname{Spin} }   %
\renewcommand{\o}{ \operatorname{O} }   %
\newcommand{\gl}{ \operatorname{GL} }   %
\renewcommand{\sl}{ \operatorname{SL} }   %

\newcommand{\la}{\langle \,}   
\newcommand{\ra}{\, \rangle}   

\def\AFourSize{\setlength{\topmargin}{-.5in}           %
	       \setlength{\textheight}{208mm}        %
	       \setlength{\oddsidemargin}{-13mm}       
	       \setlength{\evensidemargin}{-13mm}      %
	       \setlength{\textwidth}{150mm}         %
	       }                                     %
\AFourSize
\begin{document}

\title [Seiberg-Witten's Invariants I]
   {Introduction to\\ Seiberg-Witten's Invariants \\
 Part I \\
 Theory of Spinors}
\author[Jan V. Yang]{Jan Vacter Yang}
\address{Department of Mathematics\\
      Chinese University of Hong Kong\\
       Shatin, Hong Kong}
 \date{ \today}

\maketitle

\bigskip
\begin{center}
{\bf Abstract}
\end{center}
\medskip
 In 1994, Witten has defined a monopole invariant and he has
shown the
equivalence of this invariant with Donaldson's polynomial using his
result in
\( \SS \)-duality. This new invariant is very powerful because
 the gauge
 group  is abelian. By using such an invariant, many new results
are found   in the smooth, K\"{a}hler and even the symplectic
categories.
However, almost every paper in this topic write the monopole
equations in a different way. Therefore is it necessary to clarify the
basic idea behind the definition of such an invariant.

In this paper we  investigate the algebraic structure (Clifford algebra and
\(spin^{c}\) representation ) underlying this invariant and explain the
equations explicitly, especially the K\"{a}hlerian case.
Details of the computations are shown explicitly, and some minute
mistakes
in the existing papers are corrected.

\newpage
\begin{center}
{\bf Forwards}
\end{center}
\medskip
This paper was written by the author while  studying for his Ph.D.
 degree in the
Chinese University of Hong Kong under  Prof. Shing-Tung Yau.

\medskip
 As a fulltime professor at Harvard and   the Director of the Institute of
Mathematical Science at The Chinese University of Hong Kong,
Yau returns to
Hong Kong during the Summer and Winter breaks. In other period,
  the author
 communicate with him through the e-mail. On the other hand,
he also
organize seminars with his colleagues on  topics around
Donaldson's theory.
In the summer of 94,   Yau advised the author to consider monopoles.
Unfortunately, at that time, the author was mainly interested
 in instantons and
 he knew nothing about monopoles. So he ignored such a good advice.

\medskip
In early December 94, (shortly after
Witten has written his preprint), the author received a copy of
Witten's preprint
from Yau in Harvard.
 The paper is so interesting that he started working on it day after
day with his
colleagues. Two weeks later, when Yau arrives at Hong Kong, the
author  thought
 that he had already understood the paper so he started a seminar
with his
colleagues on this topic. However, on the first talk, the author
realised that he
didn't understand it at all ! The most confusing is the fact that
every paper on this
topic appears to write the monopole equations in a unique way.
 ( see \cite{eW94b} ,
 \cite{cT94} , \cite{cT95a} , \cite{KM95} )
As a result, the author spent the whole Christmas and New Year
 vacation with
Yau and his colleagues, trying to clarify everything.  Yau often
asks his students to
 provide all the details explicitely. By doing so, they can learn
 abstract ideas
concretely and correct many wrong misconceptions. After a
 few months' work,
Yau encourage the author and his colleagues to rearrange their
manuscripts and
write down everything into an expository paper.

\medskip
The present paper is the first part of the theory which
includes the basics of
spinors and representation which
is necessary in order to understand Seiberg-Witten's paper. The
 second and third part will be written
by the author's colleagues S. L. Kong and X. Zhang.

\medskip
All the results in this paper had already appeared somewhere else.
 The author's
job is to put them together
logically and consistently, with extra detail, in order to provide
 an introduction
to this new and
exciting field in geometry. This paper may serve as a graduate
course in modern
 geometry for
students with basic knowledge in differential geometry, analysis,
 linear algebra
 and topology.

 \bigskip
\begin{center}
{\bf Acknowledgements}
\end{center}
\medskip
\noindent
This paper cannot be produced without the patience, advice
and encouragement from Prof. Shing-Tung Yau.

\medskip
\noindent
The author  also benefits a lot from
 Sheng Li Kong , Conan Leung, Weiping Li, Kefeng Liu, Li Ma,
Guowu Meng, Xiaowei Peng, Chuen Li Shen, Gang Tian,
 Mu-Tao Wang, Ke Wu, Min Yan,  Yanlin Yu, Xiao Zhang,  Qin Zhou,
and many other good mathematicians.

\tableofcontents
\newpage

\indent
\allowdisplaybreaks

\part{Algebra of Clifford Modules and Spinors}
\bigskip \bigskip
\section{\bf Clifford Algebras and Modules}
To understand spinors and Dirac equations, we must first study
the linear algebra
behind it, which is the theory of Clifford algebra. Here we shall
follow the
approach of \cite{BGV91}. We provide the general theory first
 and then consider
the special case of dimension four.

\bigskip
\subsection{Basic Properties}

\begin{Def}
Let $V$ be a vector space over $\RR$ with a quadratic form $Q$
 on it. The
{\bf Clifford algebra} of $(V, Q)$, denoted by $C(V, Q)$, is the
algebra over $\RR$
generated by $V$ with the relations
\[ v_{1} \cdot v_{2} + v_{2} \cdot v_{1} = - 2 Q( v_{1} , v_{2} )
\qquad\qquad
\forall v_{1} , v_{2} \in V \]
\end{Def}

Since $Q$ is symmetric, we have \( v^2=-Q(v)  \) for all \( v \in V \).

For fixed $Q$, we may abbreviate $C(V, Q)$ and $Q(v_{1}, v_{2})$
 into $C(V)$ and
$(v_{1}, v_{2})$
respectively.

It is a basic fact in algebra that the Clifford algebra is the unique
(up to isomorphism)
solution to the following universal problem.

\begin{Thm}
If $A$ is an algebra and \( c : V \longrightarrow A \) is a linear
 map satisfying
\[ c ( v_{2} ) c ( v_{1} )+c ( v_{1} ) c (v_{2} ) = -2 Q (v_{1} , v_{2} )
\qquad \qquad
\forall  v_{1}, v_{2} \in V , \]
then there is a unique algebra homomorphism from $ C (V, Q)$ to
 $A$ extending
the map $c$. That means we have the following commutative diagram:
\[
 \begin{CD}
   V @>\text{natural}>\text{injection}>  C(V, Q )  \\
    @V{\forall c \text{ linear}}VV      @VV{\exists \,! \text{algebra
 homomorphism}}V \\
    A      @=     A
 \end{CD}
\]
The Clifford algebra may be realized as the quotient \( T(V) /
\cal{I}_{Q} \)
where
\[T(V) = \bigoplus _{k=1}^{\infty} T^{k} (V) \]
is the {\bf tensor algebra} of $V$ with $T^{k}(V)$ generated by
\[ \{\, v_{1} \otimes v_{2} \otimes \dots \otimes v_{k}\, | \,v_{1} ,
 v_{2} , \dots , v_{k}
\in V \,\} \]
and $\cal{I}_{Q}$ is generated by
\[ \{ v_{1} \otimes v_{2} + v_{2} \otimes v_{1} + 2 Q ( v_{1} ,
v_{2} ) \, \mid \, v_{1} ,
v_{2} \in V \}  \]
\end{Thm}

\begin{Rem}
The tensor algebra $T(V)$ has a $\ZZ _{2}$
-grading obtained from
the natural
$\NN$-grading after reduction mod 2:
\[ T(V) = T^{+} (V) + T^{-} (V)  \]
where
\[ T^{+}(V) = \RR \oplus T^{2}(V)
\oplus T^{4}(V) \oplus \dots
\oplus T^{2k}(V)
\oplus \dots \]
\[ T^{-}(V) = V \oplus T^{3}(V)
\oplus T^{5}(V) \oplus \dots
\oplus T^{2k+1}(V)
\oplus \dots \]
Therefore it forms a {\bf superalgebra}.

Similarly, for $k=0, 1, 2, 3, \dots $, let
\[ C^{k}(V) = T^{k}(V) / \cal{I}_{Q} \]
and let
\[ C^{+}(V) = \RR \oplus C^{2}(V)
\oplus C^{4}(V) \oplus \dots
\oplus C^{2k}(V)
\oplus \dots \]
\[ C^{-}(V) =  V  \oplus C^{3}(V) \oplus C^{5}(V) \oplus \dots
\oplus C^{2k+1}(V)
\oplus \dots \]

Since the ideal $\cal{I} _{Q}$ is generated
 by elements from
the evenly graded
subalgebra $T^{+} (V)$, $C(V)$ is itself a
superalgebra and
 we have the grading
\[  C(V) = C ^{+} (V) + C^{-} (V)  . \]
\end{Rem}

\begin{Def}
Let $E$ be a module over $\RR$ or $\CC$
 which is $\ZZ _{2}$
-graded,
\[  E = E ^{+} \oplus E ^{-}   \]
$E$ is called a {\bf Clifford module} over
a Clifford algebra $C(V)$
 if there is a {\bf Clifford action }
\[ \begin{matrix}
 C(V) \times E  & \overset{\cm}{\longrightarrow} &  E  \\
  (\,\, \, a\,\, \,\,  , \,\, \, e \,\,\,) & \longmapsto & a \cm e
\end{matrix} \]
 or equivalently, an algebra homomorphism
\[  \begin{matrix}
C(V) & \overset{c}{\longrightarrow} & \End (E) \\
   a   &    \longmapsto   &  c(a)
\end{matrix} \]
with \[   c(a)\, (e) \,= \,a \cm e \]
which is {\bf even} with respect to this grading:
\[   C ^{+} (V) \cm E ^{\pm} \subset E ^{\pm}  , \]
\[   C ^{-} (V) \cm E ^{\mp} \subset E ^{\mp}  . \]
\end{Def}

\begin{Def}
Let $ O ( V, Q ) $ be the group of linear
 transformations of $V$ which
preserve $Q$.
That means \( \forall \phi \in O(V,Q) ,
 \, \forall v_{1} , v_{2} \in V , \)
\[ Q ( \phi \, v_{1} \,  , \, \phi v_{2} )
 = Q (v_{1}\, , \, v_{2}).  \]

The action of $O(V , Q )$ on generators
 of $T(V)$ are defined by
\[ \phi(v_{1} \otimes v_{2} \otimes
\dots \otimes v_{k} ) =
  \sum _{i=1}^{k} v_{1} \otimes \dots
 \otimes \phi ( v _{i} ) \otimes \dots
\otimes v_{k}  \]
and extends to the whole $T(V)$ linearly.
\end{Def}

\begin{Rem}
$\cal{I}_{Q} $ is invariant under the action
of $O(V , Q )$.
 Hence $C(V, Q )$ carries a natural action
of $O(V,Q)$.
\end{Rem}

\begin{Def}
Let \( ^{\ast} : a \mapsto a^{\ast} \) be the
 anti-automorphism of $T(V)$
induced by
\( v \mapsto -v \) on $T$, and satisfies
\[ ( a_{1} a_{2} ) ^{\ast} \, = \, a_{2}^{\ast}
 \, a_{1}^{\ast} \]
Hence,
\[( v_{1} v_{2} \dots v_{k} )^{\ast} =
		\begin{cases}
		  (v_{k} v_{k-1} \dots v_{1}),
&\text{if $k$ is even,} \\
		 -(v_{k} v_{k-1} \dots v_{1}),
 &\text{if $k$ is odd.}
		  \end{cases}
\]
 \end{Def}

\begin{Rem}
Since
 \begin{align}
			    &v_{1}^{\ast} \otimes v_{2}^{\ast} + v_{2}^{\ast} \otimes v_{1}^{\ast} +
2 Q ( v_{1} ^{\ast} , v_{2} ^{\ast} )  \notag\\
			    =\,&(-v_{1}) \otimes (-v_{2}) + (-v_{2}) \otimes (-v_{1}) + 2 Q
 ( -v_{1} , -v_{2} )
  \notag\\
			    =\,&v_{1} \otimes v_{2} + v_{2} \otimes v_{1} + 2 Q ( v_{1} , v_{2} )
\notag
			    \end{align}
therefore $\cal{I}_{Q}$ is invariant under $ ^{\ast}$. So it induces an
 anti-isomorphism
\( a \mapsto a^{\ast} \) of $C(V)$.
\end{Rem}

\begin{Def}
If $Q$ is a positive-definite quadratic form, then a Clifford module $E$ of
$C(V)$
with an inner product is said to be {\bf self-adjoint} if
\[ c ( a ^{\ast} ) = c ( a) ^{\ast}   \]
 where $c(v)$ denote the action of $v \in V$ on a Clifford module of $C(V)$.
( This module may be $C(V)$ itself.)
\end{Def}

\begin{Rem}
The inner product on $E$ must be {\bf C(V) invariant}:
\[ (\, c(a) e_{1} \, , \, c(a) e_{2} \, ) = ( \, e_{1} \, , \, e_{2} \, )   \]
\( \forall a \in C(V) , \, e_{1} , e_{2} \in E. \)
Hence for a self-dual module $E$,
\[ ( c(a) e_{1}, e_{2} ) = (   e_{1} , c(a)^{\ast} e_{2} ) =
(   e_{1} , c(a^{\ast}) e_{2} ) \]
Especially, we have
\[ ( c(v) e_{1}, e_{2} ) = (   e_{1} , c(v^{\ast}) e_{2} ) =
 (   e_{1} , -c(v) e_{2} ) \]
\( \forall v \in  V  , \,  e_{1} , e_{2} \in E. \)
That means $c(v)$ is {\bf skew-adjoint} $\forall v \in V $.
\end{Rem}

\begin{Def}
Let $E$ be a $\ZZ _{2}$-graded Clifford module over the
 Clifford algebra $C(V)$.
We denote by $\, \End _{C(V)} (E) \,$ the algebra of
 homomorphisms of $E$
supercommuting with the action of $C(V)$.
\end{Def}

\bigskip
\subsection{The Exterior Algebra}
The first interesting example of Clifford module is the
 exterior algebra.

\begin{Def}
The {\bf exterior algebra}
$\LV$ of a vector space $V$  is defined to be
\[   T(V)/\cal{I}_{\Lambda}  \]
where $\cal{I}_{\Lambda}$ is the ideal generated by
elements of the form
\[   ( \,v_{1} \otimes \dots \otimes v_{i} \otimes v_{i+1}
 \otimes \dots
\otimes v_{k}  \, )
 \,+\,( \, v_{1} \otimes \dots \otimes v_{i+1} \otimes v_{i}
\otimes \dots
\otimes v_{k} \, )  \]
or equivalently,
\[ \LV \, = \, \bigoplus_{k=1}^{\infty} \Lambda ^{k} V  \]
where
\[ \Lambda ^{k} V = T ^{k} (V) / \cal{I}_{\Lambda}  \]
\end{Def}

\begin{Def}
 Let
\[ \epsilon : V \longrightarrow \Hom ( \,\Lambda ^{k} V \, , \,
\Lambda ^{k+1} V \,) \]
be the action  of $V$ on $\Lambda V$ by
 {\bf exterior product}, ie.
$\forall v \in V$,
\[ \begin{matrix}
 \epsilon(v) \, :& \Lambda ^{k} V  &  \longrightarrow  &
 \Lambda ^{k+1} V  \\
  & w & \longmapsto & v \wedge w
\end{matrix} \]
Explicitly,
\[ \epsilon (v) \, ( v_{1} \wedge \dots \wedge v_{k} ) \,
   = \, v \wedge v_{1} \wedge \dots \wedge v_{k}  \]
\end{Def}

\begin{Def}
 Let
\[ \iota : V \longrightarrow \Hom ( \,\Lambda ^{k} V \, ,
\,\Lambda ^{k-1} V \,) \]
be the action  of $V$ on $\LV$ by {\bf interior product} or
 {\bf contraction} ,
 ie. $\forall v \in V$,
\[ \begin{matrix}
 \iota(v) \, :& \Lambda ^{k} V  &  \longrightarrow
 & \Lambda ^{k-1} V  \\
  & w & \longmapsto &  Q ( v ,   w )
\end{matrix} \]
Explicitly,
\[ \iota (v) \, ( v_{1} \wedge \dots \wedge v_{k} ) \,
   = \, \sum_{i=1}^{k} ( -1 ) ^{i-1} \, Q (v, v_{i} )  \, v_{1}
\wedge \dots
\wedge \widehat{v_{i}} \wedge \dots  \wedge  v_{k}  \]
\end{Def}

\begin{Def}
The Clifford action of $\, v \in V\,$ on $\, w \in \LV \,
$ is given by
\[  v \cm w \, = \, c(v) \, w \, = \, \epsilon (v) \, w  \, - \,
\iota (v) \, w  \]
\end{Def}

\begin{Lem}
 For any $v , w $ in $V$,
 \[ \epsilon (v) \, \iota(w) \, + \,  \iota(w) \, \epsilon (v) \,
 = \, Q (v \, , \, w )  \]
\end{Lem}

\begin{pf}
For  any generator $ v_{1} \wedge \dots \wedge v_{k}$ of
 $ \Lambda ^{k} V $,
\begin{align}
     & \epsilon (v) \, \iota (w) \, (v_{1} \wedge \dots
 \wedge v_{k}) \notag\\
=\, & \epsilon (v) \, ( \, \sum_{i=1}^{k} ( -1 ) ^{i-1} \,
 Q (w, v_{i} )  \, \,
	v_{1} \wedge \dots \wedge \widehat{v_{i}} \wedge \dots  \wedge  v_{k} \, )
\notag\\
=\, &  \sum_{i=1}^{k} ( -1 ) ^{i-1} \, Q (w, v_{i} ) \, \, w
 \wedge v_{1} \wedge
\dots \wedge \widehat{v_{i}} \wedge \dots  \wedge
  v_{k} \notag\\
\notag\\
     & \iota(w) \, \epsilon (v) \, (v_{1} \wedge \dots \wedge v_{k})
\notag\\
 =\,& \iota(w) \, (  v \wedge v_{1} \wedge \dots  \wedge  v_{k} )
 \notag\\
 =\,& Q (w, v) \, \, v_{1} \wedge \dots \wedge v_{k} +
	 \sum_{i=1}^{k} ( -1 ) ^{i } \, Q (w, v_{i} )  \, \, w \wedge v_{1}
\wedge \dots
\wedge \widehat{v_{i}} \wedge \dots  \wedge  v_{k} \notag
\end{align}
Therefore
\[ (\epsilon (v) \, \iota(w) \, + \,  \iota(w) \, \epsilon (v))
\, ( v_{1} \wedge \dots \wedge v_{k}) \, = \, Q (v \, , \, w ) \,
 ( v_{1} \wedge
 \dots \wedge v_{k})  \]
Since this operation is an algebra homomorphism, the above

 equation holds on
the whole $\LV$.
\end{pf}

\begin{Cor}
The action $c : V \longrightarrow \End(\LV) $ extends to an action
of the Clifford algebra $C(V)$ on $\LV$.
\end{Cor}
\medskip
\begin{pf}
 Observe that
  \begin{align}
	  &  ( \, \epsilon (v) \, \epsilon (w) \, + \, \epsilon (w) \,
\epsilon (v)  \, ) \,
		  ( v_{1} \wedge \dots \wedge v_{k} ) \notag\\
     = \,&  v \wedge  w \wedge  v_{1} \wedge \dots
 \wedge v_{k}   \, +
\, w \wedge v \wedge v_{1} \wedge \dots \wedge
 v_{k} \notag\\
     = \,&   -w \wedge  v \wedge  v_{1} \wedge \dots
\wedge v_{k}   \, +
\, w \wedge v \wedge v_{1} \wedge \dots \wedge v_{k} \notag\\
     = \,&  0 \notag
\end{align}
and
\begin{align}
	  &  ( \, \iota (v) \, \iota (w) \, + \, \iota (w) \, \iota (v)  \, ) \,
		  ( v_{1} \wedge \dots \wedge v_{k} ) \notag\\
      =\,& \iota (v) \, ( \, \sum_{i=1}^{k} ( -1 ) ^{i-1} \,
 Q (w, v_{i} )  \, \,
	      v_{1} \wedge \dots \wedge \widehat{v_{i\,}} \wedge \dots
\wedge  v_{k} \, )  \notag \\
	& + \, \iota (w) \, ( \, \sum_{j=1}^{k} ( -1 ) ^{i-1} \, Q (v, v_{j} )  \, \,
	v_{1} \wedge \dots \wedge \widehat{v_{j}} \wedge \dots
\wedge  v_{k} \, )  \notag\\
      =\,& \sum_{i=1}^{k} (-1)^{i-1} \, \bigl( \, \sum_{j=1}^{i-1}
 ( -1 ) ^{j-1} \,
 Q (w, v_{i} ) \, Q ( v , v_{j} ) \, \,
	      v_{1} \wedge \dots \wedge \widehat{v_{j}} \wedge
\dots \wedge
\widehat{v_{i\,}} \wedge \dots  \wedge  v_{k}   \notag \\
	   & +\, \sum_{j=i+1}^{k} ( -1 ) ^{j-2} \, Q (w, v_{i} ) \,
 Q ( v , v_{j} ) \, \,
	      v_{1} \wedge \dots \wedge \widehat{v_{i\,}} \wedge
 \dots \wedge
 \widehat{v_{j}} \wedge \dots  \wedge  v_{k}   \bigr) \notag \\
	   & + \, \sum_{j=1}^{k} (-1)^{j-1} \, \bigl( \, \sum_{i=1}^{j-1}
 ( -1 ) ^{i-1} \,
 Q (w, v_{j} ) \, Q ( v , v_{i} ) \, \,
	      v_{1} \wedge \dots \wedge \widehat{v_{i\,}} \wedge
 \dots \wedge
\widehat{v_{j}} \wedge \dots  \wedge  v_{k}   \notag \\
	   & +\, \sum_{i=j+1}^{k} ( -1 ) ^{i-2} \, Q (w, v_{j} ) \,
 Q ( v , v_{i} ) \, \,
	      v_{1} \wedge \dots \wedge \widehat{v_{j}} \wedge
\dots \wedge
\widehat{v_{i\,}} \wedge \dots  \wedge  v_{k}   \bigr) \notag \\
    = \, & \, 0  \notag
\end{align}
Hence
\begin{align}
	  &  c(v) \, c(w) \, + \, c(w) \, c(v) \notag\\
    = \, &  ( \, \epsilon (v) - \iota (v) \, ) \, ( \, \epsilon (w)
 - \iota (w) \, )
	       \, + \, ( \, \epsilon (w) - \iota (w) \, ) \, ( \, \epsilon (v)
 - \iota (v) \, )
\notag\\
    = \, &  \epsilon (v) \, \epsilon (w) \, - \, \epsilon (v) \,
\iota (w) \, -
\, \iota (v) \, \epsilon (w) \, + \, \iota (v) \, \iota (w) \notag\\
	  &  + \, \epsilon (w) \, \epsilon (v) \, - \, \epsilon (w) \,
 \iota (v) \, -
 \, \iota (w) \, \epsilon (v) \, + \, \iota (w) \, \iota (v) \notag\\
    = \, &  - \,( \, \epsilon (v) \, \iota (w) \, + \, \iota (w)
 \, \epsilon (v) \, )
\, - \, ( \, \epsilon (w) \, \iota (v) \, + \, \iota (v) \,
\epsilon (w) \, ) \notag\\
    = \, &  -2 \, Q ( v,w) \notag
\end{align}
 by the Lemma above.

As a result, we see that the action $c$ defined above
 extends to a Clifford
action on $\LV$.
\end{pf}

\begin{Def}
The {\bf symbol map} $\sigma : C(V) \longrightarrow
\LV $ is defined by
\[ \sigma (a) \, = \, c(a) \, \ilv \]
here $ \ilv \in \Lambda ^{0} V $ is the identity in the
 exterior algebra $\LV$.
\end{Def}

\begin{Rem}
If $\, 1_{_{C(V)} }\,$ denotes the identity in $C(V)$, then
 $\sigma(1_{_C(V)})$
is the identity
$1_{_{\End (\LV )}}$ in $\End (\LV )$
\end{Rem}

\begin{Not}
Let $\{ e_{i} \}_{i=1, \dots , \dim V }$ be a orthonormal basis of $V$
with respect to the quadratic
form $Q$, ie. $ Q(e_{i} , e_{j} )\, = \, \delta _{ij} $.
Let $c_{i}$ denote the element of $C(V)$ corresponding to $e_{i}$.
\end{Not}

\begin{Thm}
The symbol map $\sigma$ has an inverse $\q : \LV
\longrightarrow C(V)$,
called the {\bf quantization map} , which is given by
\[ \q  ( \, e_{i_{1}} \wedge \dots \wedge e_{i_{k}} \, ) \, = \, c_{i_{1}}
\dots c_{i_{k}}  \]
on the basis of $\LV$ and extends linearly to the whole $\LV$.
\end{Thm}

\begin{pf}
For any generator
\(  e_{i_{1}} \wedge \dots \wedge e_{i_{k}} \)
of $\Lambda V$ , the indices
\( \, i_{j}  \, ( j=1, \dots , k ) \, \) are distinct, hence
\begin{align}
\sigma ( \, \q  ( \, e_{i_{1}} \wedge \dots \wedge e_{i_{k}} \, )  \, ) \,
= \, & \sigma  (  c_{i_{1}} \dots c_{i_{k}} ) \notag\\
= \, & c( c_{i_{1}} \dots c_{i_{k}}  ) \, \ilv \notag\\
= \, & c( c_{i_{1}} ) \dots c ( c_{i_{k}} ) \, \ilv \notag\\
= \, & c( c_{i_{1}} ) \dots c ( c_{i_{k-1}} ) \, (\, \epsilon ( e_{i_{k}} )
 \, \ilv \,
 - \, \iota (e_{i_{k}}) \, \ilv \, ) \notag\\
= \, & c( c_{i_{1}} ) \dots c ( c_{i_{k-1}} ) \, e_{i_{k}} \notag\\
= \, & c( c_{i_{1}} ) \dots c ( c_{i_{k-2}} ) \, ( \, \epsilon ( e_{i_{k-1}} )
 \, e_{i_{k}} \, - \, \iota ( e_{i_{k-1}} ) \, e_{i_{k}} \, ) \notag\\
= \, & c( c_{i_{1}} ) \dots c ( c_{i_{k-2}} ) \, e_{i_{k-1}}
\wedge e_{i_{k}} \notag\\
= \, & \,\, {\bf \cdot \,\, \cdot \,\,\cdot }  \notag \\
= \, & e_{i_{1}} \wedge \dots \wedge e_{i_{k}} \notag
\end{align}
and for the generators \( \{ c_{i} \}_{i=1, \dots , \dim V}\)
 of the algebra $C(V)$,
\begin{align}
   \q ( \sigma ( c_{i} ) ) \, = \, &  \q ( c ( c_{i} ) \, \ilv )\notag\\
    = \, &  \q ( \, \epsilon ( e_{i} ) \, \ilv \, - \, \iota (e_{i})
\, \ilv \, ) \notag\\
    = \, &  \q ( e_{i} ) \notag\\
    = \, &  c_{i} \notag
\end{align}
Hence $\q$ is the inverse of $\sigma$.
\end{pf}

\begin{Cor}
The Clifford algebra $C(V)$ is isomorphic to the tensor algebra
$\Lambda V$
 and   is
therefore  a $ 2^{_{\dim V} } $ dimensional vector
 space with generators
\[ \{ (c_{1})^{n_{1}}  (c_{2})^{n_{2}}  \dots
  (c_{_{\dim V}})^{n_{_{ \dim V}}} \,\, |
 \,\,
    ( n_{1} , n_{2} , \dots , n_{_{ \dim V}} )
 \in \{ 0 , 1 \} ^{_{\dim V}}   \}   \]
\end{Cor}

\begin{Rem}
If we consider $C(V)$ and $\LV$ as $\ZZ_{2}$-graded $O(V)$-modules,
 then $\sigma$ and
$\q$ preserve the $\ZZ_{2}$-grading and the $O(V)$ action.
Hence they are isomorphisms of
$\ZZ_{2}$-graded $O(V)$-modules.
\end{Rem}

\begin{Not}
There is a natural increasing filtration
\[ C_{0}(V) \, \subseteq \, C_{1}(V) \, \subseteq \, \dots \,
\subseteq \, C_{k}(V)
 \, \subseteq \, \dots \, \subseteq \, \bigcup _{i=0}^{\infty}
\, C_{i}(V) \, = \, C(V) \]
where
\[ C_{i}(V) \, = \, C^{0}(V)\, \oplus\, C^{1}(V) \, \oplus\,  C^{2}(v) \,
 \oplus \, \dots \, \oplus \, C^{i}(V) \]
and $C^{0}(V) \, = \, \RR$. It follows that
\[ C_{i}(V) \, = \, \sspan \{ \, v_{1} \dots v_{k} \, | \, v_{j} \in V
\hookrightarrow C(V) \, for \, j = 1 , \dots , k \leq i \,\} \]

The Clifford algebra $C(V)$ with this filtration is called the
 {\bf associated graded algebra} of $C(V)$ and is denoted
by $ \gr C(V)$. The ith grading  of $\gr C(V)$ is denoted by $\gri C(V)$.
\end{Not}

\begin{Rem}
The associated graded algebra $\gr C(V)$ is naturally isomorphic
 to the exterior algebra $\LV$,
where the isomorphism is given by sending $\gri ( v_{1}
\dots v_{i} ) \, \in \gri C(V)$ to
$v_{1} \Lambda \dots \Lambda v_{i} \, \in \Lambda ^{i} V$.
The symbol map $\sigma$ extends the symbol map
\[ \sigma _{i} : C_{i} (V) \longrightarrow \gri C(V) \cong
\Lambda ^{i} V ,\]
in the sense that if $a \in C_{i}(V)$, then $\sigma (a) _{[i]} =
 \sigma _{i} (a) $.
The filtration $C_{i} (V)$ may be written as
\[ C_{i} (V) = \sum_{j=0}^{i} \, \q ( \Lambda ^{j} V ) . \]
Hence the Clifford algabra $C(V)$ may be identified with
the exterior algebra $\LV$ with a
{\bf twisted} , or {\bf quantized} multiplication $\alpha \,
\cdot _{_{Q}} \, \beta $.
\end{Rem}

\begin{Lem} \label{L:sigcom}
If $v \in V \hookrightarrow C(V)$ and $ a \in C^{+}(V) $ , then
\[ \sigma ( \, [ v , a ] \, ) \, = \, -2 \, \iota(v) \, \sigma (a) \]
\end{Lem}

\begin{pf}
Firstly, consider the simple case when
\[ v=e_{i} \qquad\text{ and }\qquad a=c_{i_{1}} \dots c_{i_{k}}=
\q ( \, e_{i_{1}} \wedge \dots \wedge e_{i_{k}} \, ). \]
In this case, we have
\begin{align}
& \sigma  ( \, [ v , a ] \, )  \notag \\
= \, & \sigma ( \, v\, a \, -\,  a \,v \, )  \notag \\
= \, & c(\, v\, a \, -\,  a \,v \, ) \, \ilv \notag \\
= \, & c(\, v\, a\,)\, \ilv \, - \, c(\, a \, v \,) \, \ilv \notag \\
= \, & \epsilon( \, v \, a \, ) \, \ilv\, -\, \iota( \, v \, a \, ) \, \ilv\,
 - \, \epsilon( \, a \,v \, ) \, \ilv\, + \,\iota( \, a \, v \, ) \, \ilv
\notag \\
= \, & \epsilon( \, v \, a \, ) \, \ilv\,  - \, \epsilon( \, a \,v \, ) \, \ilv
\notag \\
= \, & \epsilon( \, c_{i} c_{i_{1}} \dots c_{i_{k}} \, ) \, \ilv\,
- \, \epsilon( \, c_{i_{1}} \dots c_{i_{k}} c_{i}\, ) \, \ilv \notag \\
= \, & e_{i} \wedge e_{i_{1}} \wedge \dots \wedge  e_{i_{k}} \,
- \,  e_{i_{1}} \wedge \dots \wedge  e_{i_{k}} \wedge e_{i} \notag \\
= \, & e_{i} \wedge e_{i_{1}} \wedge \dots \wedge  e_{i_{k}} \,
 - \,  e_{i_{1}} \wedge \dots \wedge  e_{i_{k}} \wedge e_{i} \notag \\
      &  + \, \sum_{j=1}^{k-1} \,  e_{i_{1}} \wedge \dots \wedge e_{i_{j}}
\wedge e_{i} \wedge e_{i_{j+1}} \wedge \dots \wedge e_{i_{k}} \notag \\
      &  - \, \sum_{j=1}^{k-1} \,  e_{i_{1}} \wedge \dots \wedge e_{i_{j}}
 \wedge e_{i} \wedge e_{i_{j+1}} \wedge \dots \wedge e_{i_{k}} \notag \\
= \, & \sum_{j=1}^{k-1} \, (-1)^{j-1} \, \bigl(  \,  e_{i_{1}} \wedge
\dots \wedge ( e_{i} \wedge e_{i_{j}} ) \wedge e_{i_{j+1}} \wedge
\dots \wedge e_{i_{k}} \notag \\
      & \qquad \qquad \,\, \, + \,  e_{i_{1}} \wedge \dots \wedge
( e_{i_{j}} \wedge e_{i} ) \wedge e_{i_{j+1}} \wedge \dots \wedge
 e_{i_{k}}  \,\bigr) \notag \\
      & + \, (-1)^{k-1} \,  e_{i_{1}} \wedge \dots \wedge e_{i_{k-1}}
\wedge e_{i} \wedge e_{i_{k}} \notag \\
      & + \, (-1)^{k-1} \,  e_{i_{1}} \wedge \dots \wedge e_{i_{k-1}}
\wedge e_{i_{k}} \wedge e_{i} \notag \\
      & - \, (-1)^{k-1} \,  e_{i_{1}} \wedge \dots \wedge e_{i_{k-1}}
\wedge e_{i_{k}} \wedge e_{i} \notag \\
      & - \,  e_{i_{1}} \wedge \dots \wedge  e_{i_{k}} \wedge e_{i} \notag \\
= \, & \sum_{j=1}^{k} \, (-1)^{j-1}  \, ( -2 ) \, Q(e_{i} , e_{i_{j}})
\, e_{i_{1}} \wedge \dots \wedge \widehat{ e_{i_{j}} } \wedge
\dots \wedge e_{i_{k}} \notag \\
      &  + \, (\, -1 \, + \, (-1)^{k} \,) \,  e_{i_{1}} \wedge \dots
\wedge  e_{i_{k}} \wedge e_{i} \notag \\
= \, &
	 \begin{cases}
	      (-1)^{j}\, 2 \, e_{i_{1}} \,\wedge \dots \wedge \,
\widehat{ e_{i_{j}} } \,\wedge \dots \wedge \, e_{i_{k}}  &
 \text{ if } i_{j} = i  \notag\\
	      \,\, 0  & \text{ if \( i_{j} \neq i \, \) for \(j=1, \dots , k\)
and $k$ is even} \notag \\
	      -2 \,  e_{i_{1}} \wedge \dots \wedge  e_{i_{k}} \wedge e_{i}
& \text{ if \( i_{j} \neq i \, \) for \(j=1, \dots , k\) and $k$ is odd}
\notag
	 \end{cases}
\end{align}
On the other hand,
\begin{align}
-2 \, \iota(v) \, \sigma(a) \, = \, & -2 \, \iota(v) \, c(a) \, \ilv \notag \\
=\,&-2 \, \iota(v) \, ( \, \epsilon(\,e_{i_{1}} \wedge \dots
\wedge e_{i_{k}}\,)
 \, \ilv \, - \,
  \iota (\,e_{i_{1}} \wedge \dots \wedge e_{i_{k}}\,)\, \ilv \, ) \notag \\
=\,&-2 \, \iota (e_{i}) \, ( \,e_{i_{1}} \wedge \dots \wedge e_{i_{k}}\,)
 \notag \\
=\,&-2 \, \sum_{j=1}^{k} \, Q ( \, e_{i} \, , \, e_{i_{j}} \, ) \, (-1)^{j-1}
\,
       e_{i_{1}} \, \wedge \dots \wedge \, \widehat{ e_{i_{j}} } \,
\wedge \dots
 \wedge \, e_{i_{k}} \notag \\
=\,&
\begin{cases}
 (-1)^{j} \, 2 \, e_{i_{1}} \, \wedge \dots \wedge \,
\widehat{ e_{i_{j}} } \,
\wedge \dots \wedge \, e_{i_{k}}
	  &\text{ if $i_{j} = i$ } \notag\\
\,\, 0 &\text{ if $ i_{j} \neq i \, $ for $j=1, \dots , k$} \notag
\end{cases}
\end{align}
Hence for $a \in C^{+}(V)$ , the above $k$  is even and we have the
required equality.
In general, the equality extends linearly to every $v \in V$
and $a \in C^{+}(V)$.
\end{pf}

\begin{Thm}
The space $C^{2} (V) = \q (\Lambda ^{2} V)$ is a Lie subalgebra of
 $C(V)$, where the Lie bracket
is just the commutator in $C(V)$. It is isomorphic to the Lie algebra
 $\lso (V) $, under the map
\[ \tau : C^{2} (V) \longrightarrow \lso (V)  \]
where any $a \in C^{2}(V) $ is mapped into $\tau (a)$ which acts
on any $v \in V \cong C^{1}(V)$ by
the {\bf adjoint action}:
\[   \tau (a) \, v \, = \, [ \, a \, , \, v \, ] . \]
\end{Thm}
\begin{Rem}
Here the bracket is the bracket of the {\bf Lie superalgebra}
$C(V)$, ie.
\[ [\, a_{1}\, , \, a_{2} \, ] \, = \, a_{1} a_{2} \, - (-1) ^{|a_{1}|
\, |a_{2}|}
\, a_{2}a_{1}  \]
for $a_{1} \, \in C^{|a_{1}|}(V) , \, a_{2} \, \in C^{|a_{2}|}(V)$.

It satisfies the following {\bf Axioms of Lie superalgebra}:
\begin{description}
\item [Supercommutativity]
\( [ \, a_{1} \, , \, a_{2} \, ] \, + \, (-1)^{|a_{1}| \, |a_{2}|} \,
[ \, a_{2} \, ,
 \, a_{1} \, ] \, = \, 0 \)
\item [Jacobi's identity]
	  \( [ \, a_{1} \, , \, [ a_{2} , a_{3} ] \, ] \, = \, [ \,
 [ a_{1}, a_{2}] \, , \, a_{3} \, ] \,
	   + \, (-1)^{|a_{1}| \, |a_{2}| } \, [ \, a_{2} \, , \,
 [a_{1}, a_{3}] \, ]    \)
\end{description}
\end {Rem}

\begin{pf}
Obviously, the action $\tau (a)$ preserves $C^{1}(V) \cong V$.
So it defines a Lie algebra homomorphism from $C^{2}(V)$ to
 $\lgl (V)$.
On the other hand, for any $\, a \, \in C^{2}(V) \,$ and any $v
\, , \, w \,
  \in V$,
\begin{align}
 Q( \, \tau (a) \, v \, , \, w  \, ) \, + \, Q( \, v \, , \, \tau(a) \, w \, )
 \, = \, &
  - \frac{1}{2}\, [\,[ a , v ] \, , \, w\, ] \, - \frac{1}{2} \,
 [ \,v \, , \, [ a,w] \,]
  \notag \\
=\,& -\frac{1}{2} \, [ \, [ v , w] \, , \, a \, ]   \notag
\end{align}
by Jacobi's identity .

Since $a \in C^{2}(V)$, we may  consider  $a=a_{1}  a_{2} $
 for $a_{1} ,
a_{2} \in V$.
In this case,
\begin{align}
[ \, [ v,w]\, , \, a \, ] \, = \, & [ \, (\, vw \, + \, wv \, ) \, ,
\, a_{1}a_{2} \, ] \notag \\
=\,&vwa_{1}a_{2} \, + \, wva_{1}a_{2}\,  -\, a_{1}a_{2}vw \,
- \, a_{1}a_{2}wv  \notag \\
=\,&vwa_{1}a_{2} \, + \, va_{1}wa_{2} \, - \, va_{1}wa_{2} \,
- \, a_{1}vwa_{2} \notag \\
& + \, a_{1}vwa_{2} \, + \, a_{1}va_{2}w \, -  \, a_{1}va_{2}w \,
 - \, a_{1}a_{2}vw \notag \\
&+wva_{1}a_{2} \, + \, wa_{1}va_{2} \, - \, wa_{1}va_{2} \,
- \, a_{1}wva_{2} \notag \\
& + \, a_{1}wva_{2} \, + \, a_{1}wa_{2}v \, -  \, a_{1}wa_{2}v \,
 - \, a_{1}a_{2}wv \notag \\
=\,&-2 \,Q(w,a_{1})\,va_{2}\, +2 \,Q(v,a_{1})\,wa_{2}\,
-2 \,Q(w,a_{2})\,a_{1}v\, +2 \,Q(v,a_{2})\,a_{1}w \notag \\
&-2 \,Q(v,a_{1})\,wa_{2}\, +2 \,Q(w,a_{1})\,va_{2}\,
-2 \,Q(v,a_{2})\,a_{1}w\, +2 \,Q(w,a_{2})\,a_{1}v  \notag \\
=\,&  0 \notag
\end{align}
Hence
\[ Q \, ( \, \tau(a) \, v \, , \, w \, ) \, = \, Q \, ( \, v \, , \, -\tau (a)
 \, w \, )   \]
Therefore
\begin{center}
$ \tau (a) ^{\ast} \, = \, - \, \tau (a)  $
\end{center}
and which means that $\tau$ maps $C^{2}(V)$ into $\lso (V)$.

Now  for linearly independent $ \ a_{1}, \, a_{2} \, \in V$,
 if for any $v \, \in V$,
\[ \tau(a) \, v \, = \, [ \, a_{1} a_{2} \, , \, v \, ] \, = \, 0 , \]
then
\begin{align}
0 \, = \, & a_{1}a_{2}v \, - \, v a_{1}a_{2} \notag \\
=\,&a_{1}a_{2}v \, +\,  a_{1}va_{2} \, - \, a_{1}va_{2} \,
- v a_{1}a_{2} \notag \\
=\,&-2 \, Q (\, v \, , \, a_{2} \, ) \, a_{1} \, +
2 \, Q ( \, v \, , \, a_{1} \, ) \, a_{2} \notag
\end{align}
for any $v \in V$. Which is a contradiction.
Therefore the kernel of $\tau$ contains only the
identity of $C(V)$.

Now both $C^{2}(V)$ and $\lso (V)$ are
$\dim V (\dim V -1) / 2$ dimensional vector
 spaces and the map $\tau$ is injective ,
 therefore $\tau$ is an isomorphism.
\end{pf}

\begin{Lem}
For any $A \, \in \lso (V)$, the corresponding Clifford element is
\[ \tau^{-1} (A)\, = \, \frac{1}{2} \, \sum_{i<j}
\, ( \, A \, e_{i} \, , \, e_{j} \, ) \, c_{i}c_{j}  \]
\end{Lem}
\begin{pf}
Pick an orthonormal basis $\{ e_{i} \}_{i=1, \dots ,  _{\dim V}}$ with
 respect to $Q$.
Any vector in $V$ may be written as $v = \sum_{k=1}^{ _{\dim V}} \,
v^{k} \, e_{k} $.
Any matrix $A \, = \, ( A_{ij})$ in $\lso (V)$ is antisymmetric
 with zero diagonal.
Therefore
\begin{align}
A \, v \, = \, &A \, \sum_{k=1}^{ _{\dim V}} \, v^{k} \, e_{k}  \notag \\
=\,& \sum_{k=1}^{ _{\dim V}} \, v^{k} \, A \, e_{k}  \notag \\
=\,& \sum_{i, \, k=1}^{ _{\dim V}} \, v^{k} \, A_{ik} \, e_{i}  \notag
\end{align}
On the other hand,
\begin{align}
[ \,  \frac{1}{2} \,\sum_{i<j} \, ( \, A \, e_{i} \, , \, e_{j} \, )
\, c_{i}c_{j} \, , \, v \, ] \, = \,
& \frac{1}{2} \,\sum_{i<j} \, ( \, A \, e_{i} \, , \, e_{j} \, )
\, [ \,  c_{i}c_{j} , \, \sum_{k=1}^{ _{\dim V}} \, v^{k} \, c_{k} \, ]
 \notag \\
=\,&\frac{1}{2} \,\sum_{k=1}^{ _{\dim V}} \, v^{k} \,\sum_{i<j}
\, ( \, A \, e_{i} \, , \, e_{j} \, ) \, [ \,  c_{i}c_{j} , \,  c_{k} \, ]
\notag \\
=\,&\frac{1}{2} \,\sum_{k=1}^{ _{\dim V}} \, v^{k} \,\sum_{i<j}
 \,  A_{ji} \, ( \,  c_{i}c_{j} c_{k} \,  - \, c_{k}c_{i}c_{j} \, ) \notag \\
=\,&\frac{1}{2} \,\sum_{k=1}^{ _{\dim V}} \, v^{k} \,\sum_{i<j}
 \,  A_{ji} \, ( \,  c_{i}c_{j} c_{k} \,  + \, c_{i}c_{k}c_{j} \, -
\, c_{i}c_{k}c_{j} \,
- \, c_{k}c_{i}c_{j} \, ) \notag \\
=\,&\frac{1}{2} \,\sum_{k=1}^{ _{\dim V}} \, v^{k} \,\sum_{i<j} \,
 A_{ji} \, ( \,  -2 \, Q ( e_{j}, e_{k} ) \, c_{i} \, + \, 2 \, Q ( e_{i} ,
e_{k} )
\, c_{j}  \, ) \notag \\
=\,&\sum_{k=1}^{ _{\dim V}} \, v^{k} \,\sum_{i<j} \,  A_{ji} \, ( \,
- \, \delta_{jk} \, c_{i} \, +  \, \delta_{ik} \, c_{j}  \, ) \notag \\
=\,& -\,\sum_{i<j} \, v^{j} \, A_{ji}  \, c_{i} \, +   \, \sum_{i<j} \, v^{i}
\, A_{ji} \, c_{j}  \notag \\
=\,& \sum_{i<j} \, v^{j} \, A_{ij}  \, c_{i} \, +   \, \sum_{i<j} \, v^{i}
\, A_{ji} \, c_{j}  \notag \\
=\,& \sum_{i<k} \, v^{k} \, A_{ik}  \, c_{i} \, +   \, \sum_{k<i} \,
 v^{k} \, A_{ik} \, c_{i}  \notag \\
=\,& \sum_{i, \, k=1}^{ _{\dim V}} \, v^{k} \, A_{ik} \, c_{i}  \notag
\end{align}
By identifying $\, e_{i} \, \in V$ with $\, c_{i} \, \in C(V)$, we get
\[ A \, = \, \tau \, ( \, \frac{1}{2} \, \sum_{i<j} \, ( \, A \, e_{i} \, ,
\, e_{j} \, ) \, c_{i}c_{j} \, ) \]
\end{pf}

\begin{Rem}
We usually identify $A \, \in \lso (V)$ with
\[ \sum_{i<j} \, ( \, A\, e_{i} \, , \, e_{j} \, ) \, e_{i} \wedge e_{j}
\qquad \in \Lambda ^{2} V \]
then we have
\[ \q \, ( A ) \, = \, \sum_{i<j} \, ( \, A\, e_{i} \, , \, e_{j} \, )
\, c_{i}c_{j}  \]
which is {\bf twice} of $\tau ^{-1} (A)$.
\end{Rem}

\bigskip
\subsection{The Spin Group}
\begin{Def}
For any $a$ in the Lie algebra $C(V)$, we may form its
{\bf exponential} in $C(V)$ by
\[ \exp \, a \, = \, 1_{_{C(V)}} \, + \, a \, + \, \frac{1}{2}
\, a^{2} \, + \frac{1}{3!} \, a^{3} \, + \,
\dots \, + \, \frac{1}{n!} \, a^{n} \, + \, \dots   \]
which is an element in the associated Lie group of $C(V)$.
\end{Def}

\begin{Thm} \label{T:expsincos}
For any $v_{1}, \, v_{2}$ in $V \hookrightarrow C(V)$ satisfying
\[ Q(\,v_{1}\, , \, v_{1} \, ) \, = \, Q(\,v_{2}\, , \, v_{2} \, ) \, = \, 1 \]
\[ Q(\,v_{1}\, , \, v_{2} \, ) \, = \, 0 \]
we have the following formula
\[ \exp \, t\,  (\,  v_{1} \, v_{2} \, )\, = \, (\cos t) \, 1_{_{C(V)}} \, +
\, ( \sin t ) \, v_{1} \, v_{2} \]
where $t \, \in \RR$. In fact, $t$ is well-defined $\mod 2\pi$.

Consequently, this formula is satisfied for some vectors in $V$
whenever $\, \dim V > 1$.
\end{Thm}

\begin{pf}
Since
\begin{align}
( \, v_{1} \, v_{2} \, ) ^{2} \, = \,& v_{1} \, v_{2} \, v_{1} \, v_{2}
 \notag \\
=\,&v_{1} \, v_{2} \, v_{1} \, v_{2} \, + \, v_{1} \, v_{1} \, v_{2}
 \, v_{2} \,
 - \, v_{1} \, v_{1} \, v_{2} \, v_{2}  \notag \\
=\,&-2 \, Q ( \, v_{1} \, , \, v_{2} \, ) \, v_{1} \, v_{2} \, - (-1)
\, Q(\, v_{1}
\, , \, v_{1} \, ) \, (-1) \, Q( \, v_{2} \, , \, v_{2} \, ) \notag \\
=\,&-1 \notag
\end{align}
Therefore
\begin{align}
 \exp \, t\,  (\,  v_{1} \, v_{2} \, )\, =\, & \sum_{k=0}^{\infty} \,
\frac{1}{k!} \, t^{k} \, (\,  v_{1} \, v_{2} \, )^{k} \notag \\
=\,&\sum_{k=0}^{\infty} \, \frac{1}{(2k)!} \, t^{2k} \, (\,  v_{1} \,
 v_{2} \, )^{2k} \, + \, \sum_{k=0}^{\infty} \, \frac{1}{(2k+1)!} \,
 t^{2k+1} \, (\,  v_{1} \, v_{2} \, )^{2k+1} \notag \\
=\,&\sum_{k=0}^{\infty} \, \frac{1}{(2k)!} \, t^{2k} \, (-1)^{k} \,
+ \, \sum_{k=0}^{\infty} \, \frac{1}{(2k+1)!} \, t^{2k+1} \, (-1)^{k}
\, (\,  v_{1} \, v_{2} \, ) \notag \\
=\,&(\cos t) \, 1_{_{C(V)}} \, + \, (\sin t) \, (\,  v_{1} \, v_{2} \, ) \notag
\end{align}
\end{pf}

\begin{Def}
Let the {\bf Spin group} of the vector space $V$ be the Lie
 group associated
to the Lie subalgebra
$C^{2}(V)$ of the Clifford algebra $C(V)$, ie.
\[ \spin (V) \, = \, \exp C^{2}(V) \]
\end{Def}

\begin{Rem}
The adjoint action $\tau$ of the Lie algebra
$C^{2}(V)$ on $V$ may be
exponentiated to an orthogonal
action of {\bf conjugation} which is  denoted by
 $\tg$ . Explicitly, for $g
\, \in \spin (V)$ and $v \, \in V$, there is a fundamental relation
\[   \tg (g) \, v  \, = \, g \, v \, g^{-1}\]

Indeed, writing $g=\exp (a)$ for some $a \in C^{2}(V)$,
then \[ [ \, a \, , \, v \, ] \, = \, \tau (a) \, v \]
implies
\[ \exp (a) \, v \, (\exp (a) )^{-1} \, = \, \exp ( \, \tau(a) \, ) \, v \]
In other words, we have the following commutative diagram:
\[
 \begin{CD}
   C^{2}(V) @> \tau >>  \lso (V)  \\
    @V \exp VV      @VV \exp V \\
    \spin (V)      @>> \tg >     \so (V)
 \end{CD}
\]
\end{Rem}

\begin{Thm}
If $\, \dim V > 1 $, then the homomorphism
\[ \tg : \spin (V) \longrightarrow \so (V)  \]
is a double covering.
\end{Thm}

\begin{pf}
Since $\tau$ is an isomorphism and the exponential map
 is surjective, therefore
$\tg$ is also surjective.
Pick any $\, g \, \in \ker ( \tg )$. Then $ \, g \, = \,
\exp (a)$ for some $\, a \,
\in C^{2}(V)$.
Now \( \tg \, ( g ) \, = \, 1_{_{\so( V)}} \) implies
\( [ \, a \, , \, v \, ] \, = \, 0 \) for all \( v \in V \).
By Lemma \ref{L:sigcom} , we have
\( \iota (v) \, \sigma (a) \, = \, 0 \)
for all $v \, \in V$. Therefore $a \, \in \RR$ and is
identified with a scalar
multiplication in $\lso (V)$.
So $g$ is also a scalar multiplication.

Now apply the anti-automorphism \( a \mapsto
 a^{\ast}\) on $C^{2}(V)$.
For orthogonal vectors $v_{1} , \, v_{2} \, \in V$,
\[ ( \, v_{1} \, , \, v_{2} \, )^{\ast} \, = \, (v_{2})^{\ast} \,
(v_{1})^{\ast} \, =
\, ( -v_{2}) \, (-v_{1}) \,
= \, v_{2} \, v_{1} \, = \, - v_{1} \, v_{2} \]
$C^{2}(V)$ may be generated ( as a vector space)
by pairs of orthogonal
vectors, therefore
\[ a^{\ast} \, = \, - a \qquad \qquad \qquad \forall a \, \in C^{2}(V) \]
Now we have
\[ (\exp (a) )^{\ast} \, = \, \exp ( a ^{\ast}) \, = \, \exp (-a)  \]
which implies that for any $g \, \in \spin (V)$,
\[ g \cdot g^{\ast} \, = \, 1_{_{\spin (V)}}    \]
Hence any $g \, \in \ker ( \tg)$ is a scalar with $g \cdot g^{\ast} = 1$
which implies $g \, = \, \pm 1$.

For $\, \dim V > 1$, we may apply Theorem \ref{T:expsincos}
with $t =
\pi$ and hence $-1$ is in $\spin (V)$.
As a result, $\ker (\tg) \, = \, \ZZ_{2}$ and we have a double cover
\[ 0 \longrightarrow \ZZ_{2} \longrightarrow \spin (V) \overset{\tg}
{\longrightarrow}
\so (V) \longrightarrow 0 \]
\end{pf}

\bigskip
\subsection{Four Dimensional Case}
Now consider the most interesting case when \( V \cong \RR ^{4} \).
Fix a basis $\{ \, e_{1} , e_{2} , e_{3} , e_{4}\, \}$ which is orthonormal
with respect to the
fixed quadratic form $Q$.
Then any vector $v \, \in V$ may be written as
\[ v \, = \, \sum_{k=1}^{4} \, v^{k} \, e_{k}   \]
and the Clifford algebra is
\[ C(\RR ^{4}) \, = \, \sspan \{ \, c_{1} ^{n_{1}} c_{2} ^{n_{2}} c_{3}
 ^{n_{3}}
 c_{4} ^{n_{4}} \, | \,
   n_{i} \in  \{ 0, 1 \} \, \}     \]
where $c_{i} \, = \, \q ( \, e_{i} \, ) $

\begin{Not}
For convenience, we may denote $c_{i}$ by $e_{i}$ without ambiguity.
\end{Not}

Especially,
\[ C^{2}(\RR ^{4}) \, = \, \sspan \{ \,e_{1}e_{2}\,, \,  e_{1}e_{3}\,, \,
 e_{1}e_{4}\,, \, e_{2}e_{3}\,, \, e_{2}e_{4}\,, \, e_{3}e_{4} \,\}  \]

Now consider the isomorphism $\tau : C ^{2} ( \RR ^{4} )
\longrightarrow
\lso (4) = \lso ( \RR ^{4})$.
and the corresponding $\tg : \spin (4) \longrightarrow \so (4)$.
We have
\begin{align}
\tau ( e_{i} e_{j} ) \cdot v \, = \, &[ \, e_{i}e_{j} \, , \, \sum_{k=1}^{4} \,
v^{k} \, e_{k}   \, ] \notag \\
=\,&\sum_{k=1}^{4} \, v^{k} \, [ \, e_{i}e_{j} \, , \, e_{k}   \, ]  \notag \\
=\,&\sum_{k=1}^{4} \, v^{k} \, ( \, e_{i}e_{j}  e_{k}  \, - \, e_{k}e_{i}e_{j}
 \, ) \notag \\
=\,&\sum_{k=1}^{4} \, v^{k} \, ( \, e_{i}e_{j}  e_{k}  \, + \, e_{i}e_{k}e_{j}
 \, - \, e_{i}e_{k}e_{j}\,- \, e_{k}e_{i}e_{j} \, )  \notag \\
=\,&\sum_{k=1}^{4} \, v^{k} \, ( \, -2 \, Q(\, e_{j}\, , \, e_{k} \, )
\,e_{i}\,
 + \, 2 \, Q(\, e_{k} \, , \, e_{i} \,) \, e_{j}\, )  \notag \\
=\,&-2 \, v^{j} \, e_{i} \, + 2 \, v^{i} \, e_{j} \notag
\end{align}
Hence $\tau \, ( \, e_{i} e_{j} \, )$ corresponds to $2 \cdot m(i,j) \, \in
\lso
(4)$
where $m(i,j)\, = \, (m(i,j)_{\alpha \beta})$ is a matrix with entries
\[ m(i,j)_{\alpha \beta} \, = \,
	      \begin{cases}
		  1 , &\text{if $\alpha =j$ and $\beta=i$,} \\
		   -1, &\text{if $\alpha =i$ and $\beta=j$,} \\
		   0 , &\text{otherwise.}
		  \end{cases}  \]
Explicitly, we have
\[
\tau ( \, e_{1}e_{2} \, ) \, = \,2 \,
\begin{pmatrix}
0 & -1 & 0 &  0 \\
1 & 0 & 0 & 0 \\
0 & 0 & 0 & 0 \\
0 & 0 & 0 & 0
\end{pmatrix}
\qquad\qquad
\tau ( \, e_{1}e_{3} \, ) \, = \,2 \,
\begin{pmatrix}
0 & 0 & -1 &  0 \\
0 & 0 & 0 & 0 \\
1 & 0 & 0 & 0 \\
0 & 0 & 0 & 0
\end{pmatrix}  \]
\[
\tau ( \, e_{1}e_{4} \, ) \, = \,2 \,
\begin{pmatrix}
0 & 0 & 0 &  -1 \\
0 & 0 & 0 & 0 \\
0 & 0 & 0 & 0 \\
1 & 0 & 0 & 0
\end{pmatrix}
\qquad\qquad
\tau ( \, e_{2}e_{3} \, ) \, = \,2 \,
\begin{pmatrix}
0 & 0 & 0 &  0 \\
0 & 0 & -1 & 0 \\
0 & 1 & 0 & 0 \\
0 & 0 & 0 & 0
\end{pmatrix}  \]
\[
\tau ( \, e_{2}e_{4} \, ) \, = \,2 \,
 \begin{pmatrix}
0 & 0 & 0 &  0 \\
0 & 0 & 0 & -1 \\
0 & 0 & 0 & 0 \\
0 & 1 & 0 & 0
\end{pmatrix}
\qquad\qquad
\tau ( \, e_{3}e_{4} \, ) \, = \,2 \,
\begin{pmatrix}
0 & 0 & 0 &  0 \\
0 & 0 & 0 & 0 \\
0 & 0 & 0 & -1 \\
0 & 0 & 1 & 0
\end{pmatrix}\]
Notice that for $t \in [0 , 2 \pi )$, $\tau ( \, \frac{t}{2} \, e_{i}e_{j}\, )$
corresponds to the matrix $ t \cdot m(i,j)$
and
\[ (t \cdot m(i,j)_{\alpha \beta})^{2} \, = \, -t^{2} \cdot \Delta(i,j) \]
\[ (t \cdot m(i,j)_{\alpha \beta})^{4} \, = \, t^{4} \cdot \Delta(i,j) \]
where $\Delta(i,j) $ is a matrix with
\[ \Delta(i,j)_{\alpha \beta} \, = \,
\begin{cases}
1, &\text{if $\alpha = \beta = i$ or $\alpha=\beta =j$,} \\
0, &\text{otherwise.}
\end{cases}  \]
therefore we have
\[ \exp \, \tau \, ( \, \frac{t}{2} \, e_{i}e_{j} \, ) \, = \,  (\cos t -1)
\cdot
 \Delta(i,j) \, + \,1_{_{\so (4)}} \, + \, ( \sin t) \, m(i,j)  \]

Since we have the commutative relation
\[ \tg \, \cdot \, \exp \, = \, \exp \, \cdot \, \tau  ,\]
therefore
\begin{align}
\tg \, \exp \, ( \, t \, e_{i}e_{j} \, ) \, = \, &\exp \, \tau \,
 ( \, t e_{i} e_{j} \, ) \notag \\
=\,&(\cos 2t -1) \cdot \Delta(i,j) \, + \,1_{_{\so (4)}} \, +
\, ( \sin 2t) \, m(i,j)  \notag
\end{align}

Explicitly, we have
\[
\tg ( \, \exp ( \, t \, e_{1}e_{2} \, ) \, ) \, = \,
\begin{pmatrix}
\cos 2t & -\sin 2t & 0 &  0 \\
\sin 2t & \cos 2t & 0 & 0 \\
0 & 0 & 1 & 0 \\
0 & 0 & 0 & 1
\end{pmatrix}  \]
\[
\tg ( \, \exp ( \, t \, e_{1}e_{3} \, ) \, ) \, = \,
\begin{pmatrix}
\cos 2t & 0 & -\sin 2t &  0 \\
0 & 1 & 0 & 0 \\
\sin 2t & 0 & \cos 2t & 0 \\
0 & 0 & 0 & 1
\end{pmatrix}  \]
\[
\tg ( \, \exp ( \, t \, e_{1}e_{4} \, ) \, ) \, = \,
\begin{pmatrix}
\cos 2t & 0 & 0 &  -\sin 2t \\
0 & 1 & 0 & 0 \\
0 & 0 & 1 & 0 \\
\sin 2t & 0 & 0 & \cos 2t
\end{pmatrix}   \]
\[
\tg ( \, \exp ( \, t \, e_{2}e_{3} \, ) \, ) \, = \,
 \begin{pmatrix}
1 & 0 & 0 &  0 \\
0 & \cos 2t & -\sin 2t & 0 \\
0 & \sin 2t & \cos 2t & 0 \\
0 & 0 & 0 & 1
\end{pmatrix} \]
\[
\tg ( \, \exp ( \, t \, e_{2}e_{4} \, ) \, ) \, = \,
\begin{pmatrix}
1 & 0 & 0 &  0 \\
0 & \cos 2t & 0 & -\sin 2t \\
0 & 0 & 1 & 0 \\
0 & \sin 2t & 0 & \cos 2t
\end{pmatrix} \]
\[
\tg ( \, \exp ( \, t \, e_{3}e_{4} \, ) \, ) \, = \,
 \begin{pmatrix}
1 & 0 & 0 &  0 \\
0 & 1 & 0 & 0 \\
0 & 0 & \cos 2t & -\sin 2t \\
0 & 0 & \sin 2t & \cos 2t
\end{pmatrix} \]

\newpage
\bigskip \bigskip
\section{\bf The Quaternions and Related Representations}
In this section, we shall look at Clifford algebra and the Spin groups
 more concretely.
To do this, we need to use the quaternion. We shall define the Spin
 group in another more concrete
way and find more relations between $\spin (4)$ another classical
 groups. Our references here are
\cite{BD88}, \cite{LM89}, \cite{rH89}, and \cite{FH91}.

\bigskip
\subsection{Basic Linear Algebra of the Quaternions}
\begin{Def}
Let $\HH$ be the {\bf quaternion algebra} which, as a vector space
over $ \RR$, has four generators
denoted by $\, \1\, , \, \ii \, , \, \jj \, , \, \kk \,$ and its
multiplication satisfies the relations of the
group $Q_{8}$, namely
\[ \1 \, q \, = \, q \, \1 \, = \, q \qquad \forall q \in \HH   \]
\[ \ii ^{2} \, = \, \jj ^{2} \, = \, \kk^{2} \, = \, -1  \]
\[ \ii \, \jj \, = \, \kk , \qquad  \jj \, \kk \, = \, \ii , \qquad  \kk \,
\ii \, = \, \jj , \]
\[ \jj \, \ii \, = \, -\kk , \qquad  \kk \, \jj \, = \, -\ii , \qquad  \ii \,
\kk \, = \, -\jj , \]
\end{Def}
Since $\1$ is the unit of the algebra, it may be omitted for
convenience.
Therefore any element $q \in \HH$ may be uniquely
written as $a\,+\,
b\ii\,+\,c\jj\,+\,d\kk$
where $a,\, b, \, c, \, d \, \in \RR$.

Just like the complex field $\CC$, we have the following
\begin{Def}
For any $q\, = \, a\,+\,b\ii\,+\,c\jj\,+\,d\kk \in \HH$, the
 {\bf conjugate}
of $q$
is $\overline{q}\, = \,a\,-\,b\ii\,+\,c\jj\,+\,d\kk$. The
 {\bf real part} of
$q$ is $\Re (q) = a$ and the {\bf imaginary part} of $q$ is
$\Im (q) = b\ii\,+\,c\jj\,+\,d\kk$.

Moerover the {\bf real part} $\Re \HH$ of $\HH$ is $\RR$ while the
{\bf imaginary part} $\Im (\HH)$
of $\HH$ is $\, \RR\ii  \, \oplus \, \RR\jj \, \oplus \, \RR\kk $.

There is a natural {\bf inner product} or {\bf quadratic form} $Q ( \,
\cdot \, , \, \cdot \, )$
on $\HH$ given by
\[ Q( \, q_{1} \, , \, q_{2} \, ) \, = \, q_{1} \, \overline{q}_{2} \]
that means
\begin{align}
& Q ( \, a_{1} + b_{1} \ii + c_{1} \jj + d_{1} \kk \, , \, a_{2} + b_{2}
\ii + c_{2} \jj + d_{2}\kk \, ) \notag \\
=\,& a_{1} a_{2} \, + \, b_{1} b_{2}\, +\, c_{1} c_{2}\, +\, d_{1} d_{2}
 \notag
\end{align}

Therefore  $\HH $ becomes a normed vector space
with the corresponding
 {\bf norm} given by
\[   \Vert  \, q \, \Vert  \, = \, \sqrt{\, Q( \, q\, , \, q \, )\, } \, =
 \, \sqrt{\, q \, \overline{q}\, }  \]
Under this norm, the basis $\{ \, \1 \, , \, \ii \, , \, \jj \, , \, \kk \, \}$
 is {\bf orthonormal}.
\end{Def}

Obviously we also have the following
\begin{Lem}
For any $q\, = \, a\,+\,b\ii\,+\,c\jj\,+\,d\kk \in \HH$, if $q \neq 0$,
then there exists a
multiplicative inverse
\[ q^{-1} \,= \, \frac{ \overline{q} }{ \, \Vert q \Vert  ^{2} }   \]
such that
\[ q \, q^{-1} \, = \, q^{-1} \, q \, = \, \1   \]
That means $\HH$ is in fact a {\bf field}
\end{Lem}

\begin{Rem}
If we consider the $n$ -dimensional {\bf quaternion vector space}
$\HH ^{n}$, then we have a
standard inner product
\[ Q( \, (q_{1}, \cdots , q_{n}) \, , \, ( q'_{1}, \cdots , q'_{n} ) \, )\,
= \, \sum_{i=1}^{n}\, q_{i} \, \overline{ q'}_{i}  \]
and a standard norm
\[ \Vert \, (q_{1}, \cdots , q_{n}) \, \Vert \, =
\,  \sum_{i=1}^{n}\, q_{i} \, \overline{q}_{i} \]
\end{Rem}

\bigskip
\subsection{Real and Complex Representations}
Firstly, we look at the real representation.
\begin{Thm}
 The action of $\HH$ on itself by left
multiplication corresponds to the action of $\End (\RR ^{4})$
on $\RR ^{4}$. Explicitly,
we have a representation
\[  \HH \longrightarrow \End (\RR ^{4}) \]
given by
\[ a \, + \, b \ii \, + \, c \jj \, + \, d \kk \, \mapsto
\begin{pmatrix}
\, a\, & -b & -c & -d \\
b & a & -d & c \\
c & d & a & -b \\
d & -c & b & a
\end{pmatrix} \]
\end{Thm}

\begin{pf}
Firstly, we have
\[ \1 \, q \, = \, q \, = \, 1_{_{\End (\RR^{4})}} \, q , \qquad
 \forall q \in \HH \]
Now, look at the action of left multiplication by $\ii$ :
\[ \ii \, \1 \, = \, \ii, \qquad \ii \, \ii \, = \, -\1, \qquad \ii \,
\jj \, = \, \kk, \qquad \ii \, \kk \, = \, -\jj, \]
Therefore, by writing
\[ a \, + \, b \ii \, + \, c \jj \, + \, d \kk \quad \text{as} \quad
\begin{pmatrix}
a\\
b\\
c\\
d
\end{pmatrix}, \]
we have
\[ \ii \, \begin{pmatrix}
a\\
b\\
c\\
d
\end{pmatrix} \, = \,
\begin{pmatrix}
-b\\
a\\
-d\\
c
\end{pmatrix} \, = \,
\begin{pmatrix}
\, 0\, & -1 &\,0\, & 0 \\
1 & 0 & 0 & 0\\
0 & 0 & 0 & -1 \\
0 & 0 & 1 & 0
\end{pmatrix} \,
\begin{pmatrix}
a\\
b\\
c\\
d
\end{pmatrix} \]
Similarly
\[ \jj \, \begin{pmatrix}
a\\
b\\
c\\
d
\end{pmatrix} \, = \,
\begin{pmatrix}
-c\\
d\\
a\\
-b
\end{pmatrix} \, = \,
\begin{pmatrix}
\, 0\, & 0 &-1 & 0 \\
0 & 0 & 0 & \, 1 \,\\
1 & 0 & 0 & 0 \\
0 & -1 & 0 & 0
\end{pmatrix} \,
\begin{pmatrix}
a\\
b\\
c\\
d
\end{pmatrix} \]
\[ \kk \, \begin{pmatrix}
a\\
b\\
c\\
d
\end{pmatrix} \, = \,
\begin{pmatrix}
-d\\
-c\\
b\\
a
\end{pmatrix} \, = \,
\begin{pmatrix}
\, 0\, & \,0\, &0 & -1 \\
0 & 0 & -1 & 0\\
0 & 1 & 0 & 0 \\
1 & 0 & 0 & 0
\end{pmatrix} \,
\begin{pmatrix}
a\\
b\\
c\\
d
\end{pmatrix} \]
Hence we have the required representation.
It is a homomorphism because the matrices corresponding to $\,
 \1,\, \ii,\, \jj$ and $ \kk  $
satisfies the multiplicative relations of the group $Q_{8}$.
It is obviously injective and linear over $\RR$.
\end{pf}

Now we consider the complex representations of $\HH$.
\begin{Thm}
There is a isomorphism from $\HH$ into the algebra
\footnote{We will see that this algebra is in fact $\RR
\cdot \su (2)$.}
 of $2 \times 2$ complex matrices of the form
\[ \begin{pmatrix}
  c_{1} & c_{2} \\
  -\overline{c}_{2} & \overline{c}_{1}
\end{pmatrix}  \]
where $c_{1}, c_{2} \in \CC $.
\end{Thm}

\begin{pf}
The required map is given by
\[  a\,+\,b\ii\,+\,c\jj\,+\,d\kk \, \mapsto
\begin{pmatrix}
  a\,+\,b\ii & c+\,d\ii \\
  -c+\,d\ii  & a\,-\,b\ii
\end{pmatrix}   \]
It is obviously bijective.

We also have
\[ \1 \, \mapsto
\begin{pmatrix}
  \,1\, & 0 \\
  0  & \,1\,
\end{pmatrix}   \]

\[ \ii \, \mapsto
\begin{pmatrix}
 \,\ii\, & 0 \\
  0  & -\ii
\end{pmatrix}   \]

\[ \jj \, \mapsto
\begin{pmatrix}
  0 & \, 1 \, \\
  -1  & 0
\end{pmatrix}   \]

\[ \kk \, \mapsto
\begin{pmatrix}
  0 & \,\ii\, \\
  \,\ii \,  & 0
\end{pmatrix}   \]
Therefore it is straight forward to see that this map
 is an isomorphism.
\end{pf}

\begin{Thm}
There is an isomorphism
\[ \CC ^{2} \, \longrightarrow \, \HH \]
given by
\[ (\,c_{1}\, , \, c_{2} \, ) \, \mapsto \, c_{1} \, + \, c_{2} \jj  \]
This isomorphism is {\bf norm preserving}.
\end{Thm}

\begin{pf}
This map is obviously an isomorphism of  vector spaces over $\CC$.
Now the norm of $\CC ^{2}$ is given by
\[ \Vert (\, c_{1} \, , \, c_{2} \, ) \Vert \, = \, (\Re c_{1})^{2} \, + \,
 (\Im c_{1})^{2} \, + \,
(\Re c_{2})^{2} \, + \, (\Im c_{2})^{2}  \]
which is equal to the norm of its image in $\HH$.
Therefore $\CC ^{2}$ is isomorphic to $\HH$ as normed $\CC$
vector spaces.
\end{pf}

\bigskip
\subsection{The Quaternion Group}
If we generalize the concept of orthogonal and unitary
group to quaternion vector spaces,
we may have the following
\begin{Def}
Let the {\bf symplectic group}, $\sp (n)$, be the group of
 norm-preserving automorphisms of
the $\HH ^{n}$, ie.
\[ \sp (n) \, = \, \{ \, \phi \in \gl (n, \HH) \, | \, \Vert \,
\phi \, q \, \Vert \, = \, \Vert \,  q \, \Vert \, \forall q \in \HH \, \} \]
\end{Def}

As the most interesting example of symplectic groups, we have
\begin{Def}
Let the {\bf quaternion group} be the symplectic group
of $\HH$, which is simply
\[ \sp (1) \, = \, \{ \, q \in \HH \, | \, \Vert \, q \, \Vert \,
 = 1 \, \} \]
\end{Def}

\begin{Rem}
Every element $q$ in $\sp (1)$ has inverse
\( q ^{-1} \, = \, \overline{q}.  \)
\end{Rem}

\begin{Lem}
The underlying manifold of the Lie group $\sp (1)$ is the
 three-sphere $S^{3}$.
\end{Lem}
\begin{pf}
Obvious from definition.
\end{pf}

\begin{Lem} \label{L:sp1su2}
We have an isomorphism
\[ \sp (1) \, \cong \, \su (2)  \]
as Lie groups.
\end{Lem}

\begin{pf}
Since $\sp (1)$ consists of elements in $\HH$ of unit norm,
\[ \sp (1) \, = \, \{ \, a \, + \, b \ii \, + \, c \jj \, + \, d \kk \,
 | \, a^{2}+\,b^{2}+\,c^{2}+\,d^{2} =\, 1 \, \} \]
By the identification
\[
a \, + \, b \ii \, + \, c \jj \, + \, d \kk \, \mapsto
\begin{pmatrix}
a+b \ii & c+d \ii \\
-c + d \ii & a -b \ii
\end{pmatrix}
\]
we have
\[
 \det \, \begin{pmatrix}
a+b \ii & c+d \ii \\
-c + d \ii & a -b \ii
\end{pmatrix}  \,
= \, a^{2}+\,b^{2}+\,c^{2}+\,d^{2} =\, 1
\]
and
\begin{align}
\begin{pmatrix}
a+b \ii & c+d \ii \\
-c + d \ii & a -b \ii
\end{pmatrix} \,
\begin{pmatrix}
a+b \ii & c+d \ii \\
-c + d \ii & a -b \ii
\end{pmatrix} ^{\ast} \,
=\,&\begin{pmatrix}
a+b \ii & c+d \ii \\
-c + d \ii & a -b \ii
\end{pmatrix} \,
\begin{pmatrix}
a-b \ii & -c-d \ii \\
c-d \ii & a+b \ii
\end{pmatrix} \notag \\
=\,&\begin{pmatrix}
a^{2}+b^{2}+c^{2}+d^{2} & 0 \\
0 & a^{2}+b^{2}+c^{2}+d^{2}
\end{pmatrix} \notag \\
=\,&\begin{pmatrix}
1 & 0 \\
0 &  1
\end{pmatrix}  \notag
\end{align}
Therefore
\[
\begin{pmatrix}
a+b \ii & c+d \ii \\
-c + d \ii & a -b \ii
\end{pmatrix}\, \in \, \su (2)
\]
\end{pf}

The significance of $\sp (1)$ lies on the following action due to Hopf.
\begin{Def}
Let the {\bf Hopf action} $\hh$ of $\sp (1)$ on $\HH$ be the
 representation
\[ \hh \, : \, \sp (1) \, \longrightarrow \, \End (\HH ) \]
such that for any $\phi \in \sp (1) $ , $ \hh ( \phi )  :  \HH
\longrightarrow  \HH $ is given by the conjugation
\[ \hh ( \phi) \, : \, q \, \mapsto \, \phi \, q \, \phi ^{-1}
\qquad \forall q \in \HH \]
\end{Def}

\begin{Thm}
There is a $\RR \, \oplus \, \RR ^{3} $ representation of the
quaternion group $\sp (1)$ through
the Hopf action $\hh$. Explicitly, we have
\[ \hh \, : \, \sp (1) \, \longrightarrow \, \End ( \Re \HH ) \,
\oplus \, \End ( \Im \HH )  \]
where for any $\phi \in \sp (1)$,
\[ \hh ( \phi ) \, | _{_{\End ( \Re \HH )}} \, = \, 1 \]
and
\[  \hh ( \phi ) \, | _{_{\End ( \Im \HH )}} \, \in \so ( \Im
\HH ) \]
\end{Thm}

\begin{pf}
To find out the required representation, we pick
\[ \phi \, = \, \phi_{1}\,+\,\phi_{2}\ii\,+\,\phi_{3}\jj\,+
\,\phi_{4}\kk  \, \in \sp (1) \]
and
\[ q \, = \, q_{1} \, + \, q_{2} \ii \, + \, q_{3} \jj \, +
 \, q_{4} \kk \in \HH. \]
By using the real representation of the left multiplication, we have
\begin{align}
   \phi  q  \phi ^{-1}
=\,&
\begin{pmatrix}
  \phi_{1} & -\phi_{2} & -\phi_{3} & -\phi_{4} \\ \\
  \phi_{2} & \phi_{1} & -\phi_{4} & \phi_{3} \\ \\
  \phi_{3} & \phi_{4} & \phi_{1} & -\phi_{2} \\ \\
  \phi_{4} & -\phi_{3} & \phi_{2} & \phi_{1}
\end{pmatrix} \,
\begin{pmatrix}
  q_{1}\\ \\
  q_{2}\\ \\
  q_{3}\\ \\
  q_{4}
\end{pmatrix}
		     \, \phi ^{-1} \notag \\
=\,&
 \begin{pmatrix}
 \phi_{1}q_{1}  -\phi_{2}q_{2}  -\phi_{3}q_{3}  -\phi_{4}q_{4} \\ \\
 \phi_{2}q_{1}   +\phi_{1}q_{2}  -\phi_{4}q_{3}    +\phi_{3}q_{4} \\ \\
 \phi_{3}q_{1}  + \phi_{4}q_{2}   +\phi_{1}q_{3}  -\phi_{2}q_{4} \\ \\
 \phi_{4}q_{1}  -\phi_{3}q_{2}   +\phi_{2}q_{3}   + \phi_{1}q_{4}
 \end{pmatrix}
		       \, \phi ^{-1} \notag \\
=\,&
\begin{pmatrix}
\begin{matrix}
\phi_{1}q_{1}-\phi_{2}q_{2} \\ \,\,\,-\phi_{3}q_{3}-\phi_{4}q_{4} \,\,\,
\end{matrix}
& \begin{matrix}
 -(\phi_{2}q_{1}+\phi_{1}q_{2}\\-\phi_{4}q_{3}+\phi_{3}q_{4})
\end{matrix}
&\begin{matrix}
 -(\phi_{3}q_{1}+\phi_{4}q_{2}\\+\phi_{1}q_{3}-\phi_{2}q_{4})
\end{matrix}
&\begin{matrix}
  -(\phi_{4}q_{1}-\phi_{3}q_{2}\\+\phi_{2}q_{3}+\phi_{1}q_{4})
\end{matrix}
\\ \\
\begin{matrix}
\phi_{2}q_{1}+\phi_{1}q_{2}\\ \,\,\,-\phi_{4}q_{3}+\phi_{3}q_{4}  \,\,\,
\end{matrix}
&\begin{matrix}
 \phi_{1}q_{1}-\phi_{2}q_{2}\\-\phi_{3}q_{3}-\phi_{4}q_{4}
\end{matrix}
&\begin{matrix}
 -(\phi_{4}q_{1}-\phi_{3}q_{2}\\+\phi_{2}q_{3}+\phi_{1}q_{4})
\end{matrix}
&\begin{matrix}
 \phi_{3}q_{1}+\phi_{4}q_{2}\\+\phi_{1}q_{3}-\phi_{2}q_{4}
 \end{matrix}
\\ \\
\begin{matrix}
\phi_{3}q_{1}+\phi_{4}q_{2}\\ \,\,\,+\phi_{1}q_{3}-\phi_{2}q_{4}  \,\,\,
\end{matrix}
&\begin{matrix}
 \phi_{4}q_{1}-\phi_{3}q_{2}\\+\phi_{2}q_{3}+\phi_{1}q_{4}
\end{matrix}
&\begin{matrix}
 \phi_{1}q_{1}-\phi_{2}q_{2}\\-\phi_{3}q_{3}-\phi_{4}q_{4}
\end{matrix}
&\begin{matrix}
 -(\phi_{2}q_{1}+\phi_{1}q_{2}\\-\phi_{4}q_{3}+\phi_{3}q_{4})
\end{matrix}
\\ \\
\begin{matrix}
\phi_{4}q_{1}-\phi_{3}q_{2}\\ \,\,\,+\phi_{2}q_{3}+\phi_{1}q_{4}  \,\,\,
\end{matrix}
&\begin{matrix}
 -(\phi_{3}q_{1}+\phi_{4}q_{2}\\+\phi_{1}q_{3}-\phi_{2}q_{4})
\end{matrix}
&\begin{matrix}
 \phi_{2}q_{1}+\phi_{1}q_{2}\\-\phi_{4}q_{3}+\phi_{3}q_{4}
\end{matrix}
&\begin{matrix}
 \phi_{1}q_{1}-\phi_{2}q_{2}\\-\phi_{3}q_{3}-\phi_{4}q_{4}
\end{matrix}
\end{pmatrix} \,
\begin{pmatrix}
  \phi_{1}\\ \\ \\
  -\phi_{2}\\ \\ \\
  -\phi_{3}\\ \\ \\
  -\phi_{4}
\end{pmatrix} \notag \\ \notag\\
=\,& \, \, \centerdot \, \, \centerdot \, \, \centerdot \notag\\
 \notag\\
=\,&
\begin{pmatrix}
1 & 0 & 0 & 0 \\ \\
0 & \phi_{1}^{2}+\phi_{2}^{2}-\phi_{3}^{2}-\phi_{4}^{2} &
-2 \phi_{1}\phi_{4} + 2 \phi_{2}\phi_{3} & 2 \phi_{1}\phi_{3}
 + 2 \phi_{2}\phi_{4} \\ \\
0 & 2 \phi_{1}\phi_{4} + 2\phi_{2}\phi_{3} & \phi_{1}^{2}-
\phi_{2}^{2}+\phi_{3}^{2}-\phi_{4}^{2} & -2\phi_{1}\phi_{2}
+2\phi_{3}\phi_{4} \\ \\
0 & -2\phi_{1}\phi_{3}+2\phi_{2}\phi_{4} & 2\phi_{1}\phi_{2}
+2\phi_{3}\phi_{4} & \phi_{1}^{2}-\phi_{2}^{2}-\phi_{3}^{2}
+\phi_{4}^{2}
\end{pmatrix} \,
\begin{pmatrix}
q_{1}\\ \\
q_{2}\\ \\
q_{3}\\ \\
q_{4}
\end{pmatrix} \notag
\end{align}

The required representation is given by mapping
\[ \phi \, = \, \phi_{1}\,+\,\phi_{2}\ii\,+\,\phi_{3}\jj\,
+\,\phi_{4}\kk  \, \in \sp (1) \]
into the matrix
\[ \begin{pmatrix}
1 & 0 & 0 & 0 \\ \\
0 & \phi_{1}^{2}+\phi_{2}^{2}-\phi_{3}^{2}-
\phi_{4}^{2} & -2 \phi_{1}\phi_{4} +
 2 \phi_{2}\phi_{3} & 2 \phi_{1}\phi_{3} + 2 \phi_{2}\phi_{4} \\ \\
0 & 2 \phi_{1}\phi_{4} + 2\phi_{2}\phi_{3} & \phi_{1}^{2}
-\phi_{2}^{2}+\phi_{3}^{2}-\phi_{4}^{2} & -2\phi_{1}\phi_{2} +
2\phi_{3}\phi_{4} \\ \\
0 & -2\phi_{1}\phi_{3}+2\phi_{2}\phi_{4} & 2\phi_{1}\phi_{2}+
2\phi_{3}\phi_{4} & \phi_{1}^{2}-\phi_{2}^{2}-\phi_{3}^{2}+\phi_{4}^{2}
\end{pmatrix} \]
It is tedious but easy to show that
\[ \begin{pmatrix}
 \phi_{1}^{2}+\phi_{2}^{2}-\phi_{3}^{2}-\phi_{4}^{2} & -
2 \phi_{1}\phi_{4} + 2 \phi_{2}\phi_{3} & 2 \phi_{1}\phi_{3} +
 2 \phi_{2}\phi_{4} \\ \\
 2 \phi_{1}\phi_{4} + 2\phi_{2}\phi_{3} & \phi_{1}^{2}-
\phi_{2}^{2}+\phi_{3}^{2}-\phi_{4}^{2} & -2\phi_{1}\phi_{2}
+2\phi_{3}\phi_{4} \\ \\
 -2\phi_{1}\phi_{3}+2\phi_{2}\phi_{4} & 2\phi_{1}\phi_{2}
+2\phi_{3}\phi_{4} & \phi_{1}^{2}-\phi_{2}^{2}-\phi_{3}^{2}+\phi_{4}^{2}
\end{pmatrix} \quad \in \o (3) \]
Since $\sp (1)$ is connected and for $\1 \, \in \sp (1)$,
\[ \det ( \, \hh( \1) \, |_{_{\End ( \Im \HH ) }} \, ) \, = \, 1 \]
therefore $\hh |_{_{\End ( \Im \HH ) }} ( \, \sp (1) \,) \,
\subset \so ( \Im \HH )$.
\end{pf}

\begin{Cor} \label{C:sp1so3}
There is a two-fold covering
\[ 0 \, \longrightarrow \, \ZZ _{2} \, \longrightarrow \,
\sp (1) \, \overset{\rho }{\longrightarrow}
 \, \so (3) \, \longrightarrow \, 0 \]
where $ \rho \, = \,\hh |_{_{\End ( \Im \HH ) }}$.
\end{Cor}
\begin{pf}
Let \( \phi \, =\, \phi_{1}\,+\,\phi_{2}\ii\,+\,\phi_{3}\jj\,+
\,\phi_{4}\kk \) be any element in
$ \ker ( \rho ) $.
Therefore $ \rho ( \, \phi \, ) \, = \, 1 _{_{\End ( \Im \HH ) }}$
 which is diagonal.
However, for any
\[ \begin{pmatrix}
 \phi_{1}^{2}+\phi_{2}^{2}-\phi_{3}^{2}-\phi_{4}^{2} & -
2 \phi_{1}\phi_{4} + 2 \phi_{2}\phi_{3} & 2 \phi_{1}\phi_{3} +
 2 \phi_{2}\phi_{4} \\ \\
 2 \phi_{1}\phi_{4} + 2\phi_{2}\phi_{3} & \phi_{1}^{2}-
\phi_{2}^{2}+\phi_{3}^{2}-\phi_{4}^{2} & -2\phi_{1}\phi_{2} +
2\phi_{3}\phi_{4} \\ \\
 -2\phi_{1}\phi_{3}+2\phi_{2}\phi_{4} & 2\phi_{1}\phi_{2}+
2\phi_{3}\phi_{4} & \phi_{1}^{2}-\phi_{2}^{2}-\phi_{3}^{2}+\phi_{4}^{2}
\end{pmatrix} \]
to be diagonal, we must have
\[ \phi_{i}\phi_{j} \, = \, 0 , \qquad \text{for } 1 \leq i < j \leq 4 \]
which implies one of the following cases
\begin{enumerate}
  \item $\phi_{1} \, \neq \, 0 , \quad \text{and} \quad \phi_{2} \, =
 \, \phi_{3} \,= \, \phi_{4} \,=\, 0  $
  \item $\phi_{2} \, \neq \, 0 , \quad \text{and} \quad \phi_{3} \, =
 \, \phi_{4} \,= \, \phi_{1} \,=\, 0  $
  \item $\phi_{3} \, \neq \, 0 , \quad \text{and} \quad \phi_{4} \, =
\, \phi_{1} \,= \, \phi_{2} \,=\, 0  $
  \item $\phi_{4} \, \neq \, 0 , \quad \text{and} \quad \phi_{1} \, =
 \, \phi_{2} \,= \, \phi_{3} \,=\, 0  $
\end{enumerate}
These cases corresponds to the following elements in $\sp (1) $:
\begin{enumerate}
  \item $\quad \Longrightarrow \quad \phi_{1} \, =
 \, \pm 1 \quad \Longrightarrow \quad \phi \, = \, \pm \1 $
  \item $\quad \Longrightarrow \quad \phi_{2} \, =
 \, \pm 1 \quad \Longrightarrow \quad \phi \, = \, \pm \ii $
  \item $\quad \Longrightarrow \quad \phi_{3} \, =
 \, \pm 1 \quad \Longrightarrow \quad \phi \, = \, \pm \jj $
  \item $\quad \Longrightarrow \quad \phi_{4} \, =
 \, \pm 1 \quad \Longrightarrow \quad \phi \, = \, \pm \kk $
\end{enumerate}
and we have the following correspondence :
\[
\rho (\, \pm \1 \,)\, = \,
\begin{pmatrix}
1 & 0 & 0  \\
0 & 1 & 0  \\
0 & 0 & 1
\end{pmatrix}
\qquad\qquad
\rho (\, \pm \ii \,)\, = \,
\begin{pmatrix}
1 & 0 & 0  \\
0 & -1 & 0  \\
0 & 0 & -1
\end{pmatrix}  \]
\[
\rho (\, \pm \jj \,)\, = \,
\begin{pmatrix}
-1 & 0 & 0  \\
0 & 1 & 0  \\
0 & 0 & -1
\end{pmatrix}
\qquad\qquad
\rho (\, \pm \kk \,)\, = \,
\begin{pmatrix}
-1 & 0 & 0  \\
0 & -1 & 0  \\
0 & 0 & 1
\end{pmatrix}  \]
Therefore we see that $\ker ( \, \rho \, ) \, = \, \ZZ _{2} $.
In fact, for any $\phi $  in $\sp (1)$ ,
\[ \rho ( \, - \phi \, ) \, = \, \rho ( \,  \phi \, ) \]
\end{pf}

\begin{Cor}
There is a two-fold covering
\[ 0 \, \longrightarrow \, \ZZ _{2} \, \longrightarrow \,
\su (2) \, \longrightarrow
 \, \so (3) \, \longrightarrow \, 0 \]
\end{Cor}

\begin{pf}
By {\bf Lemma \ref{L:sp1su2}} , we have an isomorphism
\[ \su (2) \, \longrightarrow \, \sp (1) \]
\[ \begin{pmatrix}
a + b \ii & c + d \ii \\
- c + d \ii & a - b \ii
\end{pmatrix}
\, \mapsto \,
a \, + \, b \ii \, + \, c \jj \, + \, d \kk \]
which combines with the two-fold cover
\[ \sp (1) \, \longrightarrow \, \so (3) \]
of {\bf  Corollary \ref{C:sp1so3}} gives the required two-fold cover :
\[
\begin{pmatrix}
a + b \ii & c + d \ii \\
- c + d \ii & a - b \ii
\end{pmatrix}
\, \mapsto \,
\begin{pmatrix}
 a^{2}+b^{2}-c^{2}-d^{2} & -2 ad + 2 bc & 2 ac + 2 bd \\ \\
 2 ad + 2bc & a^{2}-b^{2}+c^{2}-d^{2} & -2ab +2cd \\ \\
 -2ac+2bd & 2ab+2cd & a^{2}-b^{2}-c^{2}+d^{2}
\end{pmatrix} \]
\end{pf}

\bigskip
\subsection{Realization of Low Dimensional Clifford Algebras}
We now apply the quaternion to the representation of Clifford
 algebras in low dimensions.

\begin{Lem}
We have an isomorphism
\[ C(\,\RR\,)\,\cong \CC \]
of algebras over $\RR$.
 \end{Lem}

\begin{pf}
By definition,
\[ C(\,\RR\,)\, = \, \{ \, e_{1} \, | \, e_{1}e_{1} \, = \, -1 \}   \]
Hence $C( \, \RR \, ) $ contains elements $ a + b e_{1} $
with $e_{1} ^{2} = -1 $.
So by identifying $e_{1}$ with $\ii$, we may get the
 required isomorphism.
\end{pf}

\begin{Lem}
We have an isomorphism
\[ C(\,\RR ^{2}\,)\,\cong \, \HH \]
of algebras over $\RR$.
 \end{Lem}

\begin{pf}
By choosing orthonormal basis $e_{1}, \, e_{2}$ of  $\RR ^{2}$,
\[ C(\,\RR ^{2}\,)\, = \, \{ \, e_{1} , \, e_{2}\, | \, e_{i}e_{j} \,
+ \, e_{j}e_{i} \, = \, -2 Q( \, e_{i} \, , \, e_{j} \, ) \}   \]
Hence $C( \, \RR ^{2} \, ) $ contains elements $ a + b e_{1}
+ c e_{2} + d e_{1}e_{2} $ with
\[ e_{1} ^{2} \, = \, e_{2} ^{2} \, = \, ( e_{1}e_{2} )^{2} \,
= \, -1 \]
\[ e_{1}e_{2}= (e_{1}e_{2}) , \quad e_{2} ( e_{1}e_{2})
= e_{1} , \quad ( e_{1}e_{2}) e_{1} = e_{2} \]
\[ e_{2} e_{1} = - ( e_{1}e_{2}) , \quad ( e_{1}e_{2}) e_{2}
= - e_{1} , \quad  e_{1} ( e_{1}e_{2}) = - e_{2} \]
So by identifying $e_{1}$ with $\ii$, $e_{2}$
with $\jj$ and $e_{1}e_{2}$ with $\kk$,
we may get the required isomorphism.
\end{pf}

\begin{Lem} \label{L:c3hh}
We have an isomorphism
\[ C(\,\RR ^{3}\,)\,\cong \, \HH \, \oplus \, \HH  \]
of algebras over $\RR$.
 \end{Lem}

\begin{pf}
By choosing orthonormal basis $e_{1}, \, e_{2}, \, e_{3}$
of  $\RR ^{3}$,
\[ C(\,\RR ^{3}\,)\, = \, \{ \, e_{1} , \, e_{2}, \, e_{3}\, | \,
 e_{i}e_{j} \, + \, e_{j}e_{i} \, = \, -2 Q( \, e_{i} \, , \, e_{j} \, ) \}
\]
As a real vector space, $C(\,\RR ^{3}\,)$ has eight generators
\[ 1, \quad e_{1}, \quad e_{2},\quad e_{3},\quad
 e_{1}e_{2} ,\quad e_{1}e_{3} ,\quad e_{2}e_{3} ,
\quad e_{1}e_{2}e_{3} \]
On the other hand, $\HH \, \oplus \, \HH$ also
has eight generators
\[ \1 \oplus 0 , \quad \ii \oplus 0 , \quad \jj
\oplus 0 , \quad \kk \oplus 0 ,
 \qquad 0\oplus \1, \quad 0\oplus \ii, \quad 0\oplus
\jj, \quad 0\oplus \kk \]
The required map
\[  C(\,\RR ^{3}\,)\, \longrightarrow \, \HH \,
 \oplus \, \HH \]
is given by
\begin{alignat}{4}
  1 \,       & \mapsto & \, \1 \, &\oplus \, \1,
    & \qquad e_{1}e_{2}e_{3} \, & \mapsto &\, -\1 \, &\oplus \, \1
 \notag \\
  e_{1} \, & \mapsto &\, \kk \, &\oplus \, \kk,    &
\qquad e_{1}e_{2} \,        & \mapsto &\, -\ii \, &\oplus \,  \ii  \notag \\
  e_{2} \, & \mapsto &\, \jj \, &\oplus \, -\jj, &
\qquad e_{1}e_{3} \,       & \mapsto &\, -\jj \, &\oplus \,  -\jj  \notag \\
  e_{3} \, & \mapsto &\, -\ii \, & \oplus \, -\ii,     &
\qquad e_{2}e_{3} \,        & \mapsto &\, \kk \, &\oplus \,  -\kk  \notag
\end{alignat}
while the inverse
\[ \HH \, \oplus \, \HH \, \longrightarrow \, C(\,\RR ^{3}\,) \]
is given by
\begin{alignat}{4}
  \1 \,&\oplus\, 0 \,&  & \mapsto \, \frac{1- e_{1}e_{2}e_{3}}{2} ,
 & \qquad 0 &\oplus\, \1 \, & & \mapsto \,
\frac{1+ e_{1}e_{2}e_{3}}{2} \notag \\
  \ii \,&\oplus\, 0 \, & & \mapsto \, \frac{-e_{1}e_{2}-e_{3}}{2} ,
  & \qquad 0 &\oplus \,\ii \, & & \mapsto \,
\frac{e_{1}e_{2}-e_{3}}{2} \notag \\
  \jj \,&\oplus\, 0 \, & & \mapsto \, \frac{-e_{1}e_{3}+e_{2}}{2} ,
  & \qquad 0&\oplus \,\jj \,  & & \mapsto \,
\frac{-e_{1}e_{3}-e_{2}}{2} \notag \\
  \kk \,&\oplus \,0 \,& & \mapsto \, \frac{e_{2}e_{3}+e_{1}}{2} ,
 & \qquad 0 &\oplus \,\kk \, & & \mapsto \,
\frac{-e_{2}e_{3}+e_{1}}{2}\notag
\end{alignat}
\end{pf}
\begin{Rem}
There are other homomorphisms between $C(\RR ^{3})$
 and $\HH \oplus \HH$. The one we choose
in the above proof is the one we need in the future
which is compatible with our description of the
action of $C^{2}(\RR ^{4}) \otimes \CC$ on the spinors.
\end{Rem}

To find the representation of higher dimensional Clofford
 algebra is more difficult. However, our
main interest is the even Clifford algebra of dimension four.
 This can be found with the help of the
following observation:

\begin{Lem} \label{L:c+n+1cn}
We have a canonical isomorphism of the two  algebras
\[ C^{+}(\, \RR ^{n+1}\,  ) \, \cong \, C( \, \RR ^{n} \, )  \]
\end{Lem}
\begin{pf}
Pick an orthonormal basis $\{ e_{i} \}_{_{i = 1 , \dots , n}}$ of
 $\RR ^{n}$ and embed $\RR ^{n}$
into $\RR ^{n+1}$ canonically. Pick $e_{n+1} \in \RR ^{n+1}$
 such that $\{ e_{i} \}_{_{i = 1 , \dots , n+1}}$
is an orthonormal basis of $\RR ^{n+1}$. The embedding
$ \RR ^{n} \hookrightarrow \RR ^{n+1}$
induces the map
\[  C( \, \RR ^{n} \, ) \, \hookrightarrow \,
C^{+}(\, \RR ^{n+1}\,  ) \]
\[   e_{i_{1}} e _{i_{2}} \dots e_{i_{2k+1}} \, \mapsto
\,  e_{i_{1}} e _{i_{2}} \dots e_{i_{2k+1}} e_{n+1}  \]
\[  e_{i_{1}} e _{i_{2}} \dots e_{i_{2k}} \, \mapsto
\,  e_{i_{1}} e _{i_{2}} \dots e_{i_{2k }}  \]
for any $k = 1 , \dots , [ \frac{n}{2} ] $.
This map is obviously an isomorphism.
\end{pf}

\begin{Cor} \label{C:c+4hh}
We have an isomorphism
\[ C^{+}(\, \RR ^{4} \, ) \, \cong \, \HH \, \oplus \, \HH  \]
\end{Cor}

\begin{pf}
By {\bf Lemma \ref{L:c+n+1cn}}, we have
\[ C^{+}(\, \RR ^{4} \, ) \, \cong \, C(\, \RR ^{3} \, )  \]
while by {\bf Lemma \ref{L:c3hh}}, we have
\[ C(\, \RR ^{3} \, ) \, \cong \, \HH \, \oplus \, \HH  \]
Their composition gives the required isomorphism.

Explicitly, we have
\[  C^{+}(\,\RR ^{4}\,)\, \longrightarrow \, \HH \, \oplus \, \HH \]
by
\begin{alignat}{4}
  1 \,       & \mapsto &\, \1 \, &\oplus \, \1,      &
 \qquad e_{1}e_{2}e_{3}e_{4} \, & \mapsto &\, -\1 \,
&\oplus \, \1  \notag \\
  e_{1}e_{4} \, & \mapsto &\, \kk \, &\oplus \, \kk,
& \qquad e_{1}e_{2} \,        & \mapsto &\, -\ii \,
&\oplus \,  \ii  \notag \\
  e_{2}e_{4} \, & \mapsto &\, \jj \, &\oplus \, -\jj,
& \qquad e_{1}e_{3} \,       & \mapsto &\, -\jj \,
&\oplus \,  -\jj  \notag \\
  e_{3}e_{4} \, & \mapsto &\, -\ii \, &\oplus \, -\ii,
  & \qquad e_{2}e_{3} \,        & \mapsto &\, \kk \, &\oplus \,  -\kk
 \notag
\end{alignat}
and the inverse
\[ \HH \, \oplus \, \HH \, \longrightarrow \, C^{+}(\,\RR ^{4}\,) \]
by
\begin{alignat}{4}
  \1 \,&\oplus\, 0 \, & & \mapsto \, \frac{1-e_{1}e_{2}e_{3}e_{4}}{2} ,
 & \qquad 0 &\oplus\, \1 \, & & \mapsto \,
\frac{1+e_{1}e_{2}e_{3}e_{4}}{2} \notag \\
  \ii \,&\oplus\, 0 \, & & \mapsto \,
\frac{-e_{1}e_{2}-e_{3}e_{4}}{2} ,
  & \qquad 0 &\oplus \,\ii \, & & \mapsto \,
\frac{e_{1}e_{2}-e_{3}e_{4}}{2} \notag \\
  \jj \,&\oplus\, 0 \, & & \mapsto \,
 \frac{-e_{1}e_{3}+e_{2}e_{4}}{2} ,
  & \qquad 0&\oplus \,\jj \,  & & \mapsto \,
\frac{-e_{1}e_{3}-e_{2}e_{4}}{2} \notag \\
  \kk \,&\oplus \,0 \,& & \mapsto \,
 \frac{e_{2}e_{3}+e_{1}e_{4}}{2} ,
& \qquad 0& \oplus \,\kk \,& & \mapsto \,
\frac{-e_{2}e_{3}+e_{1}e_{4}}{2}\notag
\end{alignat}
\end{pf}

\bigskip
\subsection{Alternative Definition of the Spin Group}
\begin{Def}
Let $C ^{\times}(V)$ be the {\bf multiplicative group}
of the Clifford algebra $C(V)$ over $V$.
Let the {\bf pin group} $\pin (V)$ of a real vector space
 $V$ be the subgroup of the multiplicative
group generated by unit vectors in the unit sphere in $V$, ie.
\[ \pin (V) \, = \, \{ \, v_{1} v_{2} \cdots v_{k} \in C ^{\times}
(V) \, | \, v_{i} \, \in S^{\dim V \, -1} \, \hookrightarrow V \} \]
\end{Def}

\begin{Rem}
The $\ZZ _{2}$-grading of the Clifford algebra
\[ C(V) \, = \, C^{+} (V) \, \oplus \, C^{-} (V) \]
induces a similar $\ZZ _{2}$-grading of the pin group
\[ \pin (V) \, = \, \pin ^{+}(V) \, \cup \, \pin ^{-}(V) \]
where $\pin ^{\pm}(V) \, = \, \pin (V) \, \cap \, C^{\pm} (V) $.
\end{Rem}

\begin{Def}
The {\bf spin group} $\spin (V)$ is the subgroup  $\pin ^{+}(V)$ of
$\pin (V)$. In other words,
\[ \spin (V) \, = \, \{ \, v_{1} v_{2} \cdots v_{2k} \, | \, v_{i} \, \in
S^{\dim V \, -1} \, \hookrightarrow V \} \]
\end{Def}

\begin{Rem}
This definition is equivalent to the previous one. In fact,
\[ \cos t \, + \, \sin t \, e_{i}e_{j} \, = \, e_{i} \, ( \, - \cos t \, e_{i}
\, +
\, \sin t \, e_{j} \, ) \]
and for any $v_{1} , v_{2} \,  \in  S^{\dim V \, -1}$, $v_{1} v_{2}$ may
 be expressed as the
product of elements of the form $\cos t \, + \, \sin t \, e_{i}e_{j}$.
\end{Rem}

\begin{Thm} \label{T:spin4sp1su2}
We have the following isomorphism
\[ \spin (4) \, \cong \, \sp (1) \,\times \,\sp (1) \, \cong \, \su (2)
\,\times \,\su (2) \]
of Lie groups.
\end{Thm}

\begin{pf}
By definition,
\[ \spin (4) \, = \, \pin (4) \, \cap \, C^{+} (\RR ^{4}) \]
Now $\pin (4) $ is the multiplicative group generated by
elements in $\sp (1)$
and by {\bf Corollary \ref{C:c+4hh}}, we have
\[ \spin (4) \, = \, \sp (1) \, \times \, \sp (1)  \]
Explicitly, by our new definition, $\spin (4)$ is generated by
$v_{1} v_{2}$ where
\[ v_{i} \, = \, a_{i} e_{1} \, + \, b_{i} e_{2} \, + \, c_{i} e_{3} \,
+ \, d_{i} e_{4} \, \in S ^{3} \hookrightarrow \RR ^{4} \]
Hence
\begin{align}
v_{1}v_{2} \, = \, &- \,( \, a_{1}a_{2} \, + \, b_{1}b_{2} \,
+ \, c_{1}c_{2} \, + \, d_{1}d_{2} \, ) \notag \\
&+ \, ( \, a_{1}b_{2} \, - \, a_{2}b_{1} \, ) \, e_{1}e_{2} \,
+ \, ( \, c_{1}d_{2} \, - \, c_{2}d_{1} \, ) \, e_{3}e_{4}  \notag \\
&+ \, ( \, a_{1}d_{2} \, - \, a_{2}d_{1} \, ) \, e_{1}e_{4} \,
+ \, ( \, b_{1}c_{2} \, - \, b_{2}c_{1} \, ) \, e_{2}e_{3}  \notag \\
&+ \, ( \, a_{1}c_{2} \, - \, a_{2}c_{1} \, ) \, e_{1}e_{3} \,
+ \, ( \, b_{1}d_{2} \, - \, b_{2}d_{1} \, ) \, e_{2}e_{4}  \notag
\end{align}
which corresponds to the following element in $\sp (1) \times \sp (1) $:
\begin{align}
 &- \,( \, a_{1}a_{2} \, + \, b_{1}b_{2} \, + \, c_{1}c_{2} \,
+ \, d_{1}d_{2} \, ) \, (\1 \oplus \1) \notag \\
&+ \, ( \, a_{1}b_{2} \, - \, a_{2}b_{1} \, ) \,( -\ii \oplus \ii )\,
 \, ( \, c_{1}d_{2} \, - \, c_{2}d_{1} \, ) \,(-\ii \oplus -\ii) \notag \\
&+ \, ( \, a_{1}d_{2} \, - \, a_{2}d_{1} \, ) \, (\kk \oplus \kk) \,
+ \, ( \, b_{1}c_{2} \, - \, b_{2}c_{1} \, ) \, (\kk \oplus -\kk)  \notag \\
&+ \, ( \, a_{1}c_{2} \, - \, a_{2}c_{1} \, ) \,( - \jj \oplus - \jj) \,
+ \, ( \, b_{1}d_{2} \, - \, b_{2}d_{1} \, ) \, (\jj \oplus -\jj)  \notag \\
=\, & \bigl( \,- \,( \, a_{1}a_{2} \, + \, b_{1}b_{2} \, + \, c_{1}c_{2} \,
+ \, d_{1}d_{2} \, ) \, \1 \,
+ \, ( \,-\, a_{1}b_{2} \, + \, a_{2}b_{1} \,  -\, c_{1}d_{2} \,
+ \, c_{2}d_{1} \, ) \, \ii \notag \\
&+ \, (\, -\, a_{1}c_{2} \, + \, a_{2}c_{1}  \, + \, b_{1}d_{2} \,
- \, b_{2}d_{1} \, ) \, \jj \,
+ \, (\, a_{1}d_{2} \, - \, a_{2}d_{1} \,  +  \, b_{1}c_{2} \, -
 \, b_{2}c_{1} \, ) \, \kk  \,\bigr)\notag \\
 &  \oplus\,\bigl( \, - \,( \, a_{1}a_{2} \, + \, b_{1}b_{2} \, +
 \, c_{1}c_{2} \, + \, d_{1}d_{2} \, ) \,  \1 \,
+ \, ( \, a_{1}b_{2} \, - \, a_{2}b_{1} \, -\, c_{1}d_{2} \,
+ \, c_{2}d_{1} \, )\, \ii \notag \\
&+ \, ( \, -\, a_{1}c_{2} \, + \, a_{2}c_{1} \, -\, b_{1}d_{2} \,
+ \, b_{2}d_{1} \, ) \, \jj \,
+ \, ( \, a_{1}d_{2} \, - \, a_{2}d_{1}  \, -  \, b_{1}c_{2} \, +
 \, b_{2}c_{1} \, ) \,  \kk  \, \bigr) \notag
\end{align}

If we use the old definition, then $\spin (4)$ is generated by
 elements of the form
$ \cos t \, + \, \sin t \, e_{i}e_{j} $.
Therefore the map
\[ \spin (4) \, \cong \, \sp (1) \, \times \, \sp (1) \]
 is given by
\begin{alignat}{2}
\cos t \, + \, \sin t \, e_{1}e_{2} \, & \mapsto \, & (\,\cos t \, -
 \, \sin t \, \ii \,)\, &\oplus \, (\,\cos t \, + \, \sin t \, \ii\,)  \notag
\\
\cos t \, + \, \sin t \, e_{3}e_{4} \, & \mapsto \, & (\,\cos t \, -
 \, \sin t \, \ii \,)\, &\oplus \, (\,\cos t \, - \, \sin t \, \ii \,) \notag
\\
\cos t \, + \, \sin t \, e_{1}e_{3} \, & \mapsto \, & (\,\cos t \,
- \, \sin t \, \jj \,)\, &\oplus \,(\, \cos t \, - \, \sin t \, \jj \,) \notag
\\
\cos t \, + \, \sin t \, e_{2}e_{4} \, & \mapsto \, & (\,\cos t \,
+ \, \sin t \, \jj \,)\, &\oplus \, (\,\cos t \, - \, \sin t \, \jj \,) \notag
\\
\cos t \, + \, \sin t \, e_{1}e_{4} \, & \mapsto \, &( \,\cos t \,
+ \, \sin t \, \kk \,)\, &\oplus \,(\, \cos t \, + \, \sin t \, \kk \,)
\notag \\
\cos t \, + \, \sin t \, e_{2}e_{3} \, & \mapsto \, & (\,\cos t \,
+ \, \sin t \, \kk \,)\, &\oplus \,(\, \cos t \, - \, \sin t \, \kk \,) \notag
\end{alignat}

By applying {\bf Lemma \ref{L:sp1su2}}, we have
\[ \spin (4) \,  \cong \, \su (2) \,\times \,\su (2) \]
given by
\begin{alignat}{2}
 \cos t \, + \, \sin t \, e_{1}e_{2} \, &\mapsto \, &
\begin{pmatrix}
\cos t - \, \sin t \, \ii  & 0 \\
0 & \cos t + \sin t \, \ii
\end{pmatrix}
\, &\oplus \,
\begin{pmatrix}
\cos t + \, \sin t \, \ii  & 0 \\
0 & \cos t - \sin t \, \ii
\end{pmatrix} \notag\\
 \cos t \, + \, \sin t \, e_{3}e_{4} \, &\mapsto \, &
\begin{pmatrix}
\cos t - \, \sin t \, \ii  & 0 \\
0 & \cos t + \sin t \, \ii
\end{pmatrix}
\, &\oplus \,
\begin{pmatrix}
\cos t - \, \sin t \, \ii  & 0 \\
0 & \cos t + \sin t \, \ii
\end{pmatrix}\notag\\
 \cos t \, + \, \sin t \, e_{1}e_{3} \, &\mapsto \, &
\begin{pmatrix}
\cos t   & -\sin t \\
\sin t & \cos t
\end{pmatrix}
\, &\oplus \,
\begin{pmatrix}
\cos t   & -\sin t \\
\sin t & \cos t
\end{pmatrix} \notag\\
 \cos t \, + \, \sin t \, e_{2}e_{4} \, &\mapsto \, &
\begin{pmatrix}
\cos t   & \sin t \\
-\sin t & \cos t
\end{pmatrix}
\, &\oplus \,
\begin{pmatrix}
\cos t   & -\sin t \\
\sin t & \cos t
\end{pmatrix} \notag\\
 \cos t \, + \, \sin t \, e_{1}e_{4} \, &\mapsto \, &
\begin{pmatrix}
\cos t   & \sin t \,\ii \\
\sin t \,\ii & \cos t
\end{pmatrix}
\, &\oplus \,
\begin{pmatrix}
\cos t   & \sin t \,\ii\\
\sin t \,\ii& \cos t
\end{pmatrix} \notag\\
 \cos t \, + \, \sin t \, e_{2}e_{3} \, &\mapsto \, &
\begin{pmatrix}
\cos t   & \sin t \,\ii \\
\sin t \,\ii & \cos t
\end{pmatrix}
\, &\oplus \,
\begin{pmatrix}
\cos t   & -\sin t \,\ii\\
-\sin t \,\ii& \cos t
\end{pmatrix} \notag
\end{alignat}
\end{pf}

\newpage
\bigskip \bigskip
\section{\bf Spinors}
We now proceed to define spinors and the spinor
 representation of Clifford algebras.
\bigskip
\subsection{Basic Properties}
\begin{Def}
Let $\{ e_{i} \}_{i=1, \dots , \dim V} $ be an
{\bf oriented}, orthonormal basis of $V$,
that means there is a preferred ordering of
 the basis elements modulo even permutations.
The {\bf chirality operator}
\footnote{Also known as the {\bf volume element}
or the {\bf complex unit}.
Physicists also denote it as ${\bf \gamma^{5}}$ for
 the four dimensional case.}
is
\[ \Gamma \, = \, \ii ^{p} \, e_{1} \dots e_{_{\dim V}}
\quad \in C(V) \otimes \CC \]
where
\[ p \, = \, \bigl[ \frac{\dim V +1}{2} \bigr] \, = \,
\begin{cases}
\frac{\dim V}{2} &\text{if $\dim V$ is even,} \\
\frac{\dim V+1}{2} &\text{if $\dim V$ is odd.}
\end{cases} \]
\end{Def}

\begin{Lem}
The Chirality operator satisfies
\[ \Gamma ^{2} \, = \, 1_{_{C(V)}}  \]
and it {\bf  super-anticommutes} with elements $v \in V$, ie.
\[ \Gamma \, v \, + \, (-1)^{\dim V} \, v \, \Gamma \, = \, 0  \]
\end{Lem}

\begin{pf}
By direct computation, we have
\begin{align}
\Gamma^{2}\, = \, &  \Gamma \, = \, \ii ^{p} \, e_{1}
 \dots e_{_{\dim V}} \,  \ii ^{p} \, e_{1} \dots e_{_{\dim V}} \notag \\
=\,& (-1) ^{p}\, e_{1} \dots e_{_{\dim V}}\, e_{1} \dots
 e_{_{\dim V}} \notag \\
=\,&(-1) ^{p}(-1)^{\dim V -1} \,  e_{1} e_{1} \, e_{2}
\dots e_{_{\dim V}}\, e_{2} \dots e_{_{\dim V}} \notag \\
=\,&(-1) ^{p}(-1)^{\dim V -1} (-1)^{\dim V -2}\,  e_{1} e_{1}
\, e_{2}e_{2}\, e_{3} \dots e_{_{\dim V}}\, e_{3} \dots e_{_{\dim V}}
 \notag \\
=\,& \centerdot \,  \centerdot \,  \centerdot \notag \\
=\,&(-1) ^{p}(-1)^{\dim V -1}  (-1)^{\dim V -2} ,\dots
, (-1)^{1} \,  e_{1} e_{1} \, e_{2}e_{2}\, \dots \, e_{_{\dim V}}
e_{_{\dim V}} \notag \\
=\,&(-1) ^{p}(-1)^{\dim V -1}  (-1)^{\dim V -2} ,\dots
, (-1)^{1} \, (-1)^{\dim V}\, 1_{_{C(V)}}  \notag \\
=\,& (-1)^{ ( \, p\,+\,\frac{\dim V\cdot (\dim V -1)}{2} \, +
 \, \dim V\, )}\, 1_{_{C(V)}}   \notag \\
=\,& 1_{_{C(V)}} \notag
\end{align}
where the last equality may be checked separately for
$\dim V \equiv 1, 2, 3$ or $ 0 \,  \mod 4$.

Similarly, for any basis element $e_{i} \in V$,
\begin{align}
\Gamma\, e_{i} \, = \, & \ii ^{p} \, e_{1} \dots e_{_{\dim V}}
\, e_{i} \notag \\
=\,&(-1)^{(\dim V -i)} \, \ii ^{p} \, e_{1} \dots e_{i}e_{i}
\dots e_{_{\dim V}} \notag \\
=\,&(-1)^{(\dim V -i)}\,(-1)^{(i-1)} \, \ii ^{p} \, e_{i} \, e_{1}
\dots e_{i} \dots e_{_{\dim V}} \notag \\
=\,& (-1)^{(\dim V -1)}\, e_{i} \, \Gamma \notag
\end{align}
Hence the equation is true for any $v \in V$.
\end{pf}

\begin{Cor}
For $\dim V$ odd, $\Gamma$ is in the center of $C(V)
\otimes \CC$.
\end{Cor}

\begin{Cor}
For $\dim V$ even, every complex Clifford module $E$ has a
 $\ZZ_{2}$-grading defined by the
$\pm 1$ eigen-spaces of the chirality operator:
\[ E^{\pm} \, = \, \{ v \in E \, | \, \Gamma v \, = \, \pm v \, \} \]
Especially, for $\dim V \equiv 0 \, \mod 4, \Gamma
\in C(V)$, so in this case,
 the real Clifford modules are also $\ZZ_{2}$-graded.
\end{Cor}

\begin{Def}
A {\bf polarization} of a complex vector space $V\otimes
\CC$ is a subspace $P \subset V\otimes \CC$
which is {\bf isotropic}, ie.
\[Q ( \, v \, , \, v \, ) \,= \, 0 \qquad \forall v \in P , \]
and we have a spliting
\[ V \,\otimes \,\CC \, = \, P \, \oplus \, \overline{P}.  \]
Here the quadratic form $Q$ extends from $V$ to
$V\otimes \CC$ {\bf complex linearly}, ie.
\[ Q(\, a+\ii b\, , \, c+\ii d\, )\, =\, Q(\, a \, , \, c\, )\, + \, \ii \,
\bigl( \, Q(\, a \, , \, d\, )\, + \,
Q(\, b \, , \, c\, )\, \bigr)\, - \, Q(\, b \, , \, d\, ) \]
\end{Def}
\begin{Def}
A polarization is called {\bf oriented}, if there is an oriented
orthonormal basis $\{ e_{i} \}$
of $V$, such that $P$ is spanned by the vectors
\[ \{ \,w_{i}\, = \,\frac{(\, e_{2i-1} \, - \, \ii e_{2i} \,)}{\sqrt{2\,}}
\, | \, 1 \leq i \leq \frac{\dim V}{2}  \} \]
and therefore the complement $\overline{P}$ is
 spanned by the vectors
\[ \{ \,\overline{w}_{i}\, = \,\frac{(\, e_{2i-1} \, + \,
 \ii e_{2i} \,)}{\sqrt{2\,}}\, | \, 1 \leq i \leq \frac{\dim V}{2}  \} \]
\end{Def}

\begin{Lem} \label{L:wiwj}
The basis $\{ w_{i} \}_{i=1, \dots , \frac{_{\dim V}}{2} }$
and the corresponding conjugate
$ \{ \overline{w}_{i} \}_{i=1, \dots , \frac{_{\dim V}}{2} }$
satisfy the folloowing equations:
For $1 \,\leq \, i \,\leq \, \frac{\dim V}{2}, $
\[ w_{i}w_{i} \, = \, \overline{w}_{i}\overline{w}_{i}\, = \, 0 \]
\[ w_{i} \overline{w_{i}} \, + \, \overline{w_{i}} w_{i} \, = \, -2  \]
For $1 \,\leq \, i \, \neq \, j \,\leq \, \frac{\dim V}{2}, $
\[ w_{i}w_{j} \, = \, - \, w_{j}w_{i} \]
\[ w_{i}\overline{w}_{j} \, = \, - \, \overline{w}_{j}w_{i} \]
\[  \overline{w}_{i}\overline{w}_{j}  \, = \, - \, \overline{w}_{j}
\overline{w}_{i}  \]
\end{Lem}

\begin{pf}
For $1 \,\leq \, i \,\leq \, \frac{\dim V}{2}, $
\begin{align}
w_{i}w_{i} \, = \,&\frac{(\, e_{2i-1} \, - \, \ii e_{2i} \,)}{\sqrt{2\,}}\,
 \frac{(\, e_{2i-1} \, - \, \ii e_{2i} \,)}{\sqrt{2\,}} \notag \\
=\,&\frac{1}{2}\, ( \, e_{2i-1}e_{2i-1} \, -\, e_{2i}e_{2i}\, -
 \, \ii \, e_{2i}e_{2i-1}\, - \, \ii \, e_{2i-1}e_{2i} \, ) \notag \\
=\,&\frac{1}{2}\, ( \,(-1) \, - \, (-1)\,  + \, \ii \, e_{2i-1}e_{2i}\,
- \, \ii \, e_{2i-1}e_{2i} \, ) \notag \\
=\,&\, 0 \notag
\end{align}
and
\[ \overline{w}_{i}\, \overline{w}_{i}\, =
\, (\overline{ w_{i}w_{i} }) \, = \, 0 \]
Also,
\begin{align}
w_{i} \overline{w}_{i} \, + \, \overline{w}_{i} w_{i} \, =\, &
\frac{e_{2i-1} \, - \, \ii \, e_{2i}}{\sqrt{2\,}}\,\frac{e_{2i-1} \, +
\, \ii \, e_{2i}}{\sqrt{2\,}}\,+ \,
\frac{e_{2i-1} \, + \, \ii \, e_{2i}}{\sqrt{2\,}}\,\frac{e_{2i-1} \,
- \, \ii \, e_{2i}}{\sqrt{2\,}} \notag \\
=\,&\frac{1}{2}\, \bigl( 2\, e_{2i-1}e_{2i-1} \, +
 \,2\, e_{2i}e_{2i} bigr) \notag \\
=\,& -2 \notag
\end{align}
For $1 \,\leq \, i \, \neq \, j \,\leq \, \frac{\dim V}{2}, $
\begin{align}
w_{i}w_{j} \, = \,&\frac{(\, e_{2i-1} \, -
 \, \ii e_{2i} \,)}{\sqrt{2\,}}\, \frac{(\, e_{2j-1} \, -
 \, \ii e_{2j} \,)}{\sqrt{2\,}} \notag \\
=\,&\frac{1}{2}\, ( \, e_{2i-1}e_{2j-1} \,-
 \,  e_{2i}e_{2j}\, - \, \ii \, e_{2i}e_{2j-1}\, -
 \, \ii \, e_{2i-1}e_{2j} \, ) \notag \\
=\,&-\frac{1}{2}\, ( \, e_{2j-1}e_{2i-1} \, -
 \, e_{2j}e_{2i}\, - \, \ii \, e_{2j}e_{2i-1}\, -
 \, \ii \, e_{2j-1}e_{2i} \, ) \notag \\
=\,&-\,\frac{(\, e_{2j-1} \, -
 \, \ii e_{2j} \,)}{\sqrt{2\,}}\, \frac{(\, e_{2i-1} \,
- \, \ii e_{2i} \,)}{\sqrt{2\,}} \notag \\
=\,&-\, w_{j}w_{i} \notag
\end{align}
Similarly for the other equalities.
\end{pf}

\begin{Thm}
If $V$ is an {\bf even}-dimensional oriented Euclidean
 vector space, then there is a unique $\ZZ_{2}$-graded
Clifford module
\[S \, = \, S ^{+} \, \oplus \, S^{-} \]
called the {\bf spinor module}, such that
\[ C(V) \otimes \CC \, \cong \, \End (S) \]
Elements in $S^{+}$ and $S^{-}$ are called {\bf positive}
 and {\bf negative spinors} respectively.

In particular, we have
\[\dim_{_{\CC}} (S) \, = \, 2^{ (\frac{\dim V}{2})}\]
and
\[ \dim _{_{\CC}} (S^{+}) \, = \, \dim _{_{\CC}} (S_{-}) \,
= \, 2^{(\frac{\dim V}{2} -1 )}\]
\end{Thm}

\begin{pf}
Given a polarization $P$ of $V\otimes \CC$,
we may define the spinor module $S$  to be equal to the
exterior algebra $\Lambda P$.
Any element $v \in V \otimes \CC  \hookrightarrow
C(V) \otimes \CC$ splits into two parts
\[ v \, = \, v_{1} \, \oplus \, v_{2} \quad \in P \,
\oplus \, \overline{P}   \]
The component in $P$ acts on $\Lambda P$ by exterior product:
\[ c( \, v_{1} \, ) \, w \, = \, \sqrt{2\, }\, \epsilon(\, v_{1} \, )\,
w\, =\, \sqrt{2\, }\,  v_{1} \, \wedge\, w , \qquad  \forall w \in
 \Lambda P \]

On the other hand the quadratic form $Q$ induce a duality
between $P$ and $\overline{P}$, ie.
\[ w_{i} ^{\ast} \, = \, Q ( \, \overline{w}_{i} \, , \, \cdot\, ) \]
So, the component in $\overline{P} \cong P ^{\ast}$ acts  on
$\Lambda P$ by contraction:
\[ c( \, v_{2} \, ) \, w \, = \, -\sqrt{2\, }\, \iota(\, v_{2} \, )\,
w \qquad  \forall w \in \Lambda P \]
Therefore the action of any $v = v_{1} \oplus v_{2} \in V
\otimes \CC$ on $\Lambda P$ is
\[ c( \, v_{1} \oplus v_{2} \, ) \, w \, = \, \sqrt{2\, }\,
\epsilon(\, v_{1} \, )\, w\,  - \, \sqrt{2\, }\, \iota(\, v_{2} \, )\, w
   \qquad  \forall w \in \Lambda P \]

Observe that the map  $C(V) \otimes \CC \longrightarrow
 \End ( S )$ is injective and
\[ \dim_{_{\CC}}  (\, C(V) \otimes \CC\, )\, = \, 2^{\dim V}\,
 = \, \dim_{_{\CC}} ( \, S \, ) ^{2} \, = \, \dim_{_{\CC}} ( \, \End S \, ) \]
therefore
\[ C(V) \otimes \CC\, \cong \, \End (\, S \, ) \]

The uniqueness of the spinor representation is a consequence
 of the fact that the algebra $\End (S)$
is simple , and it has unique irreducible module.

Now, assume the polarization $P$ is oriented, then for
\[ w_{i} \, = \, \frac{e_{2i-1} \, - \, \ii \, e_{2i}}{\sqrt{2\,}}  \]
we have
\[e_{2i-1} \, = \, \frac{ (\, w_{i} \, +
 \, \overline{w}_{i}\,)}{\sqrt{2\,}} \qquad
e_{2i} \, = \, \frac{ (\, - w_{i} \, +
 \, \overline{w}_{i}\,)}{\sqrt{2\,} \ii} \]
therefore from the above Lemma, we get
\[ e_{2i-1}e_{2i} \, = \,  \frac{(\,w_{i}\overline{w}_{i} \, -
 \, \overline{w}_{i}w_{i}\,)}{2 \ii} \]
and hence

\[ \Gamma \, = \, \frac{(\,w_{1}\overline{w}_{1} \,
- \, \overline{w}_{1}w_{1}\,) \, (\,w_{2}\overline{w}_{2} \, -
 \, \overline{w}_{2}w_{2}\,) \dots
(\,w_{\frac{_{\dim V}}{2}}\overline{w}_{\frac{_{\dim V}}{2}} \,
- \, \overline{w}_{\frac{_{\dim V}}{2}}
w_{\frac{_{\dim V}}{2}}\,) }{\sqrt{2\,}^{_{\dim V}}} \]
Therefore, $\Gamma$ acts on $\Lambda P$ by
\[ (-1)^{\frac{_{\dim V}}{2}} \, ( \,\epsilon(w_{1})
\iota (\overline{w}_{1}) \, - \, \iota (\overline{w}_{1})
\epsilon(w_{1})\, ) \, \dots \,  ( \,\epsilon(w_{1})
\iota (\overline{w}_{1}) \, - \, \iota (\overline{w}_{1})
\epsilon(w_{1})\,) \]
Now
\[  ( \,\epsilon(w_{i}) \iota (\overline{w}_{i}) \, -
 \, \iota (\overline{w}_{i})\epsilon(w_{i})\, )\, 1_{_{\Lambda P}} \, =
 \, -1_{_{\Lambda P}} \]
So $\Gamma$ acts on $\Lambda ^{0} P$ by $1_{_{C(V)}}$.

On the other hand, for $\Lambda ^{1} P$,
\begin{align}
  ( \,\epsilon(w_{i}) \iota (\overline{w}_{i}) \, -
 \, \iota (\overline{w}_{i})\epsilon(w_{i})\, )\, w_{j} \,
= \,& (-1)^{1+ \delta _{ij}} \,w_{j} \notag \\
= \,&
\begin{cases}
w_{j} &\text{if $i=j$}\\
-w_{j} &\text{if $i \neq j$}
\end{cases} \notag
\end{align}
So  $\Gamma$ acts on $\Lambda ^{1} P$ by
\[ (-1)^{\frac{_{\dim V}}{2}}\, (-1)^{\frac{_{\dim V}}{2} -1 }
\,1_{_{C(V)}}\, = \, -1_{_{C(V)}}. \]

On $\Lambda ^{2} P$, we have, for any $1 \leq j < k \leq
\frac{_{\dim V}}{2}$,
\begin{align}
  ( \,\epsilon(w_{i}) \iota (\overline{w}_{i}) \, - \, \iota
 (\overline{w}_{i})\epsilon(w_{i})\, )\, w_{j} \wedge w_{k} \,
= \,& (-1)^{1+ \delta _{ij} + \delta _{ik}} \, w_{j} \wedge w_{k}
\notag \\
= \,&
\begin{cases}
w_{j}\wedge w_{k} &\text{if $i = j$ or $i=k$}\\
-w_{j}\wedge w_{k} &\text{if $j \neq i \neq k$}
\end{cases} \notag
\end{align}
So  $\Gamma$ acts on $\Lambda ^{2} P$ by
\[ (-1)^{\frac{_{\dim V}}{2}}\, (-1)^{\frac{_{\dim V}}{2} -2 }
\,1_{_{C(V)}}\, = \, 1_{_{C(V)}}. \]

In general, on $\Lambda ^{k} P$, we have,  for any $1
 \leq j_{1} < j_{2} < \dots < j_{k} \leq \frac{_{\dim V}}{2}$,
\begin{align}
( \,\epsilon(w_{i}) \iota (\overline{w}_{i}) \, -
 \, \iota (\overline{w}_{i})\epsilon(w_{i})\, )\,
w_{j_{1}} \wedge \dots \wedge w_{j_{k}} \, = \,&
 (-1)^{1+ \delta _{ij_{1}} + \dots +
\delta _{ij_{k}}} \, w_{j_{1}} \wedge \dots \wedge w_{j_{k}}
 \notag \\
=\,&
\begin{cases}
w_{j_{1}} \wedge \dots \wedge w_{j_{k}} &
\text{if $i = j_{r}$ for some $j_{r}$}\\
-w_{j_{1}} \wedge \dots \wedge w_{j_{k}} &
\text{if $i \neq j_{r}$ for any $j_{r}$}
\end{cases} \notag
\end{align}
So  $\Gamma$ acts on $\Lambda ^{k} P$ by
\[ (-1)^{\frac{_{\dim V}}{2}}\, (-1)^{\frac{_{\dim V}}{2} -k }\,
1_{_{C(V)}}\, = \, (-1)^{k}\,1_{_{C(V)}}. \]
 Therefore we have the splitting
\[ S^{\pm} \, = \, \Lambda ^{\pm} P \]
\end{pf}

\bigskip
\subsection{Four Dimensional Case}
Now we look at the spinors in four dimensional vector space.

Consider an oriented orthonormal basis $\{ e_{1}, e_{2}, e_{3},
e_{4} \}$ of the Euclidean vector
space $\RR ^{4}$. Then the standard oriented polarization $P$
 of $\CC ^{4} \cong \RR ^{4} \otimes \CC$
is generated by
\[\{ \,  w_{1} \, = \, \frac{e_{1} \, - \, \ii \, e_{2}}{\sqrt{2\,}} ,
\qquad w_{2} \, = \, \frac{e_{3} \, - \, \ii \, e_{4}}{\sqrt{2\,}} \, \} \]
and we have
\[ S^{+} \, = \, \sspan _{\CC} \langle \, 1_{_{\Lambda P}}\, ,
\, w_{1}\wedge w_{2}\, \rangle \]
\[ S^{-} \, = \, \sspan _{\CC} \langle \, w_{1}\, , \, w_{2}\, \rangle \]
Together, we have the following {\bf standard basis} of $S = S^{+}
 \oplus S^{-}$ :
\[ \{ \,  1_{_{\Lambda P}}\, , \, w_{1}\wedge w_{2}\, , \,  w_{1}\, ,
\, w_{2}\,\} \]

Under this basis, any spinor $s \in S$ may be written as a column
 vector with components
\[ s \, = \, \begin{pmatrix}
s_{1} \\
s_{2} \\
s_{3}\\
s_{4}
\end{pmatrix}  \]
and we have a spliting
\[ s \, = \, s^{+} \, \oplus \, s^{-} \, = \,  \begin{pmatrix}
s^{+}_{1} \\ \\
s^{+}_{2}
\end{pmatrix} \,  \oplus \,
 \begin{pmatrix}
s^{-}_{1} \\ \\
s^{-}_{2}
\end{pmatrix}  \, = \,  \begin{pmatrix}
s_{1} \\ \\
s_{2}
\end{pmatrix} \,  \oplus \,
 \begin{pmatrix}
s_{3} \\ \\
s_{4}
\end{pmatrix} \]

Since $S \, \cong \, \CC ^{4}$, we have a representation
\[ C(V) \, \otimes \, \CC \, \hookrightarrow \, \End (\CC ^{4})  \]
To find out the exact correspondencs, let's consider the Clifford
action of $w_{i}$ and $\overline{w}_{i}$
on $S$:
\begin{align}
 c( \, w_{1} \, ) \, : \, \begin{cases}
1_{_{\Lambda P}} \\
w_{1}\wedge w_{2} \\
w_{1} \\
w_{2}
\end{cases}
 \, &\mapsto \,
\begin{cases}
\sqrt{2\,}\, w_{1} \\
 0 \\
 0 \\
\sqrt{2\,}\, w_{1}\wedge w_{2}
\end{cases}\notag \\
 c( \, w_{2} \, ) \, : \, \begin{cases}
1_{_{\Lambda P}} \\
w_{1}\wedge w_{2} \\
w_{1} \\
w_{2}
\end{cases}
 \, &\mapsto \,
\begin{cases}
\sqrt{2\,}\, w_{2} \\
 0 \\
 - \,\sqrt{2\,}\,w_{1}\wedge w_{2} \\
0
\end{cases} \notag \\
 c( \, \overline{w}_{1} \, ) \, : \, \begin{cases}
1_{_{\Lambda P}} \\
w_{1}\wedge w_{2} \\
w_{1} \\
w_{2}
\end{cases}
 \, &\mapsto \,
\begin{cases}
 0 \\
- \,\sqrt{2\,}\, w_{2} \\
- \,\sqrt{2\,}\, 1_{_{\Lambda P}} \\
 0 \\
\end{cases} \notag \\
c( \, \overline{w}_{2} \, ) \, : \, \begin{cases}
1_{_{\Lambda P}} \\
w_{1}\wedge w_{2} \\
w_{1} \\
w_{2}
\end{cases}
 \, &\mapsto \,
\begin{cases}
 0 \\
\sqrt{2\,}\,w_{1} \\
 0 \\
- \,\sqrt{2\,}\, 1_{_{\Lambda P}}
\end{cases} \notag
\end{align}

Therefore we have the following correspondence
\begin{align}
 w_{1} \, &= \, \sqrt{2\,}\,
\begin{pmatrix}
\,0\, & \,0\, & \,0\, & \,0\, \\
0 & 0 & 0 & 1 \\
1 & 0 & 0 & 0 \\
0 & 0 & 0 & 0
\end{pmatrix} \notag \\
 w_{2} \, &= \, \sqrt{2\,}\,
\begin{pmatrix}
\,0\, & \,0\, & 0 & \,0\, \\
0 & 0 & -1 & 0 \\
0 & 0 & 0 & 0 \\
1  & 0 & 0 & 0
\end{pmatrix} \notag \\
\overline{w}_{1} \, &= \, \sqrt{2\,}\,
\begin{pmatrix}
\,0\, & 0 & -1 & \,0\, \\
0 & 0 & 0 & 0 \\
0 & 0 & 0 & 0 \\
0 & -1 & 0 & 0
\end{pmatrix} \notag \\
\overline{w}_{2} \, &= \, \sqrt{2\,}\,
\begin{pmatrix}
\,0\, & \,0\, &\,0\, & -1 \\
0 & 0 & 0 & 0 \\
0 & 1 & 0 & 0 \\
0 & 0 & 0 & 0
\end{pmatrix} \notag
\end{align}

Now by the relation
\[e_{2i-1} \, = \, \frac{ (\, w_{i} \, + \,
\overline{w}_{i}\,)}{\sqrt{2\,}} \qquad
e_{2i} \, = \, \frac{ (\, - w_{i} \, + \, \overline{w}_{i}\,)}{\sqrt{2\,} \ii}
\]
we can deduce the corresponding matrices :
\begin{align}
 e_{1} \, &= \,
\begin{pmatrix}
\,0\, & \,0\, &-1 & \, 0 \, \\
0 & 0 & 0 & 1 \\
1 & 0 & 0 & 0 \\
0 & -1 & 0 & 0
\end{pmatrix} \notag \\
 e_{2} \, &= \,
\begin{pmatrix}
\,0\, & \,0\, &\ii & \, 0 \, \\
0 & 0 & 0 & \ii \\
\ii & 0 & 0 & 0 \\
0 & \ii & 0 & 0
\end{pmatrix} \notag \\
 e_{3} \, &= \,
\begin{pmatrix}
\,0\, & \,0\, &\,0\, & -1 \\
0 & 0 & -1 & 0 \\
0 & 1 & 0 & 0 \\
1 & 0 & 0 & 0
\end{pmatrix} \notag \\
 e_{4} \, &= \,
\begin{pmatrix}
\,0\, & \,0\, &\,0\, & \ii \\
0 & 0 & -\ii & 0 \\
0 & -\ii & 0 & 0 \\
\ii & 0 & 0 & 0
\end{pmatrix} \notag
\end{align}

As a result, we have the following

\begin{Thm}    \label{T:s+s-lambda}
There are isomorphisms
\[ \Hom ( \, S ^{+} \, , \, S ^{-} \, ) \, \cong \, \Lambda ^{1}_{\CC}  \]
\[ \Hom ( \, S ^{-} \, , \, S ^{+} \, ) \, \cong \, \Lambda ^{1}_{\CC}  \]
\end{Thm}

\begin{pf}
With respect to the standard basis
\[ \{ \, 1_{_{\Lambda P}} \, , \, w_{1} \wedge w_{2} \, \} \]
of $S^{+}$ and the standard basis
\[ \{ \, w_{1} \, , \, w_{2} \, \} \]
of $S^{-}$ ,  we have
\begin{align}
 e_{1} | _{_{ \Hom ( S^{+} , S^{-} ) }} \, &= \,
\begin{pmatrix}
1 & 0  \\
0 & -1
\end{pmatrix} \notag \\
 e_{2} | _{_{ \Hom ( S^{+} , S^{-} ) }}\, &= \,
\begin{pmatrix}
\ii & 0  \\
0 & \ii
\end{pmatrix} \notag \\
 e_{3} | _{_{ \Hom ( S^{+} , S^{-} ) }}\, &= \,
\begin{pmatrix}
0 & 1  \\
1 & 0
\end{pmatrix} \notag \\
 e_{4} | _{_{ \Hom ( S^{+} , S^{-} ) }}\, &= \,
\begin{pmatrix}
0 & -\ii  \\
\ii & 0
\end{pmatrix} \notag
\end{align}
and any element in $\Hom ( \, S^{+} \, , \, S^{-} \, )$ is a matrix
\[ \begin{pmatrix}
A & B \\
C & D
\end{pmatrix} \]
which may be splited into
\begin{align}
\begin{pmatrix}
A & B \\
C & D
\end{pmatrix} \, = \, &\frac{(A-D)}{2}\,
\begin{pmatrix}
1 & 0 \\
0 & -1
\end{pmatrix} \,
- \, \ii\, \frac{(A+D)}{2} \,
\begin{pmatrix}
\ii & 0 \\
0 & \ii
\end{pmatrix} \notag \\
&+ \, \frac{(B+C)}{2}\,
\begin{pmatrix}
0 & 1 \\
1 & 0
\end{pmatrix} \,
+\, \ii \,\frac{(B-C)}{2} \,
\begin{pmatrix}
0 & -\ii \\
\ii & 0
\end{pmatrix} \notag \\
= \, &\frac{(A-D)}{2}\,
e_{1}| _{_{ \Hom ( S^{+} , S^{-} ) }} \,
- \, \ii\, \frac{(A+D)}{2} \,
e_{2}| _{_{ \Hom ( S^{+} , S^{-} ) }} \notag \\
&+ \, \frac{(B+C)}{2}\,
e_{3}| _{_{ \Hom ( S^{+} , S^{-} ) }} \,
+\, \ii \,\frac{(B-C)}{2} \,
e_{4}| _{_{ \Hom ( S^{+} , S^{-} ) }} \notag
\end{align}

On the other hand, with respect to the above standard basis, we have
\begin{align}
 e_{1} | _{_{ \Hom ( S^{-} , S^{+} ) }} \, &= \,
\begin{pmatrix}
-1 & 0  \\
0 & 1
\end{pmatrix} \notag \\
 e_{2} | _{_{ \Hom ( S^{-} , S^{+} ) }}\, &= \,
\begin{pmatrix}
\ii & 0  \\
0 & \ii
\end{pmatrix} \notag \\
 e_{3} | _{_{ \Hom ( S^{-} , S^{+} ) }}\, &= \,
\begin{pmatrix}
0 & -1  \\
-1 & 0
\end{pmatrix} \notag \\
 e_{4} | _{_{ \Hom ( S^{-} , S^{+} ) }}\, &= \,
\begin{pmatrix}
0 & \ii  \\
-\ii & 0
\end{pmatrix} \notag
\end{align}
and any element in $\Hom ( \, S^{-} \, , \, S^{+} \, )$ is a matrix
\[ \begin{pmatrix}
A & B \\
C & D
\end{pmatrix} \]
which may be splited into
\begin{align}
\begin{pmatrix}
A & B \\
C & D
\end{pmatrix} \, = \, &\frac{(D-A)}{2}\,
\begin{pmatrix}
-1 & 0 \\
0 & 1
\end{pmatrix} \,
- \, \ii\, \frac{(A+D)}{2} \,
\begin{pmatrix}
\ii & 0 \\
0 & \ii
\end{pmatrix} \notag \\
&- \, \frac{(B+C)}{2}\,
\begin{pmatrix}
0 & -1 \\
-1 & 0
\end{pmatrix} \,
+\, \ii \,\frac{(C-B)}{2} \,
\begin{pmatrix}
0 & \ii \\
-\ii & 0
\end{pmatrix} \notag \\
= \, &\frac{(D-A)}{2}\,
e_{1}| _{_{ \Hom ( S^{-} , S^{+} ) }} \,
- \, \ii\, \frac{(A+D)}{2} \,
e_{2}| _{_{ \Hom ( S^{-} , S^{+} ) }} \notag \\
&- \, \frac{(B+C)}{2}\,
e_{3}| _{_{ \Hom ( S^{-} , S^{+} ) }} \,
+\, \ii \,\frac{(C-B)}{2} \,
e_{4}| _{_{ \Hom ( S^{-} , S^{+} ) }} \notag
\end{align}
\end{pf}

The elements in $C^{2}(V)$ corresponds to the matrices
\begin{alignat}{2}
e_{1}e_{2} \, &= \,
\begin{pmatrix}
-\ii & 0 & 0 & 0 \\
0 & \,\ii\, & 0 & 0 \\
0 & 0 & \,\ii\, & 0 \\
0 & 0 & 0 & -\ii
\end{pmatrix} &
e_{3}e_{4} \, &= \,
\begin{pmatrix}
-\ii & 0 & 0 & 0 \\
0 & \,\ii\, & 0 & 0 \\
0 & 0 & -\ii & 0 \\
0 & 0 & 0 & \,\ii\,
\end{pmatrix} \notag\\
e_{1}e_{3} \, &= \,
\begin{pmatrix}
\,0\, & -1 & \,0\, & 0 \\
1 & 0 & 0 & 0 \\
0 & 0 & 0 & -1 \\
0 & 0 & 1 & 0
\end{pmatrix} &\qquad
e_{4}e_{2} \, &= \,
\begin{pmatrix}
\,0\, & -1 & 0 & \,0\,\\
1 & 0 & 0 & 0 \\
0 & 0 & 0 & 1 \\
0 & 0 & -1 & 0
\end{pmatrix} \notag\\
e_{1}e_{4} \,& = \,
\begin{pmatrix}
0 & \,\ii\, & 0 & 0 \\
\,\ii\, & 0 & 0 & 0 \\
0 & 0 & 0 & \,\ii\, \\
0 & 0 & \,\ii\, & 0
\end{pmatrix} &
e_{2}e_{3} \, &= \,
\begin{pmatrix}
0 & \,\ii\, & 0 & 0 \\
\,\ii\, & 0 & 0 & 0 \\
0 & 0 & 0 & -\ii \\
0 & 0 & -\ii & 0
\end{pmatrix} \notag
\end{alignat}

\begin{Thm}
There is a representation
\[ \spin (4) \, \longrightarrow \, \End (\CC ^{4}) \]
which splits into
\[ \su (2) \, \times \, \su (2) \]
\end{Thm}

\begin{pf}
{}From the above, we see that the action of elements in $\spin (4)$
 is given by
\begin{align}
&\cos t \,  +\, \sin t \, e_{1}e_{2} \, = \,
\begin{pmatrix}
\cos t\, -\,\ii \,\sin t & 0 & 0 & 0 \\
0 & \cos t\,+ \,\ii\,\sin t & 0 & 0 \\
0 & 0 & \cos t\,+ \,\ii\,\sin t & 0 \\
0 & 0 & 0 & \cos t\,-\,\ii \,\sin t
\end{pmatrix} \notag\\
&\cos t \,  +\, \sin t \,e_{3}e_{4} \, = \,
\begin{pmatrix}
\cos t\,-\,\ii\,\sin t & 0 & 0 & 0 \\
0 & \cos t\,+\,\ii\,\sin t & 0 & 0 \\
0 & 0 & \cos t\,-\,\ii\, \sin t& 0 \\
0 & 0 & 0 & \cos t\,+\,\ii\,\sin t
\end{pmatrix} \notag\\
&\cos t \,  +\, \sin t \,e_{1}e_{3} \, = \,
\begin{pmatrix}
\cos t & -\sin t & \,0\, & 0 \\
\sin t & \cos t & 0 & 0 \\
0 & 0 & \cos t & -\sin t \\
0 & 0 & \sin t & \cos t
\end{pmatrix} \notag\\
&\cos t \,  +\, \sin t \,e_{4}e_{2} \, = \,
\begin{pmatrix}
\cos t & -\sin t & 0 & \,0\,\\
\sin t & \cos t& 0 & 0 \\
0 & 0 & \cos t & \sin t \\
0 & 0 & -\sin t & \cos t
\end{pmatrix} \notag\\
&\cos t \,  +\, \sin t \,e_{1}e_{4} \, = \,
\begin{pmatrix}
\cos t & \ii\,\sin t & 0 & 0 \\
\ii\,\sin t & \cos t & 0 & 0 \\
0 & 0 & \cos t & \ii\,\sin t \\
0 & 0 & \ii\,\sin t & \cos t
\end{pmatrix} \notag\\
&\cos t \,  +\, \sin t \,e_{2}e_{3} \, = \,
\begin{pmatrix}
\cos t & \ii\,\sin t & 0 & 0 \\
\ii\,\sin t & \cos t & 0 & 0 \\
0 & 0 & \cos t & -\ii\,\sin t \\
0 & 0 & -\ii\,\sin t & \cos t
\end{pmatrix} \notag
\end{align}
Therefore the image of this representation splits into $\su (2)
\times \su (2)$.
In fact,  this map is compatible with the one we have in
{\bf Theorem \ref{T:spin4sp1su2}}.
\end{pf}
\bigskip
\subsection{(Anti) Self Duality}
There is an operator on the exterior algebra which is similar
to the chirality operator.

\begin{Def}
Let the {\bf  Hodge star operator}
\[\star : \Lambda V \longrightarrow \Lambda V\]
with respect to $Q$, be given by the following relation :
\[   e_{i_{1}} \wedge \dots \wedge e_{i_{k}}
       \wedge \star (  e_{i_{1}} \wedge \dots \wedge e_{i_{k}} )
      =   \epsilon _{i_{1}\dots i_{k}}        e_{1} \wedge \dots \wedge
 e_{_{\dim V}}    \]
where
\( \{ e_{i} \} _{i=1, \dots , \dim V} \) is orthonormal with respect to
$Q$ and
\[   \epsilon _{i_{1}\dots i_{k}} = sgn
    \begin{pmatrix}
     1 & 2 & \hdotsfor{3} & \dim V \\
    i_{1}  &  \dots  & i_{k}  &  i_{k+1}  &  \dots  &  i_{_{\dim V}}
    \end{pmatrix}   \]
and
\[  (  i_{k+1} ,\dots , i_{_{\dim V}} ) = (1, 2 , \dots , \widehat{i_{1}}  ,
\dots
    , \widehat{i_{k}} ,\dots , \dim V )    \]
\end{Def}

Now we consider the case when $\dim V = 4$.

\begin{Lem}
For a four dimensional vector space $V$, the restriction of the
 square of the star operator
satisfies:
\[ \star ^{2} | _{_{C^{\pm}(V)}} \, = \, \pm \,
1_{_{\End ( \, C^{ \pm}(V) \, )}}  \]
\end{Lem}

\begin{pf}
The star operator interchanges the following pair of elements:
\begin{align}
 1_{_{\Lambda V}}\, &\leftrightarrow \, e_{1}e_{2}e_{3}e_{4}
\notag\\
e_{1}e_{2} \,& \leftrightarrow e_{3}e_{4} \notag \\
e_{1}e_{3} \,& \leftrightarrow e_{4}e_{2} \notag \\
e_{1}e_{4} \,& \leftrightarrow e_{2}e_{3} \notag
\end{align}
therefore
\[ \star ^{2} | _{_{C^{+}(V)}} \, =  \, 1_{_{\End ( \, C^{+}(V) \, )}}  \]
On the other hand, on $C^{-}(V)$, we have
\begin{alignat}{3}
e_{1} \, & \overset{\star}{\mapsto} \, & e_{2}e_{3}e_{4} \,
 &\overset{\star}{\mapsto} \, -e_{1} \, \overset{\star}{\mapsto} \,&
-&e_{2}e_{3}e_{4} \notag \\
e_{2} \, & \overset{\star}{\mapsto} \, &- e_{1}e_{3}e_{4} \,
&\overset{\star}{\mapsto} \, -e_{2} \, \overset{\star}{\mapsto} \,&
&e_{1}e_{3}e_{4} \notag \\
e_{3} \, & \overset{\star}{\mapsto} \, & e_{1}e_{2}e_{4} \,
&\overset{\star}{\mapsto} \, -e_{3} \, \overset{\star}{\mapsto} \,&
-&e_{1}e_{2}e_{4} \notag \\
e_{4} \, & \overset{\star}{\mapsto} \, &- e_{1}e_{2}e_{3} \,
&\overset{\star}{\mapsto} \, -e_{4} \, \overset{\star}{\mapsto} \,&
&e_{1}e_{2}e_{3} \notag
\end{alignat}
\[ \star ^{2} | _{_{C^{-}(V)}} \, =  \, -1_{_{\End ( \, C^{-}(V) \, )}}  \]
\end{pf}

\begin{Cor}
For a four dimensional vector space $V$, $\Lambda^{2} V$
 splits into $\pm 1$ eigen-spaces of $\star$:
\[  \Lambda^{2}_{+} V \, = \, \{ \, \frac{(1+ \star )}{2} \, w \, | \,
w \in \Lambda V \, \} \]
and
\[  \Lambda^{2}_{-} V \, = \, \{ \, \frac{(1- \star )}{2} \, w \,
| \, w \in \Lambda V \, \} \]
which are called the space of {\bf self-dual ( SD )} or {\bf
anti-self-dual ( ASD ) } two-forms and are
 simplified as $\Lambda_{+}$ and $\Lambda_{-}$ respectively.
\end{Cor}

\begin{Def}
The standard basis for $\Lambda_{+}$ and $\Lambda_{-}$ are
 the {\bf self-dual  basis} :
\[ \{ \,\wwp_{1} \, = \, e_{1} \wedge e_{2}  +  e_{3} \wedge e_{4}\,
, \quad\wwp_{2} \, = \, e_{1} \wedge e_{3}  + e_{4} \wedge e_{2}\,
, \quad\wwp_{3} \, = \,e_{1} \wedge e_{4}  +  e_{2} \wedge e_{3}\, \} \]
and the {\bf anti-self-dual basis} :
\[ \{ \,\wwm_{1} \, = \, e_{1} \wedge e_{2}  -  e_{3} \wedge e_{4}\,
 , \quad \wwm_{2} \, = \,e_{1} \wedge e_{3}  -  e_{4} \wedge e_{2}\,
, \quad\wwm_{3} \, = \,e_{1} \wedge e_{4}  -  e_{2} \wedge e_{3}\, \} \]
respectively.
\end{Def}

\begin{Rem}
The corresponding elements $\q ( \, \text{\bf{w}}^{\pm}_{i} \, )$
 in $C(V)$
 are also denoted by the same symbols $\text{\bf{w}}^{\pm}_{i}$.
The push-forward of the star-operators
\[ \q _{\ast} ( \star ) : C^{2}(V) \longrightarrow C^{2}(V)  \]
induce the similar spliting
\[ C^{2}(V) \, = \, C^{2}_{+}(V) \, \oplus \, C^{2}_{-}(V) \]
and the standard basis for $C^{2}_{+}(V)$ and $C^{2}_{-}(V)$ are
\[ \{ \,\wwp_{1} \, = \, e_{1} e_{2}  +  e_{3}e_{4}\, , \quad\wwp_{2} \,
 = \, e_{1}e_{3}  + e_{4}e_{2}\, , \quad\wwp_{3} \, = \,e_{1}e_{4}  +
 e_{2}e_{3}\, \} \]
and
\[ \{ \,\wwm_{1} \, = \, e_{1}e_{2}  -  e_{3}e_{4}\, , \quad \wwm_{2} \,
 = \,e_{1}e_{3}  -  e_{4}e_{2}\, , \quad\wwm_{3} \, = \,e_{1}e_{4}  -
 e_{2}e_{3}\, \} \]
respectively.
\end{Rem}

\begin{Rem}
The corresponding elements $\q( \,\text{\bf{w}}^{\pm}_{i}\, ) $ in $C(V)
\otimes \CC$
 are also denoted by the same symbols $\text{\bf{w}}^{\pm}_{i} $.
We have the spliting
\[ C^{2}(V) \otimes \CC\, = \, C^{2}_{+}(V) \otimes \CC\, \oplus \,
 C^{2}_{-}(V)
 \otimes \CC \]
where the standard basis for $C^{2}_{+}(V)\otimes \CC$
and $C^{2}_{-}(V)\otimes \CC$ are also
\[ \{ \,\wwp_{1} \, = \, e_{1} e_{2}  +  e_{3}e_{4}\, ,
\quad\wwp_{2} \, = \, e_{1}e_{3}  + e_{4}e_{2}\, ,
\quad\wwp_{3} \, = \,e_{1}e_{4}  +  e_{2}e_{3}\, \} \]
and
\[ \{ \,\wwm_{1} \, = \, e_{1}e_{2}  -  e_{3}e_{4}\, ,
\quad \wwm_{2} \, = \,e_{1}e_{3}  -  e_{4}e_{2}\, ,
\quad\wwm_{3} \, = \,e_{1}e_{4}  -  e_{2}e_{3}\, \} \]
respectively.
\end{Rem}

\begin{Cor}
Written in terms of the basis $w_{i}$ and $\overline{w}_{i}$,
 the elements $\q( \,\text{\normalshape{\bf{w}}}^{\pm}_{i}\, )$
in $C(V) \otimes \CC$ are given by
\begin{align}
\wwp_{1} \, = \, &\ii \, ( \, 2 \, + \, \overline{w}_{1} w_{1} \, + \,
\overline{w}_{2} w_{2} \, ) \notag \\
\wwp_{2} \, = \, & w_{1} w_{2} \, + \, \overline{w}_{1}
\overline{w}_{2} \notag \\
\wwp_{3} \, = \, &\ii \, ( \, w_{1} w_{2} \, - \,
\overline{w}_{1} \overline{w}_{2} \, ) \notag \\
& \notag \\
\wwm_{1} \, = \, &\ii \, ( \, \overline{w}_{1} w_{1} \, - \,
\overline{w}_{2} w_{2} \, ) \notag \\
\wwm_{2} \, = \, & \overline{w}_{1} w_{2} \, + \, w_{1}
\overline{w}_{2} \notag \\
\wwm_{3} \, = \, &\ii \, ( \, \overline{w}_{1} w_{2} \, - \,
w_{1} \overline{w}_{2} \, ) \notag
\end{align}
\end{Cor}

\begin{pf}
By the relation
\[e_{1} \, = \, \frac{ (\, w_{1} \, + \, \overline{w}_{1}\,)}{\sqrt{2\,}}
\qquad
e_{2} \, = \, \frac{ (\, - w_{1} \, + \, \overline{w}_{1}\,)}{\sqrt{2\,} \ii}
\]
\[e_{3} \, = \, \frac{ (\, w_{2} \, + \, \overline{w}_{2}\,)}{\sqrt{2\,}}
\qquad
e_{4} \, = \, \frac{ (\, - w_{2} \, + \, \overline{w}_{2}\,)}{\sqrt{2\,} \ii}
\]
and  {\bf Lemma \ref{L:wiwj}}, we have
\begin{align}
e_{1}e_{2} \, = \, & \frac{1}{\,2 \ii \,}\, (\, w_{1} \, +
\, \overline{w}_{1}\,)\, (\, -w_{1} \, + \, \overline{w}_{1}\,) \notag \\
=\, & \frac{1}{\,2 \ii \,}\, (\, - {w_{1}} ^{2} \, +\, {\overline{w}_{1}}
^{2}\,
- \,\overline{w}_{1} w_{1} \, + \,w_{1} \overline{w}_{1} \,) \notag \\
=\,& \frac{\ii}{\,2 \,}\, ( \,\overline{w}_{1} w_{1} \, -
\,w_{1} \overline{w}_{1} \,) \notag \\
=\,& \frac{\ii}{\,2 \,}\, ( \,2\, \overline{w}_{1} w_{1} \, +\, 2 \, )
\notag \\
=\,& \ii\, ( \, \overline{w}_{1} w_{1} \, +\, 1 \, ) \notag
\end{align}
and similarly,
\begin{align}
e_{3}e_{4} \, = \, & \frac{1}{\,2 \ii \,}\, (\, w_{2} \, +
 \, \overline{w}_{2}\,)\, (\, -w_{2} \, + \, \overline{w}_{2}\,) \notag \\
=\,& \ii\, ( \, \overline{w}_{2} w_{2} \, +\, 1 \, ) \notag
\end{align}
So we have
\begin{align}
\wwp_{1} \, = \, &\ii \, ( \, 2 \, + \, \overline{w}_{1} w_{1} \,
+ \, \overline{w}_{2} w_{2} \, ) \notag \\
\wwm_{1} \, = \, &\ii \, ( \, \overline{w}_{1} w_{1} \, -
 \, \overline{w}_{2} w_{2} \, ) \notag
\end{align}
On the other hand, we have
\begin{align}
e_{1}e_{3} \, = \, & \frac{1}{\,2 \,}\, (\, w_{1} \,
+ \, \overline{w}_{1}\,)\, (\,w_{2} \, + \, \overline{w}_{2}\,) \notag \\
=\, & \frac{1}{\,2 \,}\, (\, w_{1}w_{2} \, +
\, \overline{w}_{1} \overline{w}_{2}\, +
 \,\overline{w}_{1} w_{2} \, + \,w_{1} \overline{w}_{2} \,) \notag \\
\notag \\
e_{4}e_{2} \, = \, & -\frac{1}{\,2 \,}\, (\, -w_{2} \, +
 \, \overline{w}_{2}\,)\, (\,-w_{1} \, + \, \overline{w}_{1}\,) \notag \\
=\, & -\frac{1}{\,2 \,}\, (\, w_{2}w_{1} \, +
\, \overline{w}_{2} \overline{w}_{1}\, - \,w_{2}\overline{w}_{1}  \, -
 \,\overline{w}_{2}w_{1}  \,) \notag \\
=\, & \frac{1}{\,2 \,}\, (\, w_{1}w_{2} \, +
\,  \overline{w}_{1}\overline{w}_{2}\, - \,\overline{w}_{1} w_{2} \,
- \,w_{1}\overline{w}_{2}  \,) \notag \\
\notag \\
e_{1}e_{4} \, = \, & \frac{1}{\,2 \ii\,}\, (\, w_{1} \, + \,
 \overline{w}_{1}\,)\,
 (\,-w_{2} \, + \, \overline{w}_{2}\,) \notag \\
=\, & \frac{1}{\,2 \ii\,}\, (\, -w_{1}w_{2} \, +\, \overline{w}_{1}
 \overline{w}_{2}\, - \,\overline{w}_{1} w_{2} \, +
 \,w_{1} \overline{w}_{2} \,) \notag \\
\notag \\
e_{2}e_{3} \, = \, & \frac{1}{\,2\ii \,}\, (\,- w_{1} \, +
\, \overline{w}_{1}\,)\, (\,w_{2} \, + \, \overline{w}_{2}\,) \notag \\
=\, & \frac{1}{\,2 \ii\,}\, (\, -w_{1}w_{2} \, +\, \overline{w}_{1}
\overline{w}_{2}\, + \,\overline{w}_{1} w_{2} \, -
\,w_{1} \overline{w}_{2} \,) \notag
\end{align}
So we have
\begin{align}
\wwp_{2} \, = \, & w_{1} w_{2} \, + \, \overline{w}_{1}
\overline{w}_{2} \notag \\
\wwp_{3} \, = \, &\ii \, ( \, w_{1} w_{2} \, -
\, \overline{w}_{1} \overline{w}_{2} \, ) \notag \\
\wwm_{2} \, = \, & \overline{w}_{1} w_{2} \, +
\, w_{1} \overline{w}_{2} \notag \\
\wwm_{3} \, = \, &\ii \, ( \, \overline{w}_{1} w_{2} \,
- \, w_{1} \overline{w}_{2} \, ) \notag
\end{align}
\end{pf}

\begin{Rem}
Notice that in the above corollary, the notation
$\text{\bf{w}}^{\pm}_{i}$ refers to elements
of $C(V) \otimes \CC$. If we consider $\text{\bf{w}}^{\pm}_{i}$
as elements of $\Lambda _{\CC}  V = \Lambda V \, \otimes \CC$,
then we have
\begin{align}
\wwp_{1} \, = \, &\ii \, ( \, \overline{w}_{1} \wedge w_{1} \, +
 \, \overline{w}_{2}  \wedge w_{2} \, ) \notag \\
\wwp_{2} \, = \, & w_{1} \wedge  w_{2} \, + \, \overline{w}_{1}
 \wedge \overline{w}_{2} \notag \\
\wwp_{3} \, = \, &\ii \, ( \, w_{1} \wedge  w_{2} \, -
 \, \overline{w}_{1}  \wedge \overline{w}_{2} \, ) \notag \\
& \notag \\
\wwm_{1} \, = \, &\ii \, ( \, \overline{w}_{1}  \wedge w_{1} \, -
 \, \overline{w}_{2}  \wedge w_{2} \, ) \notag \\
\wwm_{2} \, = \, & \overline{w}_{1} \wedge  w_{2} \, +
 \, w_{1}  \wedge \overline{w}_{2} \notag \\
\wwm_{3} \, = \, &\ii \, ( \, \overline{w}_{1} \wedge  w_{2} \,
- \, w_{1}  \wedge \overline{w}_{2} \, ) \notag
\end{align}
The difference in the two interpretations of $\wwp_{1}$
arises from the fact that
\[ w_{i}w_{i} \, = \, \overline{w}_{i} \overline{w}_{i} \, =
 \, -1_{_{C(V)}} \]
in $C(V) \otimes \CC$, but
\[ w_{i} \wedge w_{i} \, = \, \overline{w}_{i} \wedge
\overline{w}_{i} \, = \,0 \]
in $\Lambda_{\CC} V$.
\end{Rem}

\begin{Rem}
By using the self-dual and anti-self dual basis, we can
express our previous results again in a
``better'' way:

The map $ \tg : \spin (4) \longrightarrow \so (4)$ has the following
 images :
\[
\tg ( \, \exp ( \, t \, \wwp _{1} \, ) \, ) \, = \,
\begin{pmatrix}
\cos 2t & -\sin 2t & 0 &  0 \\
\sin 2t & \cos 2t & 0 & 0 \\
0 & 0 & \cos 2t & -\sin 2t \\
0 & 0 & \sin 2t & \cos 2t
\end{pmatrix}  \]
\[
\tg ( \, \exp ( \, t \, \wwp _{2} \, ) \, ) \, = \,
\begin{pmatrix}
\cos 2t & 0 & -\sin 2t &  0 \\
0 & \cos 2t & 0 & \sin 2t \\
\sin 2t & 0 & \cos 2t & 0 \\
0 & -\sin 2t & 0 & \cos 2t
\end{pmatrix}  \]
\[
\tg ( \, \exp ( \, t \, \wwp _{3} \, ) \, ) \, = \,
\begin{pmatrix}
\cos 2t & 0 & 0 &  -\sin 2t \\
0 & \cos 2t & -\sin 2t & 0 \\
0 & \sin 2t & \cos 2t & 0 \\
\sin 2t & 0 & 0 & \cos 2t
\end{pmatrix}   \]
\[
\tg ( \, \exp ( \, t \, \wwm _{1} \, ) \, ) \, = \,
\begin{pmatrix}
\cos 2t & -\sin 2t & 0 &  0 \\
\sin 2t & \cos 2t & 0 & 0 \\
0 & 0 & \cos 2t & \sin 2t \\
0 & 0 & -\sin 2t & \cos 2t
\end{pmatrix}  \]
\[
\tg ( \, \exp ( \, t \, \wwm _{2} \, ) \, ) \, = \,
\begin{pmatrix}
\cos 2t & 0 & -\sin 2t &  0 \\
0 & \cos 2t & 0 & -\sin 2t \\
\sin 2t & 0 & \cos 2t & 0 \\
0 & \sin 2t & 0 & \cos 2t
\end{pmatrix}  \]
\[
\tg ( \, \exp ( \, t \, \wwm _{3} \, ) \, ) \, = \,
\begin{pmatrix}
\cos 2t & 0 & 0 &  -\sin 2t \\
0 & \cos 2t & \sin 2t & 0 \\
0 & -\sin 2t & \cos 2t & 0 \\
\sin 2t & 0 & 0 & \cos 2t
\end{pmatrix}   \]

The map $ C^{+} (\RR ^{4}) \longrightarrow \HH \oplus \HH$ has images
\begin{alignat}{4}
  1 \,       & \mapsto &\, \1 \, &\oplus \, \1,      & \qquad \Gamma \, =
 \, - e_{1}e_{2}e_{3}e_{4} \, & \mapsto &\, \1 \, &\oplus \, -\1
\notag \\
 \wwp _{1} \, & \mapsto &\, -2 \, \ii \, &\oplus \, 0,
& \qquad \wwm _{1} \,        & \mapsto &\, 0\, &\oplus \,  2\,\ii
\notag \\
 \wwp _{2} \, & \mapsto &\, -2\,\jj \, &\oplus \, 0, & \qquad
 \wwm_{2}\,
 & \mapsto &\, 0 \, &\oplus \,  -2\,\jj  \notag \\
 \wwp _{3} \, & \mapsto &\, 2 \, \kk \, &\oplus \, 0,
 & \qquad \wwm _{3} \,        & \mapsto &\, 0 \, &\oplus \,  2\, \kk
\notag
\end{alignat}
and the inverse
$ \HH \, \oplus \, \HH \, \longrightarrow \, C^{+}(\RR ^{4}) $ has images
\begin{alignat}{4}
  \1 \,&\oplus\, 0 \, & & \mapsto \, \frac{1\,+\, \Gamma}{2} ,
& \qquad 0 &\oplus\, \1 \, & & \mapsto \, \frac{1 \, - \, \Gamma}{2}
 \notag \\
  \ii \,&\oplus\, 0 \, & & \mapsto \, \frac{-\, \wwp _{1}}{2} ,
& \qquad 0 &\oplus \,\ii \, & & \mapsto \, \frac{\wwm _{1}}{2}
\notag \\
  \jj \,&\oplus\, 0 \, & & \mapsto \, \frac{- \, \wwp _{2}}{2} ,
 & \qquad 0&\oplus \,\jj \,  & & \mapsto \, \frac{-\, \wwm_{2}}{2}
\notag \\
  \kk \,&\oplus \,0 \,& & \mapsto \, \frac{\wwp _{3}}{2} ,
& \qquad 0& \oplus \,\kk \,& & \mapsto \, \frac{\wwm _{3}}{2}\notag
\end{alignat}

So the map $\spin (4) \longrightarrow \su (2) \times \su (2)$ has
\[ \exp ( \, t \, \Lambda _{+} \, ) \, \longrightarrow \su (2)
\times 1_{_{\su (2)}} \]
\[ \exp ( \, t \, \Lambda _{-} \, ) \, \longrightarrow
1_{_{\su (2)}} \times \su (2) \]
given by
\begin{alignat}{2}
 \exp ( \, t \, \wwp _{1} \, ) \, &\mapsto \, &
\begin{pmatrix}
\exp (-\ii 2t)  & 0 \\
0 & \exp (\ii 2t)
\end{pmatrix}
\, &\oplus \,
\begin{pmatrix}
1  & 0 \\
0 & 1
\end{pmatrix} \notag\\
 \exp ( \, t \, \wwp _{2} \, ) \, &\mapsto \, &
\begin{pmatrix}
\cos 2t &- \sin 2t \\
\sin 2t & \cos 2t
\end{pmatrix}
\, &\oplus \,
\begin{pmatrix}
1  & 0 \\
0 & 1
\end{pmatrix} \notag\\
 \exp ( \, t \, \wwp _{3} \, ) \, &\mapsto \, &
\begin{pmatrix}
\cos 2t &\sin 2t \, \ii\\
\sin 2t \, \ii & \cos 2t
\end{pmatrix}
\, &\oplus \,
\begin{pmatrix}
1  & 0 \\
0 & 1
\end{pmatrix} \notag\\
 \exp ( \, t \, \wwm _{1} \, ) \, &\mapsto \, &
\begin{pmatrix}
1  & 0 \\
0 & 1
\end{pmatrix}
\, &\oplus \,
\begin{pmatrix}
\exp (\ii 2t)  & 0 \\
0 & \exp (-\ii 2t)
\end{pmatrix} \notag\\
 \exp ( \, t \, \wwm _{2} \, ) \, &\mapsto \, &
\begin{pmatrix}
1  & 0 \\
0 & 1
\end{pmatrix}
\, &\oplus \,
\begin{pmatrix}
\cos 2t &- \sin 2t \\
\sin 2t & \cos 2t
\end{pmatrix} \notag\\
 \exp ( \, t \, \wwm _{3} \, ) \, &\mapsto \, &
\begin{pmatrix}
1  & 0 \\
0 & 1
\end{pmatrix}
\, &\oplus \,
\begin{pmatrix}
\cos 2t &\sin 2t \, \ii\\
\sin 2t \, \ii & \cos 2t
\end{pmatrix} \notag
\end{alignat}

The action of $ \text{\normalshape{\bf{w}}}^{\pm}_{i}$ in $C(V)
 \otimes \CC$ on $ S \cong \Lambda P$
is given by
\begin{alignat}{2}
\wwp_{1} \, &= \, 2\,
\begin{pmatrix}
-\ii & 0 & 0 & 0 \\
0 & \,\ii\, & 0 & 0 \\
0 & 0 & \,0\, & 0 \\
0 & 0 & 0 & \,0\,
\end{pmatrix} &
\wwm_{1} \, &= \,  2\,
\begin{pmatrix}
\,0\, & 0 & 0 & 0 \\
0 & \,0\, & 0 & 0 \\
0 & 0 & \,\ii\, & 0 \\
0 & 0 & 0 & -\ii
\end{pmatrix} \notag\\
\wwp_{2} \, &= \,  2\,
\begin{pmatrix}
\,0\, & -1 & 0 & 0 \\
1 & \,0\, & 0 & 0 \\
0 & 0 & \,0\, & 0 \\
0 & 0 & 0 & \,0\,
\end{pmatrix} &\qquad
\wwm_{2} \, &= \,  2\,
\begin{pmatrix}
\,0\, & 0 & 0 & 0 \\
0 & \,0\, & 0 & 0 \\
0 & 0 & \,0\, & -1 \\
0 & 0 & 1 & \,0\,
\end{pmatrix} \notag\\
\wwp_{3} \, &= \,  2\,
\begin{pmatrix}
\,0\, & \ii & 0 & 0 \\
\ii & \,0\, & 0 & 0 \\
0 & 0 & \,0\, & 0 \\
0 & 0 & 0 & \,0\,
\end{pmatrix} &
\wwm_{3} \, &= \,  2\,
\begin{pmatrix}
\,0\, & 0 & 0 & 0 \\
0 & \,0\, & 0 & 0 \\
0 & 0 & \,0\, & \ii \\
0 & 0 & \ii & \,0\,
\end{pmatrix} \notag
\end{alignat}
therefore $\Lambda_{+}$ acts on $S^{+}$ while $\Lambda_{-}$
acts on $S^{-}$.
\end{Rem}

\begin{Rem}
Now consider the complexified algebra of self-dual and anti-self
dual two forms
\[ \Lambda_{\pm \CC} \, = \, \Lambda_{\pm} \otimes \CC \]

Usually, we use a basis
\begin{align}
\overline{w}_{1} \wedge w_{1} \, + \, \overline{w}_{2}
\wedge w_{2}\, = \, &- \, \ii \, \wwp_{1} \notag \\
w_{1} \wedge w_{2} \, = \,&\frac{\wwp_{2} \, - \, \ii \,
\wwp_{3}}{2} \notag \\
\overline{w}_{1}  \wedge \overline{w}_{2} \, = \,
&\frac{\wwp_{2} \, + \, \ii \, \wwp_{3}}{2} \notag
\end{align}
for $\Lambda_{+ \CC}$ and
\begin{align}
\overline{w}_{1}  \wedge w_{1} \, -
 \, \overline{w}_{2} \wedge  w_{2}\, =
 \, &- \, \ii \, \wwm_{1} \notag \\
\overline{w}_{1} \wedge w_{2} \, =
\,&\frac{\wwm_{2} \, - \, \ii \, \wwm_{3}}{2} \notag \\
w_{1}  \wedge \overline{w}_{2} \, =
 \,&\frac{\wwm_{2} \, + \, \ii \, \wwm_{3}}{2} \notag
\end{align}
for $\Lambda_{- \CC}$.

Their quantization acts on $S$ as
\begin{align}
\q(\overline{w}_{1} \wedge w_{1} \, + \, \overline{w}_{2}
 \wedge w_{2})\, = \, & 2\,
\begin{pmatrix}
-1 & 0 & 0 & 0 \\
0 & 1 & 0 & 0 \\
0 & 0 & 0 & 0 \\
0 & 0 & 0 & 0
\end{pmatrix}
 \notag \\
\q(w_{1} \wedge w_{2}) \, = \,&2 \,
\begin{pmatrix}
0 & 0 & 0 & 0 \\
1 & 0 & 0 & 0 \\
0 & 0 & 0 & 0 \\
0 & 0 & 0 & 0
\end{pmatrix}
 \notag \\
\q(\overline{w}_{1}  \wedge \overline{w}_{2}) \, = \,&2 \,
\begin{pmatrix}
0 & -1 & 0 & 0 \\
0 & 0 & 0 & 0 \\
0 & 0 & 0 & 0 \\
0 & 0 & 0 & 0
\end{pmatrix}
\notag \\
\notag \\
\q(\overline{w}_{1}  \wedge w_{1} \, - \, \overline{w}_{2}
\wedge  w_{2})\, = \, &2\,
\begin{pmatrix}
0 & 0 & 0 & 0 \\
0 & 0 & 0 & 0 \\
0 & 0 & 1 & 0 \\
0 & 0 & 0 & -1
\end{pmatrix}
 \notag \\
\q(\overline{w}_{1} \wedge w_{2}) \, = \,&2\,
\begin{pmatrix}
0 & 0 & 0 & 0 \\
0 & 0 & 0 & 0 \\
0 & 0 & 0 & 0 \\
0 & 0 & 1 & 0
\end{pmatrix}
\notag \\
\q(w_{1}  \wedge \overline{w}_{2}) \, = \,&2\,
\begin{pmatrix}
0 & 0 & 0 & 0 \\
0 & 0 & 0 & 0 \\
0 & 0 & 0 & -1 \\
0 & 0 & 0 & 0
\end{pmatrix}
\notag
\end{align}
\end{Rem}

As a result, we have the following:
\begin{Thm} \label{T:ends+sd}
There are isomorphisms
\[ \End ( \, S^{+} \, ) \, \cong \, \Lambda ^{0}_{\CC} \,
\oplus \, \Lambda _{+ \CC} \]
\[ \End ( \, S^{-} \, ) \, \cong \, \Lambda ^{0}_{\CC} \,
 \oplus \, \Lambda _{- \CC} \]
\end{Thm}

\begin{pf}
With respect to the standard basis
\[ \{ \, 1_{_{\Lambda P}} \, , \, w_{1} \wedge w_{2} \, \} \]
of $S^{+}$, any element in $\End ( S^{+} )$ is a matrix
\[ \begin{pmatrix}
A & B \\
C & D
\end{pmatrix} \]
which may be splited into
\begin{align}
\begin{pmatrix}
A & B \\
C & D
\end{pmatrix} \, = \, &\frac{(A+D)}{2}\,
\begin{pmatrix}
1 & 0 \\
0 & 1
\end{pmatrix} \,
+ \, \frac{(D-A)}{4} \,
\begin{pmatrix}
-2 & 0 \\
0 & 2
\end{pmatrix} \notag \\
&+ \, \frac{C}{2}\,
\begin{pmatrix}
0 & 0 \\
2 & 0
\end{pmatrix} \,
- \, \frac{B}{2} \,
\begin{pmatrix}
0 & -2 \\
0 & 0
\end{pmatrix} \notag \\
=\, &\frac{(A+D)}{2}\,
1_{_{\End (S^{+})}} \,
+ \, \frac{(D-A)}{4} \,
\q(\, \overline{w}_{1}\wedge w_{1} \, + \, \overline{w}_{2}
\wedge w_{2} \, )\, | _{ _{S^{+}}} \notag \\
&+ \, \frac{C}{2}\,
\q(w_{1}\wedge w_{2})\, | _{ _{S^{+}}} \,
- \, \frac{B}{2} \,
\q(\overline{w}_{1}\wedge \overline{w}_{2} )\, |  _{_{S^{+}}}\notag
\end{align}
Since $1_{_{\End (S^{+})}}\, = \,\frac{ 1 + \Gamma}{2}|_{_{S^{+}}}$,
we have
\[ \End (S^{+}) \, \cong \, \CC (\frac{ 1 + \Gamma}{2}) \, \oplus \,
\Lambda _{+ \CC} \]

On the other hand,with respect to the standard basis
\[ \{  \, w_{1} \, , \, w_{2} \, \} \]
of $S^{-}$, any element
\[ \begin{pmatrix}
A & B \\
C & D
\end{pmatrix} \]
in $\End ( S^{-} )$  may be splited into
\begin{align}
\begin{pmatrix}
A & B \\
C & D
\end{pmatrix} \, = \, &\frac{(A+D)}{2}\,
\begin{pmatrix}
1 & 0 \\
0 & 1
\end{pmatrix} \,
+ \, \frac{(A-D)}{4} \,
\begin{pmatrix}
2 & 0 \\
0 & -2
\end{pmatrix} \notag \\
&+ \, \frac{C}{2}\,
\begin{pmatrix}
0 & 0 \\
2 & 0
\end{pmatrix} \,
- \, \frac{B}{2} \,
\begin{pmatrix}
0 & -2 \\
0 & 0
\end{pmatrix} \notag \\
=\, &\frac{(A+D)}{2}\,
1_{_{\End (S^{-})}} \,
+ \, \frac{(A-D)}{4} \,
\q(\, \overline{w}_{1}\wedge w_{1} \, - \, \overline{w}_{2}
\wedge w_{2} \, )\, | _{_{S^{-}}} \notag \\
&+ \, \frac{C}{2}\,
\q(\overline{w}_{1}\wedge w_{2})\, | _{_{S^{-}}} \,
- \, \frac{B}{2} \,
\q(w_{1}\wedge \overline{w}_{2} )\, | _{_{S^{-}}}\notag
\end{align}
Since $1_{_{\End (S^{-})}}\, = \,\frac{ 1 - \Gamma}{2}|_{_{S^{-}}}$,
 we have
\[ \End (S^{-}) \, \cong \, \CC (\frac{ 1 - \Gamma}{2}) \, \oplus \,
\Lambda _{- \CC} \]
\end{pf}
\bigskip
\subsection{Hermitian Structure on the Spinors}
\begin{Def}
There is a canonical {\bf Hermitian structure} on the space of
positive spinors $S^{+}$ given by the
{\bf Hermitian inner product} $\la \cdot \, , \, \cdot \ra$ which
takes the value
\[ \la s^{+} \, , \, t^{+} \ra \, = \, \overline{s}^{+}_{1} t^{+}_{1} \, +
 \, \overline{s}^{+}_{2} t^{+}_{2} \]
on the spinors
\[ s^{+} \, = \, \begin{pmatrix}
s^{+}_{1} \\ \\
s^{+}_{2} \end{pmatrix}\quad \text{and} \quad
 t^{+} \, = \, \begin{pmatrix}
t^{+}_{1} \\ \\
t^{+}_{2} \end{pmatrix} \quad \in S^{+} \]
\end{Def}

\begin{Lem}
The above Hermitian form is $\spin (4)$ invariant.
\end{Lem}

\begin{pf}
The action of $\spin (4)$ on $S^{+}$ is the same as  $\su (2)$
action. Since $\su (2)$
preserves Hermitian inner products, therefore our Hermitian
product is invariant under $\spin (4)$.
\end{pf}

\begin{Def}
The dual vector space $S^{+\ast}$ consists of complex linear functionals
\[ \phi : S^{+} \longrightarrow \CC .\]
$S^{+\ast}$ is generated by the dual complex basis
\[ \{ \, 1_{_{C(V)}}^{\ast}\, , \,( w_{1} \wedge w_{2})^{\ast}\, \}\]
which satisfies
\begin{align}
1_{_{C(V)}}^{\ast} (\,1_{_{C(V)}}\,)\, &= \,1\notag \\
1_{_{C(V)}}^{\ast} ( \,w_{1} \wedge w_{2}\,)\, &= \,0\notag \\
( w_{1} \wedge w_{2})^{\ast} (\,1_{_{C(V)}}\,)\, &= \,0\notag \\
( w_{1} \wedge w_{2})^{\ast} ( \,w_{1} \wedge w_{2}\,)\, &= \,1\notag
\end{align}
\end{Def}

\begin{Def}
There is a {\bf Hermitian Riesz representation}
\[ S^{+} \, \overset{\cong}{\longrightarrow} S^{+\ast} \]
with the following identification
\[ \begin{pmatrix}
s^{+}_{1} \\ \\
s^{+}_{2} \end{pmatrix} \, \mapsto \, \la \begin{pmatrix}
s^{+}_{1} \\ \\
s^{+}_{2} \end{pmatrix}\, , \, \cdot \ra  \]
\end{Def}

\begin{Thm}
The Hermitian Riesz representation
\[ S^{+} \longrightarrow S^{+\ast}\]
is given by
\[ \begin{pmatrix}
s^{+}_{1} \\ \\
s^{+}_{2} \end{pmatrix} \, \mapsto \,
\begin{pmatrix}
\overline{s}^{+}_{1}\\ \\
\overline{s}^{+}_{2}\end{pmatrix}^{\ast} \]
That means
\[
\begin{pmatrix}
s^{+}_{1}\\ \\
s^{+}_{2}\end{pmatrix}^{\ast} \, = \,
\la \begin{pmatrix}
\overline{s}^{+}_{1} \\ \\
\overline{s}^{+}_{2} \end{pmatrix}\, , \, \cdot \ra  \]
\end{Thm}

\begin{pf}
Now assume
\[ 1_{_{C(V)}}^{\ast}\, = \, \, \la
 \begin{pmatrix}
s^{+}_{1} \\ \\
s^{+}_{2} \end{pmatrix}
\, , \, \cdot \ra  \]
We have
\begin{alignat}{4}
1\, = \,& 1_{_{C(V)}}^{\ast} (\,1_{_{C(V)}}\,)&\,= \,&  \la
 \begin{pmatrix}
s^{+}_{1} \\ \\
s^{+}_{2} \end{pmatrix}
\, , \, 1_{_{C(V)}}\ra &\, = \,&
\la
 \begin{pmatrix}
s^{+}_{1} \\ \\
s^{+}_{2} \end{pmatrix}
\, , \,  \begin{pmatrix}
1 \\ \\
0 \end{pmatrix}\ra &\, = \,& \overline{s}^{+}_{1} \notag \\
0\, = \,& 1_{_{C(V)}}^{\ast} (\,w_{1} \wedge w_{2}\,)&\,= \,&  \la
 \begin{pmatrix}
s^{+}_{1} \\ \\
s^{+}_{2} \end{pmatrix}
\, , \, w_{1} \wedge w_{2}\ra &\, = \,&
\la
 \begin{pmatrix}
s^{+}_{1} \\ \\
s^{+}_{2} \end{pmatrix}
\, , \,  \begin{pmatrix}
0 \\ \\
1 \end{pmatrix}\ra &\, = \,& \overline{s}^{+}_{2} \notag
\end{alignat}
Therefore
\[ 1_{_{C(V)}}^{\ast}\, = \, \, \la
 \begin{pmatrix}
1 \\  \\
0\end{pmatrix}
\, , \, \cdot \ra  \]

Similarly, assume
\[ (w_{1} \wedge w_{2})^{\ast}\, = \, \, \la
 \begin{pmatrix}
s^{+}_{1} \\ \\
s^{+}_{2} \end{pmatrix}
\, , \, \cdot \ra  \]
and we have
\begin{alignat}{4}
0\, = \,& (w_{1} \wedge w_{2})^{\ast} (\,1_{_{C(V)}}\,)&\,= \,&  \la
 \begin{pmatrix}
s^{+}_{1} \\ \\
s^{+}_{2} \end{pmatrix}
\, , \, 1_{_{C(V)}}\ra &\, = \,&
\la
 \begin{pmatrix}
s^{+}_{1} \\ \\
s^{+}_{2} \end{pmatrix}
\, , \,  \begin{pmatrix}
1 \\ \\
0 \end{pmatrix}\ra &\, = \,& \overline{s}^{+}_{1} \notag \\
1\, = \,& (w_{1} \wedge w_{2})^{\ast} (\,w_{1} \wedge w_{2}\,)&\,=
\,&  \la
 \begin{pmatrix}
s^{+}_{1} \\ \\
s^{+}_{2} \end{pmatrix}
\, , \, w_{1} \wedge w_{2}\ra &\, = \,&
\la
 \begin{pmatrix}
s^{+}_{1} \\ \\
s^{+}_{2} \end{pmatrix}
\, , \,  \begin{pmatrix}
0 \\ \\
1 \end{pmatrix}\ra &\, = \,& \overline{s}^{+}_{2} \notag
\end{alignat}
Therefore
\[ (w_{1} \wedge w_{2})^{\ast}\, = \, \, \la
 \begin{pmatrix}
0 \\ \\
1\end{pmatrix}
\, , \, \cdot \ra  \]
As a result,
\begin{align}
 (\, s^{+}_{1}1_{_{C(V)}}\, + \, s^{+}_{2}(w_{1} \wedge w_{2})\, )^{\ast} \,
= \,&s^{+}_{1}\, 1_{_{C(V)}}^{\ast}\, + \, s^{+}_{2}\, (w_{1}
 \wedge w_{2})^{\ast} \notag\\
= \,&s^{+}_{1}\, \la
 \begin{pmatrix}
1\\ \\
0\end{pmatrix}
\, , \, \cdot \ra  \, + \, s^{+}_{2}\, \la
 \begin{pmatrix}
0\\ \\
1\end{pmatrix}
\, , \, \cdot \ra   \notag\\
 = \,&
\la
 \begin{pmatrix}
\overline{s}^{+}_{1}\\ \\
\overline{s}^{+}_{2}\end{pmatrix}
\, , \, \cdot \ra  \notag
\end{align}
\end{pf}

\begin{Thm} \label{T:s+s+her}
By using the Hermitian Riesz representation, we have
\[ S^{+} \, \otimes \, S^{+} \, \cong \, \End ( \, S^{+} \, ) \]
\end{Thm}

\begin{pf}
The image of the Hermitian Riesz representation
\[  \begin{pmatrix}
a \\ \\
b \end{pmatrix} \, \mapsto \, \langle \,
 \begin{pmatrix}
a \\ \\
b \end{pmatrix} \, , \, \cdot \, \rangle \, = \,
\begin{pmatrix}
\overline{a} \\ \\
\overline{b} \end{pmatrix}^{\ast}
\]
induces the following map
\[  \begin{pmatrix}
a \\ \\
b \end{pmatrix} \, \otimes \,
\begin{pmatrix}
c \\ \\
d \end{pmatrix} \, \mapsto \,
\langle \,
\begin{pmatrix}
a \\ \\
b \end{pmatrix} \,  , \cdot \, \rangle \, \otimes \,
\begin{pmatrix}
c \\ \\
d \end{pmatrix}
\qquad  \in  \End ( \, S^{+} \, ) \]
such that for any
\[ \begin{pmatrix}
s_{1}^{+} \\ \\
s_{2}^{+} \end{pmatrix} \,
\in S^{+}, \]
\begin{align}
\langle \,
\begin{pmatrix}
a \\ \\
b \end{pmatrix} \,  , \cdot \, \rangle \, \otimes \,
\begin{pmatrix}
c \\ \\
d \end{pmatrix} \, : \,
 \begin{pmatrix}
s_{1}^{+} \\ \\
s_{2}^{+} \end{pmatrix} \,\mapsto \, &
\langle \,
\begin{pmatrix}
a \\ \\
b \end{pmatrix} \, ,\,
\begin{pmatrix}
s_{1}^{+} \\ \\
s_{2}^{+} \end{pmatrix}
\, \rangle \, \otimes \,
\begin{pmatrix}
c \\ \\
d \end{pmatrix} \notag \\
=\,&( \, \overline{a} s_{1}^{+} \, + \, \overline{b} s_{2}^{+} \, ) \,
\begin{pmatrix}
c \\ \\
d \end{pmatrix}\notag \\
=\,&
\begin{pmatrix}
 \overline{a}c \,  s_{1}^{+}  + \overline{b}c\, s_{2}^{+}  \\ \\
 \overline{a}d\, s_{1}^{+}  +  \overline{b}d\, s_{2}^{+} \end{pmatrix}
\notag \\
=\,&
\begin{pmatrix}
 \overline{a}c & &\overline{b}c   \\ \\
 \overline{a}d & &\overline{b}d \end{pmatrix} \,
\begin{pmatrix}
s_{1}^{+} \\ \\
s_{2}^{+} \end{pmatrix} \notag
\end{align}
Therefore
\[  \langle \,
\begin{pmatrix}
a \\ \\
b \end{pmatrix} \,  , \cdot \, \rangle \, \otimes \,
\begin{pmatrix}
c \\ \\
d \end{pmatrix} \, = \,
\begin{pmatrix}
 \overline{a}c & &\overline{b}c   \\ \\
 \overline{a}d\, & & \overline{b}d \end{pmatrix} \, = \,
\begin{pmatrix}
c \\ \\
d \end{pmatrix} \, ( \, \overline{a} \qquad \overline{b} \, )
\]
Explicitly, take a canonical basis of $S^{+} \otimes S^{+}$ consisting
of the following
four elements:
\begin{align}
1_{_{C(V)}} \, \otimes \, 1_{_{C(V)}} \, = \,&
\begin{pmatrix}
1 \\ \\
0 \end{pmatrix} \, \otimes \,
\begin{pmatrix}
1 \\ \\
0 \end{pmatrix}  \notag\\
1_{_{C(V)}} \, \otimes \,(w_{1} \wedge w_{2}) \, = \,&
\begin{pmatrix}
1 \\ \\
0 \end{pmatrix} \, \otimes \,
\begin{pmatrix}
0 \\ \\
1 \end{pmatrix}  \notag\\
(w_{1} \wedge w_{2}) \, \otimes \, 1_{_{C(V)}} \, = \,&
\begin{pmatrix}
0 \\ \\
1 \end{pmatrix} \, \otimes \,
\begin{pmatrix}
1 \\ \\
0 \end{pmatrix}  \notag\\
(w_{1} \wedge w_{2}) \, \otimes \, (w_{1} \wedge w_{2}) \, = \,&
\begin{pmatrix}
0 \\ \\
1 \end{pmatrix} \, \otimes \,
\begin{pmatrix}
0 \\ \\
1 \end{pmatrix}  \notag \end{align}
The isomorphism maps this basis to the following basis matrices of
$\End ( S^{+} )$:
\begin{align}
1_{_{C(V)}} \, \otimes \, 1_{_{C(V)}} \, \mapsto \,&
\begin{pmatrix}
1 & & 0  \\ \\
0 & & 0 \end{pmatrix} \notag \\
 1_{_{C(V)}} \, \otimes \,(w_{1} \wedge w_{2}) \, \mapsto \,&
 \begin{pmatrix}
0 & & 0 \\ \\
1 & & 0 \end{pmatrix}  \notag\\
(w_{1} \wedge w_{2}) \, \otimes \, 1_{_{C(V)}} \, \mapsto \,&
 \begin{pmatrix}
0 & & 1 \\ \\
0 & & 0 \end{pmatrix}  \notag\\
(w_{1} \wedge w_{2}) \, \otimes \, (w_{1} \wedge w_{2}) \, \mapsto \,&
 \begin{pmatrix}
0 & & 0 \\ \\
0 & & 1 \end{pmatrix}  \notag \end{align}
This map extends to the whole $S^{+} \otimes S^{+}$, {\bf anti-linear}
 in the 1st
$S^{+}$ and linear in the second $S^{+}$.
 \end{pf}

 \begin{Cor}
 By using the Hermitian Riesz representation, we have
 \[ S^{+} \, \otimes \, S^{+} \, \cong \, \Lambda ^{0}_{\CC} \, \oplus \,
\Lambda _{+ \CC} \]
\end{Cor}

\begin{pf}
By {\bf Theorem \ref{T:ends+sd}} and {\bf Theorem \ref{T:s+s+her} }
 we have the following
correspondence:
\begin{align}
1_{_{C(V)}} \, \otimes \, 1_{_{C(V)}} \, \mapsto \,&
\begin{pmatrix}
1 & & 0  \\ \\
0 &  & 0 \end{pmatrix} \notag \\
\mapsto \,&\frac{1}{2} \, 1_{_{\End ( S^{+} )}} \,
- \, \frac{1}{4} \, \q( \overline{w}_{1} \wedge w_{1} \, +
\, \overline{w}_{2} \wedge w_{2} \, )|_{_{S^{+}}}\notag \\
 1_{_{C(V)}} \, \otimes \,(w_{1} \wedge w_{2}) \, \mapsto \,&
 \begin{pmatrix}
0 & & 0 \\ \\
1 & & 0 \end{pmatrix}  \notag\\
\mapsto \,& \frac{1}{2} \, \q (w_{1} \wedge w_{2})|_{_{S^{+}}}\notag \\
(w_{1} \wedge w_{2}) \, \otimes \, 1_{_{C(V)}} \, \mapsto \,&
 \begin{pmatrix}
0 & & 1 \\ \\
0 & & 0 \end{pmatrix}  \notag\\
\mapsto \,&- \frac{1}{2} \, \q (\overline{w}_{1} \wedge
\overline{w}_{2})|_{_{S^{+}}}\notag \\
(w_{1} \wedge w_{2}) \, \otimes \, (w_{1} \wedge w_{2}) \, \mapsto \,&
 \begin{pmatrix}
0 & &  0 \\ \\
0 & & 1 \end{pmatrix}  \notag \\
\mapsto \,&\frac{1}{2} \, 1_{_{\End ( S^{+} )}} \,
+ \, \frac{1}{4} \, \q( \overline{w}_{1} \wedge w_{1} \, +
\, \overline{w}_{2} \wedge w_{2} \, )|_{_{S^{+}}}\notag
\end{align}
This map extends to the whole $S^{+} \otimes S^{+}$,
{\bf anti-linear} in the 1st
$S^{+}$ and linear in the second $S^{+}$.
\end{pf}
\bigskip
Now consider the negative spinors $S^{-}$. With respect to the
standard basis
\[ \{ \, w_{1}\, , \, w_{2} \, \} \]
we can define a similar {\bf $\spin (4)$ invariant
 Hermitian inner product} which takes the value
\[ \la s^{-} \, , \, t^{-} \ra \, = \, \overline{s}^{-}_{1} t^{-}_{1} \,
+ \, \overline{s}^{-}_{2} t^{-}_{2} \]
on the spinors
\[ s^{-} \, = \, \begin{pmatrix}
s^{-}_{1} \\ \\
s^{-}_{2} \end{pmatrix}\quad \text{and} \quad
 t^{-} \, = \, \begin{pmatrix}
t^{-}_{1} \\ \\
t^{-}_{2} \end{pmatrix} \quad \in S^{-} \]

We have a similar {\bf Hermitian Riesz representation} on $S^{-}$ :
\[ S^{-} \, \overset{\cong}{\longrightarrow} S^{-\ast} \]
with the following identification
\[ \begin{pmatrix}
s^{-}_{1} \\ \\
s^{-}_{2} \end{pmatrix} \, \mapsto \, \la \begin{pmatrix}
s^{-}_{1} \\ \\
s^{-}_{2} \end{pmatrix}\, , \, \cdot \ra  \]
which , as in $S^{+}$, satisfies
\[
\begin{pmatrix}
s^{-}_{1}\\ \\
s^{-}_{2}\end{pmatrix}^{\ast} \, = \,
\la \begin{pmatrix}
\overline{s}^{-}_{1} \\ \\
\overline{s}^{-}_{2} \end{pmatrix}\, , \, \cdot \ra  \]

As in $S^{+}$, we have

\begin{Thm} \label{T:s-s-her}
By using the Hermitian Riesz representation, we have
\[ S^{-} \, \otimes \, S^{-} \, \cong \, \End ( \, S^{-} \, ) \]
\end{Thm}

\begin{pf}
The proof is similar to the one for $S^{+}$ and we get
\[  \begin{pmatrix}
a \\ \\
b \end{pmatrix} \, \otimes \,
\begin{pmatrix}
c \\ \\
d \end{pmatrix} \, \mapsto \,
 \langle \,
\begin{pmatrix}
a \\ \\
b \end{pmatrix} \,  , \cdot \, \rangle \, \otimes \,
\begin{pmatrix}
c \\ \\
d \end{pmatrix} \, = \,
\begin{pmatrix}
 \overline{a}c & &\overline{b}c   \\ \\
 \overline{a}d\, & & \overline{b}d \end{pmatrix}
\]
Explicitly, we have
\begin{align}
w_{1} \, \otimes \, w_{1} \, \mapsto \,&
\begin{pmatrix}
1 & & 0  \\ \\
0 & & 0 \end{pmatrix} \notag \\
 w_{1} \, \otimes \,w_{2} \, \mapsto \,&
 \begin{pmatrix}
0 & & 0 \\ \\
1 & & 0 \end{pmatrix}  \notag\\
w_{2} \, \otimes \, w_{1} \, \mapsto \,&
 \begin{pmatrix}
0 & & 1 \\ \\
0 & & 0 \end{pmatrix}  \notag\\
w_{2} \, \otimes \, w_{2} \, \mapsto \,&
 \begin{pmatrix}
0 & & 0 \\ \\
0 & & 1 \end{pmatrix}  \notag \end{align}
This map extends to the whole $S^{-} \otimes S^{-}$,
{\bf anti-linear} in the 1st
$S^{-}$ and linear in the second $S^{-}$.
\end{pf}

 \begin{Cor}
 By using the Hermitian Riesz representation, we have
 \[ S^{-} \, \otimes \, S^{-} \, \cong \, \Lambda ^{0}_{\CC} \,
 \oplus \, \Lambda _{- \CC} \]
\end{Cor}

\begin{pf}
By {\bf Theorem \ref{T:ends+sd}} and {\bf Theorem
\ref{T:s-s-her} } we have the following
correspondence:
\begin{align}
w_{1} \, \otimes \, w_{1} \, \mapsto \,&
  \begin{pmatrix}
  1 & & 0  \\ \\
  0 &  & 0 \end{pmatrix} \notag \\
 \mapsto \,&\frac{1}{2} \, 1_{_{\End ( S^{-} )}} \,
  + \, \frac{1}{4} \, \q( \overline{w}_{1} \wedge w_{1} \, -
 \, \overline{w}_{2} \wedge w_{2} \, )|_{_{S^{-}}}\notag \\
w_{1} \, \otimes \,w_{2} \, \mapsto \,&
  \begin{pmatrix}
  0 & & 0 \\ \\
  1 & & 0 \end{pmatrix}  \notag\\
 \mapsto \,& \frac{1}{2} \, \q (\overline{w}_{1}
\wedge w_{2} ) | _{_{S^{-}}} \notag \\
w_{2} \, \otimes \, w_{1} \, \mapsto \,&
  \begin{pmatrix}
  0 & & 1 \\ \\
  0 & & 0 \end{pmatrix}  \notag\\
 \mapsto \,&- \frac{1}{2} \, \q (w_{1} \wedge
\overline{w}_{2})|_{_{S^{-}}}\notag \\
w_{2} \, \otimes \,  w_{2} \, \mapsto \,&
  \begin{pmatrix}
  0 & & 0 \\ \\
  0 & &1 \end{pmatrix}  \notag \\
 \mapsto \,&\frac{1}{2} \, 1_{_{\End ( S^{-} )}} \,
- \, \frac{1}{4} \, \q( \overline{w}_{1} \wedge w_{1} \, -
 \, \overline{w}_{2} \wedge w_{2} \, )|_{_{S^{-}}}\notag
\end{align}
This map extends to the whole $S^{-} \otimes S^{-}$,
 {\bf anti-linear} in the 1st
$S^{-}$ and linear in the second $S^{-}$.
\end{pf}

\begin{Thm} \label{T:s+s-her}
By using the Hermitian Riesz representation, we have
\[ S^{+} \, \otimes \, S^{-} \, \overset{\cong}{\longrightarrow}
\, \Hom ( \, S^{+} \, , \, S^{-} \, ) \]
\end{Thm}

\begin{pf}
With respect to the standard basis of $S^{+}$ and $S^{-}$, we have a map
\[ S^{+} \, \otimes \, S^{-} \, \cong \, \Hom ( \, S^{+} \, , \, S^{-} \, ) \]
by
\[  \begin{pmatrix}
a \\ \\
b \end{pmatrix} \, \otimes \,
\begin{pmatrix}
c \\ \\
d \end{pmatrix} \, \mapsto \,
 \langle \,
\begin{pmatrix}
a \\ \\
b \end{pmatrix} \,  , \cdot \, \rangle \, \otimes \,
\begin{pmatrix}
c \\ \\
d \end{pmatrix} \, = \,
\begin{pmatrix}
 \overline{a}c & &\overline{b}c   \\ \\
 \overline{a}d\, & & \overline{b}d \end{pmatrix}
\]
where
\[  \begin{pmatrix}
a \\ \\
b \end{pmatrix} \, \in S^{+} , \quad \text{and} \quad
\begin{pmatrix}
c \\ \\
d \end{pmatrix} \, \in S^{-}  \]
and $\la \, \cdot \, , \, \cdot \, \ra$ is the Hermitian product of $S^{+}$.

Explicitly, we have
\begin{align}
1_{_{C(V)}} \, \otimes \, w_{1} \, \mapsto \,&
\begin{pmatrix}
1 & & 0  \\ \\
0 & & 0 \end{pmatrix} \notag \\
1_{_{C(V)}} \, \otimes \,w_{2} \, \mapsto \,&
 \begin{pmatrix}
0 & & 0 \\ \\
1 & & 0 \end{pmatrix}  \notag\\
(w_{1} \wedge w_{2}) \, \otimes \, w_{1} \, \mapsto \,&
 \begin{pmatrix}
0 & & 1 \\ \\
0 & & 0 \end{pmatrix}  \notag\\
(w_{1} \wedge w_{2}) \, \otimes \, w_{2} \, \mapsto \,&
 \begin{pmatrix}
0 & & 0 \\ \\
0 & & 1 \end{pmatrix}  \notag \end{align}
This map extends to the whole $S^{+} \otimes S^{-}$, {\bf anti-linear} in
$S^{+} \otimes 0$ and linear in   $0 \otimes S^{-}$.
\end{pf}

 \begin{Cor}
 By using the Hermitian Riesz representation, we have
 \[ S^{+} \, \otimes \, S^{-} \, \cong \, \Lambda ^{1}_{\CC}  \]
\end{Cor}

\begin{pf}
By {\bf Theorem \ref{T:s+s-lambda}} and {\bf Theorem \ref{T:s+s-her} }
 we have the following
correspondence:
\begin{alignat}{4}
1_{_{C(V)}} \, &\otimes \, w_{1} \,& \mapsto \,&
  \begin{pmatrix}
  1 & & 0  \\ \\
  0 &  & 0 \end{pmatrix} \,&
 \mapsto \,&\frac{1}{2} \, e_{1} \, - \, \frac{\ii}{2} \,  e_{2} &\,=\,&w_{1}
\notag \\
1_{_{C(V)}} \, &\otimes \,w_{2} \, &\mapsto \,&
  \begin{pmatrix}
  0 & & 0 \\ \\
  1 & & 0 \end{pmatrix}  \,&
 \mapsto \,&\frac{1}{2} \, e_{3} \, - \, \frac{\ii}{2} \,  e_{4} &\,=\,& w_{2}
\notag \\
(w_{1} \wedge w_{2} )\, &\otimes \, w_{1} \, &\mapsto \,&
  \begin{pmatrix}
  0 & & 1 \\ \\
  0 & & 0 \end{pmatrix}  \,&
 \mapsto \,&\frac{1}{2} \, e_{3} \, + \, \frac{\ii}{2} \,  e_{4} &\,
=\,&\overline{w}_{2}
\notag \\
(w_{1} \wedge w_{2}) \, &\otimes \,  w_{2} \, &\mapsto \,&
  \begin{pmatrix}
  0 & & 0 \\ \\
  0 & &1 \end{pmatrix}  \,&
 \mapsto \,&-\frac{1}{2} \, e_{1} \, - \, \frac{\ii}{2} \,  e_{2} &\,
=\,&-\overline{w}_{1}
\notag
\end{alignat}
This map extends to the whole $S^{+} \otimes S^{-}$, {\bf anti-linear} in
$S^{+} \otimes 0$ and linear in   $0 \otimes S^{-}$.
\end{pf}

Similarly, by interchanging $S^{+}$ and $S^{-}$, we get the following
\begin{Thm} \label{T:s-s+her}
By using the Hermitian Riesz representation, we have
\[ S^{-} \, \otimes \, S^{+} \, \cong \, \Hom ( \, S^{-} \, , \, S^{+} \, ) \]
\end{Thm}

\begin{pf}
With respect to the standard basis of $S^{-}$ and $S^{+}$, we have a map
\[ S^{-} \, \otimes \, S^{+} \, \overset{\cong}{\longrightarrow}
\, \Hom ( \, S^{-} \, , \, S^{+} \, ) \]
by
\[  \begin{pmatrix}
a \\ \\
b \end{pmatrix} \, \otimes \,
\begin{pmatrix}
c \\ \\
d \end{pmatrix} \, \mapsto \,
 \langle \,
\begin{pmatrix}
a \\ \\
b \end{pmatrix} \,  , \cdot \, \rangle \, \otimes \,
\begin{pmatrix}
c \\ \\
d \end{pmatrix} \, = \,
\begin{pmatrix}
 \overline{a}c & &\overline{b}c   \\ \\
 \overline{a}d\, & & \overline{b}d \end{pmatrix}
\]
where
\[  \begin{pmatrix}
a \\ \\
b \end{pmatrix} \, \in S^{-} , \quad \text{and} \quad
\begin{pmatrix}
c \\ \\
d \end{pmatrix} \, \in S^{+}  \]
and $\la \, \cdot \, , \, \cdot \, \ra$ is the Hermitian product of $S^{-}$.

Explicitly, we have
\begin{align}
w_{1} \, \otimes \, 1_{_{C(V)}} \, \mapsto \,&
\begin{pmatrix}
1 & & 0  \\ \\
0 & & 0 \end{pmatrix} \notag \\
w_{1}\, \otimes (\,w_{1} \wedge w_{2}) \, \mapsto \,&
 \begin{pmatrix}
0 & & 0 \\ \\
1 & & 0 \end{pmatrix}  \notag\\
w_{2} \, \otimes \, 1_{_{C(V)}} \, \mapsto \,&
 \begin{pmatrix}
0 & & 1 \\ \\
0 & & 0 \end{pmatrix}  \notag\\
w_{2} \, \otimes \, (w_{1} \wedge w_{2}) \, \mapsto \,&
 \begin{pmatrix}
0 & & 0 \\ \\
0 & & 1 \end{pmatrix}  \notag \end{align}
This map extends to the whole $S^{-} \otimes S^{+}$, {\bf anti-linear} in
$S^{-} \otimes 0$ and linear in   $0 \otimes S^{+}$.
\end{pf}

 \begin{Cor}
 By using the Hermitian Riesz representation, we have
 \[ S^{-} \, \otimes \, S^{+} \, \cong \, \Lambda ^{1}_{\CC}  \]
\end{Cor}

\begin{pf}
By {\bf Theorem \ref{T:s+s-lambda}} and {\bf Theorem \ref{T:s-s+her} }
 we have the following
correspondence:
\begin{alignat}{4}
w_{1} \, &\otimes \, 1_{_{C(V)}} &\, \mapsto \,&
  \begin{pmatrix}
  1 & & 0  \\ \\
  0 &  & 0 \end{pmatrix} \,&
 \mapsto \,&-\frac{1}{2} \, e_{1} \,  - \, \frac{\ii}{2} \,  e_{2} &\,=
\,&-\overline{w}_{1}
\notag \\
w_{1} \, &\otimes \,(w_{1} \wedge w_{2}) &\, \mapsto \,&
  \begin{pmatrix}
  0 & & 0 \\ \\
  1 & & 0 \end{pmatrix}  \,&
 \mapsto \,&-\frac{1}{2} \, e_{3} \,  + \, \frac{\ii}{2} \,  e_{4} &\,
=\,&-w_{2}
\notag \\
 w_{2} \, &\otimes \, 1_{_{C(V)}} &\, \mapsto \,&
  \begin{pmatrix}
  0 & & 1 \\ \\
  0 & & 0 \end{pmatrix}  \,&
 \mapsto \,&-\frac{1}{2} \, e_{3} \,  - \, \frac{\ii}{2} \,  e_{4} &\,
=\,&-\overline{w}_{2}
\notag \\
 w_{2} \, &\otimes \, (w_{1} \wedge w_{2}) &\, \mapsto \,&
  \begin{pmatrix}
  0 & & 0 \\ \\
  0 & &1 \end{pmatrix}  \,&
 \mapsto \,&\frac{1}{2} \, e_{1} \, - \, \frac{\ii}{2} \,  e_{2} &\,=\,&w_{1}
\notag
\end{alignat}
This map extends to the whole $S^{-} \otimes S^{+}$, {\bf anti-linear} in
$S^{-} \otimes 0$ and linear in   $0 \otimes S^{+}$.
\end{pf}

\bigskip
\subsection{Symplectic Structure on the Spinors}

\begin{Def}
There is a canonical {\bf symplectic structure} on the space of positive
 spinors $S^{+}$ given by the
{\bf symplectic form} $\{\, \cdot \, , \, \cdot\, \}$ which takes the value
\[ \{ s^{+} \, , \, t^{+} \} \, = \, s^{+}_{1} t^{+}_{2} \, - \, s^{+}_{2}
t^{+}_{1} \]
on the spinors
\[ s^{+} \, = \, \begin{pmatrix}
s^{+}_{1} \\ \\
s^{+}_{2} \end{pmatrix}\quad \text{and} \quad
 t^{+} \, = \, \begin{pmatrix}
t^{+}_{1} \\ \\
t^{+}_{2} \end{pmatrix} \quad \in S^{+} \]
\end{Def}

\begin{Lem}
The above symplectic form is $\spin (4)$ invariant.
\end{Lem}

\begin{pf}
The action of $\spin (4)$ on $S^{+}$ is represented as a
$\sl (2 , \CC)$ action. Since $\sl (2, \CC)$
preserves simplectic forms, therefore our symplectic
 form is invariant under $\spin (4)$.
\end{pf}

\begin{Def}
There is a {\bf symplectic Riesz representation}
\[ S^{+} \, \overset{\cong}{\longrightarrow} S^{+\ast} \]
with the following identification
\[ \begin{pmatrix}
s^{+}_{1} \\ \\
s^{+}_{2} \end{pmatrix} \, \mapsto \, \{ \begin{pmatrix}
s^{+}_{1} \\ \\
s^{+}_{2} \end{pmatrix}\, , \, \cdot \}  \]
\end{Def}

\begin{Thm}
The symplectic Riesz representation
\[ S^{+} \longrightarrow S^{+\ast}\]
is given by
\[ \begin{pmatrix}
s^{+}_{1} \\ \\
s^{+}_{2} \end{pmatrix} \, \mapsto \,
\begin{pmatrix}
-s^{+}_{2} \\ \\
s^{+}_{1} \end{pmatrix}^{\ast} \]
That means
\[
\begin{pmatrix}
s^{+}_{1} \\ \\
s^{+}_{2} \end{pmatrix}^{\ast} \, = \,
\{ \begin{pmatrix}
s^{+}_{2} \\ \\
-s^{+}_{1} \end{pmatrix}\, , \, \cdot \}  \]
\end{Thm}

\begin{pf}
Now assume
\[ 1_{_{C(V)}}^{\ast}\, = \, \, \{
 \begin{pmatrix}
s^{+}_{1} \\ \\
s^{+}_{2} \end{pmatrix}
\, , \, \cdot \}  \]
We have
\begin{alignat}{4}
1\, = \,& 1_{_{C(V)}}^{\ast} (\,1_{_{C(V)}}\,)&\,= \,&  \{
 \begin{pmatrix}
s^{+}_{1} \\ \\
s^{+}_{2} \end{pmatrix}
\, , \, 1_{_{C(V)}}\} &\, = \,&
\{
 \begin{pmatrix}
s^{+}_{1} \\ \\
s^{+}_{2} \end{pmatrix}
\, , \,  \begin{pmatrix}
1 \\ \\
0 \end{pmatrix}\} &\, = \,& -s^{+}_{2} \notag \\
0\, = \,& 1_{_{C(V)}}^{\ast} (\,w_{1} \wedge w_{2}\,)&\,= \,&  \{
 \begin{pmatrix}
s^{+}_{1} \\ \\
s^{+}_{2} \end{pmatrix}
\, , \, w_{1} \wedge w_{2}\} &\, = \,&
\{
 \begin{pmatrix}
s^{+}_{1} \\ \\
s^{+}_{2} \end{pmatrix}
\, , \,  \begin{pmatrix}
0 \\ \\
1 \end{pmatrix}\} &\, = \,& s^{+}_{1} \notag
\end{alignat}
Therefore
\[ 1_{_{C(V)}}^{\ast}\, = \, \, \{
 \begin{pmatrix}
0 \\ \\
-1\end{pmatrix}
\, , \, \cdot \}  \]

Similarly, assume
\[ (w_{1} \wedge w_{2})^{\ast}\, = \, \, \{
 \begin{pmatrix}
s^{+}_{1} \\ \\
s^{+}_{2} \end{pmatrix}
\, , \, \cdot \}  \]
and we have
\begin{alignat}{4}
0\, = \,& (w_{1} \wedge w_{2})^{\ast} (\,1_{_{C(V)}}\,)&\,= \,&  \{
 \begin{pmatrix}
s^{+}_{1} \\ \\
s^{+}_{2} \end{pmatrix}
\, , \, 1_{_{C(V)}}\} &\, = \,&
\{
 \begin{pmatrix}
s^{+}_{1} \\ \\
s^{+}_{2} \end{pmatrix}
\, , \,  \begin{pmatrix}
1 \\ \\
0 \end{pmatrix}\} &\, = \,& -s^{+}_{2} \notag \\
1\, = \,& (w_{1} \wedge w_{2})^{\ast} (\,w_{1} \wedge w_{2}\,)&\,=
 \,&  \{
 \begin{pmatrix}
s^{+}_{1} \\ \\
s^{+}_{2} \end{pmatrix}
\, , \, w_{1} \wedge w_{2}\} &\, = \,&
\{
 \begin{pmatrix}
s^{+}_{1} \\ \\
s^{+}_{2} \end{pmatrix}
\, , \,  \begin{pmatrix}
0 \\ \\
1 \end{pmatrix}\} &\, = \,& s^{+}_{1} \notag
\end{alignat}
Therefore
\[ (w_{1} \wedge w_{2})^{\ast}\, = \, \, \{
 \begin{pmatrix}
1 \\ \\
0 \end{pmatrix}
\, , \, \cdot \}  \]
As a result,
\[ (\, s^{+}_{1}1_{_{C(V)}}\, + \, s^{+}_{2}(w_{1} \wedge w_{2})\, )^{\ast}
\, = \,  s^{+}_{1}\, 1_{_{C(V)}}^{\ast}\, + \, s^{+}_{2}\, (w_{1}
\wedge w_{2})^{\ast} \, = \,
\{
 \begin{pmatrix}
s^{+}_{2}\\ \\
-s^{+}_{1}\end{pmatrix}
\, , \, \cdot \}  \]
\end{pf}

Just like the Hermitian case, we also have
\begin{Thm} \label{T:s+s+sym}
By using the symplectic Riesz representation, we have
\[ S^{+} \, \otimes \, S^{+} \, \cong \, \End ( \, S^{+} \, ) \]
\end{Thm}

\begin{pf}
The image of the symplectic
 Riesz representation
\[  \begin{pmatrix}
a \\ \\
b \end{pmatrix} \, \mapsto \, \{ \,
 \begin{pmatrix}
a \\ \\
b \end{pmatrix} \, , \, \cdot \, \} \, = \,
\begin{pmatrix}
-b \\ \\
a \end{pmatrix}^{\ast}
\]
induces the following map
\[  \begin{pmatrix}
a \\ \\
b \end{pmatrix} \, \otimes \,
\begin{pmatrix}
c \\ \\
d \end{pmatrix} \, \mapsto \,
\{  \,
\begin{pmatrix}
a \\ \\
b \end{pmatrix} \,  , \cdot \, \} \, \otimes \,
\begin{pmatrix}
c \\ \\
d \end{pmatrix}
\qquad  \in  \End ( \, S^{+} \, ) \]
such that for any
\[ \begin{pmatrix}
s_{1}^{+} \\ \\
s_{2}^{+} \end{pmatrix} \,
\in S^{+}, \]
\begin{align}
\{  \,
\begin{pmatrix}
a \\ \\
b \end{pmatrix} \,  , \cdot \, \} \, \otimes \,
\begin{pmatrix}
c \\ \\
d \end{pmatrix} \, : \,
 \begin{pmatrix}
s_{1}^{+} \\ \\
s_{2}^{+} \end{pmatrix} \,\mapsto \, &
\{ \,
\begin{pmatrix}
a \\ \\
b \end{pmatrix} \, ,\,
\begin{pmatrix}
s_{1}^{+} \\ \\
s_{2}^{+} \end{pmatrix}
\, \} \, \otimes \,
\begin{pmatrix}
c \\ \\
d \end{pmatrix} \notag \\
=\,&( \, -b s_{1}^{+} \, + \, a s_{2}^{+} \, ) \,
\begin{pmatrix}
c \\ \\
d \end{pmatrix}\notag \\
=\,&
\begin{pmatrix}
 \-bc \,  s_{1}^{+} \, + \, ac\, s_{2}^{+}  \\ \\
 \-bd\, s_{1}^{+} \, + \, ad\, s_{2}^{+} \end{pmatrix}\notag \\
=\,&
\begin{pmatrix}
-bc & ac   \\ \\
-bd & ad \end{pmatrix} \,
\begin{pmatrix}
s_{1}^{+} \\ \\
s_{2}^{+} \end{pmatrix} \notag
\end{align}
Therefore
\[  \{ \,
\begin{pmatrix}
a \\ \\
b \end{pmatrix} \,  , \cdot \, \} \, \otimes \,
\begin{pmatrix}
c \\ \\
d \end{pmatrix} \, = \,
\begin{pmatrix}
-bc & ac   \\ \\
-bd & ad \end{pmatrix} \, = \,
\begin{pmatrix}
c \\ \\
d \end{pmatrix} \, ( \, -b \qquad a \, )
\]
Explicitly, the isomorphism maps this basis to the following basis
matrices of $\End ( S^{+} )$:
\begin{align}
1_{_{C(V)}} \, \otimes \, 1_{_{C(V)}} \, \mapsto \,&
\begin{pmatrix}
0 & & 1  \\ \\
0 & & 0 \end{pmatrix} \notag \\
 1_{_{C(V)}} \, \otimes \,(w_{1} \wedge w_{2}) \, \mapsto \,&
 \begin{pmatrix}
0 & & 0 \\ \\
0 & & 1 \end{pmatrix}  \notag\\
(w_{1} \wedge w_{2}) \, \otimes \, 1_{_{C(V)}} \, \mapsto \,&
 \begin{pmatrix}
-1 & & 0 \\  \\
0 & & 0 \end{pmatrix}  \notag\\
(w_{1} \wedge w_{2}) \, \otimes \, (w_{1} \wedge w_{2}) \, \mapsto \,&
 \begin{pmatrix}
0 & & 0 \\ \\
-1 & & 0 \end{pmatrix}  \notag \end{align}
 and extends linearly to the whole $S^{+} \otimes S^{+}$.
\end{pf}

 \begin{Cor}
 By using the symplectic Riesz representation, we have
 \[ S^{+} \, \otimes \, S^{+} \, \cong \, \Lambda ^{0}_{\CC} \, \oplus \,
\Lambda _{+ \CC} \]
\end{Cor}

\begin{pf}
By {\bf Theorem \ref{T:ends+sd}} and {\bf Theorem \ref{T:s+s+sym} }
 we have the following
correspondence:
\begin{align}
1_{_{C(V)}} \, \otimes \, 1_{_{C(V)}} \, \mapsto \,&
  \begin{pmatrix}
  0 & & 1  \\ \\
  0 & & 0 \end{pmatrix} \notag \\
  \mapsto \,&- \frac{1}{2} \, \q (\overline{w}_{1} \wedge
\overline{w}_{2})|_{_{S^{+}}}\notag \\
1_{_{C(V)}} \, \otimes \,(w_{1} \wedge w_{2}) \, \mapsto \,&
  \begin{pmatrix}
  0 & & 0 \\ \\
  0 & & 1 \end{pmatrix}  \notag\\
 \mapsto \,&\frac{1}{2} \, 1_{_{\End ( S^{+} )}} \,
  + \, \frac{1}{4} \, \q( \overline{w}_{1} \wedge w_{1} \, + \,
\overline{w}_{2} \wedge w_{2} \, )|_{_{S^{+}}}\notag \\
(w_{1} \wedge w_{2}) \, \otimes \, 1_{_{C(V)}} \, \mapsto \,&
  \begin{pmatrix}
  -1 & & 0 \\ \\
  0 & & 0 \end{pmatrix}  \notag\\
  \mapsto \,&-\frac{1}{2} \, 1_{_{\End ( S^{+} )}} \,
   + \, \frac{1}{4} \, \q( \overline{w}_{1} \wedge w_{1} \, +
\, \overline{w}_{2} \wedge w_{2} \, )|_{_{S^{+}}}\notag \\
(w_{1} \wedge w_{2}) \, \otimes \, (w_{1} \wedge w_{2}) \, \mapsto \,&
  \begin{pmatrix}
  0 & & 0 \\ \\
  -1 & & 0 \end{pmatrix}  \notag \\
  \mapsto \,&- \frac{1}{2} \, \q (w_{1} \wedge w_{2})|_{_{S^{+}}}\notag
\end{align}
\end{pf}

\bigskip
Now consider the negative spinors $S^{-}$. With respect to the
 standard basis
\[ \{ \, w_{1}\, , \, w_{2} \, \} \]
we can define a similar {\bf $\spin (4)$ invariant
 symplectic form} which
takes the value
\[ \{ s^{-} \, , \, t^{-} \} \, = \, s^{-}_{1} t^{-}_{2} \, - \, s^{-}_{2}
t^{-}_{1} \]
on the spinors
\[ s^{-} \, = \, \begin{pmatrix}
s^{-}_{1} \\ \\
s^{-}_{2} \end{pmatrix}\quad \text{and} \quad
 t^{-} \, = \, \begin{pmatrix}
t^{-}_{1} \\ \\
t^{-}_{2} \end{pmatrix} \quad \in S^{-} \]

We have a similar {\bf symplectic Riesz representation} on $S^{-}$ :
\[ S^{-} \, \overset{\cong}{\longrightarrow} S^{-\ast} \]
with the following identification
\[ \begin{pmatrix}
s^{-}_{1} \\ \\
s^{-}_{2} \end{pmatrix} \, \mapsto \, \{ \begin{pmatrix}
s^{-}_{1} \\ \\
s^{-}_{2} \end{pmatrix}\, , \, \cdot \}  \]
which , as in $S^{+}$, satisfies
\[
\begin{pmatrix}
s^{-}_{1}\\ \\
s^{-}_{2}\end{pmatrix}^{\ast} \, = \,
\{ \begin{pmatrix}
s^{-}_{2} \\ \\
-s^{-}_{1} \end{pmatrix}\, , \, \cdot \}  \]

As in $S^{+}$, we have
\begin{Thm} \label{T:s-s-sym}
By using the symplectic Riesz representation, we have
\[ S^{-} \, \otimes \, S^{-} \, \cong \, \End ( \, S^{-} \, ) \]
\end{Thm}

\begin{pf}
The proof is similar to the one for $S^{+}$ and we get
\[  \begin{pmatrix}
a \\ \\
b \end{pmatrix} \, \otimes \,
\begin{pmatrix}
c \\ \\
d \end{pmatrix} \, \mapsto \,
 \{ \,
\begin{pmatrix}
a \\ \\
b \end{pmatrix} \,  , \cdot \, \} \, \otimes \,
\begin{pmatrix}
c \\ \\
d \end{pmatrix} \, = \,
\begin{pmatrix}
-bc & & ac   \\ \\
-bd\, & & ad \end{pmatrix}
\]
Explicitly, we have
\begin{align}
w_{1} \, \otimes \, w_{1} \, \mapsto \,&
\begin{pmatrix}
0 & & 1  \\ \\
0 & & 0 \end{pmatrix} \notag \\
 w_{1} \, \otimes \,w_{2} \, \mapsto \,&
 \begin{pmatrix}
0 & & 0 \\ \\
0 & & 1 \end{pmatrix}  \notag\\
w_{2} \, \otimes \, w_{1} \, \mapsto \,&
 \begin{pmatrix}
-1 & & 0 \\ \\
0 & & 0 \end{pmatrix}  \notag\\
w_{2} \, \otimes \, w_{2} \, \mapsto \,&
 \begin{pmatrix}
0 & & 0 \\ \\
-1 & & 0 \end{pmatrix}  \notag \end{align}
This map extends to the whole $S^{-} \otimes S^{-}$, {\bf anti-linear}
in the 1st
$S^{-}$ and linear in the second $S^{-}$.
\end{pf}

 \begin{Cor}
 By using the symplectic Riesz representation, we have
 \[ S^{-} \, \otimes \, S^{-} \, \cong \, \Lambda ^{0}_{\CC} \,
\oplus \, \Lambda _{- \CC} \]
\end{Cor}

\begin{pf}
By {\bf Theorem \ref{T:ends+sd}} and {\bf Theorem
\ref{T:s-s-sym} } we have the following
correspondence:
\begin{align}
w_{1} \, \otimes \, w_{1} \, \mapsto \,&
  \begin{pmatrix}
   0 & & 1  \\ \\
   0 &  & 0 \end{pmatrix} \notag \\
  \mapsto \,&- \frac{1}{2} \, \q (w_{1} \wedge
\overline{w}_{2})|_{_{S^{-}}}\notag \\
w_{1} \, \otimes \,w_{2} \, \mapsto \,&
  \begin{pmatrix}
   0 & & 0 \\ \\
   0 & & 1 \end{pmatrix}  \notag\\
  \mapsto \,&\frac{1}{2} \, 1_{_{\End ( S^{-} )}} \,
   - \, \frac{1}{4} \, \q( \overline{w}_{1} \wedge w_{1} \, -
 \, \overline{w}_{2} \wedge w_{2} \, )|_{_{S^{-}}}\notag \\
w_{2} \, \otimes \, w_{1} \, \mapsto \,&
  \begin{pmatrix}
  -1 & & 0 \\ \\
  0 & & 0 \end{pmatrix}  \notag\\
 \mapsto \,&-\frac{1}{2} \, 1_{_{\End ( S^{-} )}} \,
   - \, \frac{1}{4} \, \q( \overline{w}_{1} \wedge w_{1} \,
- \, \overline{w}_{2} \wedge w_{2} \, )|_{_{S^{-}}}\notag \\
w_{2} \, \otimes \, w_{2} \, \mapsto \,&
  \begin{pmatrix}
  0 & & 0 \\ \\
  -1 &  & 0 \end{pmatrix}  \notag \\
 \mapsto \,& -\frac{1}{2} \, \q (\overline{w}_{1} \wedge w_{2} )
 | _{_{S^{-}}} \notag
\end{align}
\end{pf}

\begin{Thm} \label{T:s+s-sym}
By using the symplectic Riesz representation, we have
\[ S^{+} \, \otimes \, S^{-} \, \overset{\cong}{\longrightarrow} \,
 \Hom ( \, S^{+} \, , \, S^{-} \, ) \]
\end{Thm}

\begin{pf}
With respect to the standard basis of $S^{+}$ and $S^{-}$, we have a map
\[ S^{+} \, \otimes \, S^{-} \, \cong \, \Hom ( \, S^{+} \, , \, S^{-} \, ) \]
by
\[  \begin{pmatrix}
a \\ \\
b \end{pmatrix} \, \otimes \,
\begin{pmatrix}
c \\ \\
d \end{pmatrix} \, \mapsto \,
 \{ \,
\begin{pmatrix}
a \\ \\
b \end{pmatrix} \,  , \cdot \, \} \, \otimes \,
\begin{pmatrix}
c \\ \\
d \end{pmatrix} \, = \,
\begin{pmatrix}
-bc & & ac   \\ \\
-bd & &  ad \end{pmatrix}
\]
where
\[  \begin{pmatrix}
a \\ \\
b \end{pmatrix} \, \in S^{+} , \quad \text{and} \quad
\begin{pmatrix}
c \\ \\
d \end{pmatrix} \, \in S^{-}  \]
and $\{ \, \cdot \, , \, \cdot \, \}$ is the symplectic form of $S^{+}$.

Explicitly, we have
\begin{align}
1_{_{C(V)}} \, \otimes \, w_{1} \, \mapsto \,&
\begin{pmatrix}
0 & & 1  \\ \\
0 & & 0 \end{pmatrix} \notag \\
1_{_{C(V)}} \, \otimes \,w_{2} \, \mapsto \,&
 \begin{pmatrix}
0 & & 0 \\ \\
0 & & 1 \end{pmatrix}  \notag\\
(w_{1} \wedge w_{2}) \, \otimes \, w_{1} \, \mapsto \,&
 \begin{pmatrix}
-1 & & 0 \\ \\
0 & & 0 \end{pmatrix}  \notag\\
(w_{1} \wedge w_{2}) \, \otimes \, w_{2} \, \mapsto \,&
 \begin{pmatrix}
0 & & 0 \\ \\
-1 & & 0 \end{pmatrix}  \notag \end{align}
This map extends to the whole $S^{+} \otimes S^{-}$ complex linearly.
\end{pf}

\begin{Cor}
 By using the symplectic Riesz representation, we have
 \[ S^{+} \, \otimes \, S^{-} \, \cong \, \Lambda ^{1}_{\CC}  \]
\end{Cor}

\begin{pf}
By {\bf Theorem \ref{T:s+s-lambda}} and {\bf Theorem \ref{T:s+s-sym} }
 we have the following
correspondence:
\begin{alignat}{4}
1_{_{C(V)}} \, &\otimes \, w_{1} \,& \mapsto \,&
  \begin{pmatrix}
  0 & & 1  \\ \\
  0 &  & 0 \end{pmatrix} \,&
 \mapsto \,&\frac{1}{2} \, e_{3} \, + \, \frac{\ii}{2} \,  e_{4}&&\,=
\,\overline{w}_{2}
 \notag \\
1_{_{C(V)}} \, &\otimes \,w_{2} \, &\mapsto \,&
  \begin{pmatrix}
  0 & & 0 \\ \\
  0 & & 1 \end{pmatrix}  \,&
 \mapsto \,&-\frac{1}{2} \, e_{1} \,  - \, \frac{\ii}{2} \,  e_{2}&&\,=
\,-\overline{w}_{1}
\notag \\
(w_{1} \wedge w_{2}) \, &\otimes \, w_{1} \,& \mapsto \,&
  \begin{pmatrix}
  -1 & & 0 \\ \\
  0 & & 0 \end{pmatrix}  \,&
 \mapsto \,&-\frac{1}{2} \, e_{1} \,  + \, \frac{\ii}{2} \,  e_{2}&&\,=
\,-w_{1}
\notag \\
(w_{1} \wedge w_{2}) \, &\otimes \,  w_{2} \, &\mapsto \,&
  \begin{pmatrix}
  0 & & 0 \\ \\
  -1 & &0 \end{pmatrix}  \,&
 \mapsto \,&-\frac{1}{2} \, e_{3} \,  + \, \frac{\ii}{2} \,  e_{4}&&\,=
\,-w_{2}
\notag
\end{alignat}
This map extends to the whole $S^{+} \otimes S^{-}$ complex linearly.
\end{pf}

Similarly, by interchanging $S^{+}$ and $S^{-}$, we get the following
\begin{Thm} \label{T:s-s+sym}
By using the symplectic Riesz representation, we have
\[ S^{-} \, \otimes \, S^{+} \, \cong \, \Hom ( \, S^{-} \, , \, S^{+} \, ) \]
\end{Thm}

\begin{pf}
With respect to the standard basis of $S^{-}$ and $S^{+}$, we have a map
\[ S^{-} \, \otimes \, S^{+} \, \overset{\cong}{\longrightarrow} \,
\Hom ( \, S^{-} \, , \, S^{+} \, ) \]
by
\[  \begin{pmatrix}
a \\ \\
b \end{pmatrix} \, \otimes \,
\begin{pmatrix}
c \\ \\
d \end{pmatrix} \, \mapsto \,
 \{ \,
\begin{pmatrix}
a \\ \\
b \end{pmatrix} \,  , \cdot \, \} \, \otimes \,
\begin{pmatrix}
c \\ \\
d \end{pmatrix} \, = \,
\begin{pmatrix}
-bc & & ac   \\ \\
-bd\, & & ad \end{pmatrix}
\]
where
\[  \begin{pmatrix}
a \\ \\
b \end{pmatrix} \, \in S^{-} , \quad \text{and} \quad
\begin{pmatrix}
c \\ \\
d \end{pmatrix} \, \in S^{+}  \]
and $\{ \, \cdot \, , \, \cdot \, \}$ is the symplectic form of $S^{-}$.

Explicitly, we have
\begin{align}
w_{1} \, \otimes \, 1_{_{C(V)}} \, \mapsto \,&
\begin{pmatrix}
0 & & 1  \\ \\
0 & & 0 \end{pmatrix} \notag \\
w_{1}\, \otimes \,(w_{1} \wedge w_{2}) \, \mapsto \,&
 \begin{pmatrix}
0 & & 0 \\ \\
0 & & 1 \end{pmatrix}  \notag\\
w_{2} \, \otimes \, 1_{_{C(V)}} \, \mapsto \,&
 \begin{pmatrix}
-1 & & 0 \\ \\
0 & & 0 \end{pmatrix}  \notag\\
w_{2} \, \otimes \, (w_{1} \wedge w_{2}) \, \mapsto \,&
 \begin{pmatrix}
0 & & 0 \\ \\
-1 & & 0 \end{pmatrix}  \notag \end{align}
This map extends to the whole $S^{-} \otimes S^{+}$ complex linearly.
\end{pf}

 \begin{Cor}
 By using the symplectic Riesz representation, we have
 \[ S^{-} \, \otimes \, S^{+} \, \cong \, \Lambda ^{1}_{\CC}  \]
\end{Cor}

\begin{pf}
By {\bf Theorem \ref{T:s+s-lambda}} and {\bf Theorem \ref{T:s-s+sym} }
we have the following
correspondence:
\begin{alignat}{4}
w_{1} \, &\otimes \, 1_{_{C(V)}} &\, \mapsto \,&
  \begin{pmatrix}
  0 & & 1  \\ \\
  0 &  & 0 \end{pmatrix} \, &&
 \mapsto \,-\frac{1}{2} \, e_{3} \,  - \, \frac{\ii}{2} \,  e_{4} &&\,=
\,-\overline{w}_{2}
\notag \\
w_{1} \, &\otimes \,(w_{1} \wedge w_{2}) &\, \mapsto \,&
  \begin{pmatrix}
  0 & & 0 \\ \\
  0 & & 1 \end{pmatrix}  \, &&
 \mapsto \,\frac{1}{2} \, e_{1} \,  - \, \frac{\ii}{2} \,  e_{2} &&\,=\,w_{1}
\notag \\
 w_{2} \, &\otimes \, 1_{_{C(V)}}& \, \mapsto \,&
  \begin{pmatrix}
  -1 & & 0 \\ \\
  0 & & 0 \end{pmatrix}  \,&&
 \mapsto \,\frac{1}{2} \, e_{1} \, + \, \frac{\ii}{2} \,  e_{2} &&\,=
\,\overline{w}_{1}
\notag \\
 w_{2} \, &\otimes \, (w_{1} \wedge w_{2}) & \, \mapsto \,&
  \begin{pmatrix}
  0 & & 0 \\ \\
  -1 & & 0 \end{pmatrix}  \,&&
 \mapsto \,\frac{1}{2} \, e_{3} \, - \, \frac{\ii}{2} \,  e_{4} &&\,=\,w_{2}
\notag
\end{alignat}
This map extends to the whole $S^{-} \otimes S^{+}$ complex linearly.
\end{pf}

\bigskip
\subsection{Quaternionic Structure on the Spinors}
Now we shall combine the Hermitian and symplectic structures
 together.
\begin{Lem}
We have an isomorphism
\[ \jj ^{\ast} : S^{+\ast} \overset{\cong}{\longrightarrow} S^{+\ast} \]
given by the following identification
\[ \la
 \begin{pmatrix}
s^{+}_{1}\\ \\
s^{+}_{2}\end{pmatrix}
\, , \, \cdot \ra \, \longmapsto \,
 \begin{pmatrix}
\overline{s}^{+}_{1}\\ \\
\overline{s}^{+}_{2}\end{pmatrix}^{\ast} \, \longmapsto \,
 \{
 \begin{pmatrix}
-\overline{s}^{+}_{2}\\ \\
\overline{s}^{+}_{1}\end{pmatrix}
\, , \, \cdot \}  \]
where the first isomorphism is from the Hermition Riesz representation
 and the second is from
the symplectic Riesz representation.
\end{Lem}

\begin{Cor}
We have an {\bf anti-linear} isomorphism
\[ \jj : S^{+} \longrightarrow S^{+} \]
which is the composition of the following maps
\[  \begin{pmatrix}
s^{+}_{1}\\ \\
s^{+}_{2}\end{pmatrix} \, \underset{\text{Riesz}}
{\overset{\text{Hermitian}}{\longmapsto}} \, \la \,
 \begin{pmatrix}
s^{+}_{1}\\ \\
s^{+}_{2}\end{pmatrix}\, , \, \cdot \, \ra \,
\overset{\jj ^{\ast}}{\longmapsto} \, \{ \,
 \begin{pmatrix}
-\overline{s}^{+}_{2}\\ \\
\overline{s}^{+}_{1}\end{pmatrix} \, , \, \cdot \, \} \,
\underset{\text{Riesz}}{\overset{\text{symplectic}}{\longmapsto}} \,
 \begin{pmatrix}
-\overline{s}^{+}_{2}\\ \\
\overline{s}^{+}_{1}\end{pmatrix} \]
such that the following diagram commutes:
\[
\begin{CD}
S^{+}                        @>\jj>>                    S^{+} \\
@V\text{Hermitian Riesz}V{\cong}V @V{\cong}V\text{symplectic Riesz}V \\
S^{+\ast}              @>>{\jj^{\ast}}>                 S^{+\ast}
\end{CD}
\]
\end{Cor}

\begin{pf}
We only need to prove the anti-linearity of $\jj$
which is straight forward from the following
computation:
\[  \jj \, ( \, c \, \begin{pmatrix}
 s^{+}_{1}\\ \\
s^{+}_{2}\end{pmatrix}\, ) \, =\,\jj\,
\begin{pmatrix}
c \, s^{+}_{1}\\ \\
c \, s^{+}_{2}\end{pmatrix}\,
=\,
\begin{pmatrix}
-\overline{c} \, \overline{s}^{+}_{2}\\ \\
\overline{c} \, \overline{s}^{+}_{1}\end{pmatrix}\,
=\,\overline{c}\,
\begin{pmatrix}
- \overline{s}^{+}_{2}\\ \\
  \overline{s}^{+}_{1}\end{pmatrix}\,
=\,\overline{c}\, ( \, \jj \,
\begin{pmatrix}
s^{+}_{1}\\ \\
s^{+}_{2}\end{pmatrix} \, ) \]
\end{pf}

\begin{Thm}
The isomorphism $\jj$ above and the isomorphism
\[ \ii : S^{+} \longrightarrow S^{+} \]
\[ \ii \,  \begin{pmatrix}
s^{+}_{1}\\ \\
s^{+}_{2}\end{pmatrix}\,= \, \begin{pmatrix}
\ii \,s^{+}_{1}\\ \\
\ii \,s^{+}_{2}\end{pmatrix}  \]
together gives a {\bf quaternionic structure} on $S^{+}$.
\end{Thm}

\begin{pf}
It is obvious that $ \ii ^{2} = -1_{_{\End ( S^{+})}}$. Now we also have
\[ \jj ^{2} \, \begin{pmatrix}
s^{+}_{1}\\ \\
s^{+}_{2}\end{pmatrix}\,= \, \jj \, \begin{pmatrix}
-\overline{s}^{+}_{2}\\ \\
\overline{s}^{+}_{1}\end{pmatrix}\,= \, \begin{pmatrix}
-\overline{\overline{s}}^{+}_{1}\\ \\
-\overline{\overline{s}}^{+}_{2}\end{pmatrix}\,= \, - \begin{pmatrix}
s^{+}_{1}\\ \\
s^{+}_{2}\end{pmatrix}\]

Define
\[ \kk \, = \, \ii \jj \, : S^{+} \longrightarrow S^{+} \]
therefore
\[ \kk \, \begin{pmatrix}
s^{+}_{1}\\ \\
s^{+}_{2}\end{pmatrix}\,= \, \ii \jj \,\begin{pmatrix}
s^{+}_{1}\\ \\
s^{+}_{2}\end{pmatrix}\, = \,\ii \, \begin{pmatrix}
-\overline{s}^{+}_{2}\\ \\
\overline{s}^{+}_{1}\end{pmatrix}\,= \begin{pmatrix}
-\ii \, \overline{s}^{+}_{2}\\ \\
\ii \,\overline{s}^{+}_{1}\end{pmatrix} \]
Hence we have
\[ \kk ^{2} \, \begin{pmatrix}
s^{+}_{1}\\ \\
s^{+}_{2}\end{pmatrix}\,= \, \kk\, \begin{pmatrix}
-\ii \, \overline{s}^{+}_{2}\\ \\
\ii \,\overline{s}^{+}_{1}\end{pmatrix}\, = \, \begin{pmatrix}
-\ii \, (\overline{\ii \,\overline{s}^{+}_{1}})\\ \\
\ii \,(\overline{-\ii \, \overline{s}^{+}_{2}})\end{pmatrix}\, = \,
\begin{pmatrix}
- s^{+}_{1}\\ \\
-s^{+}_{2}\end{pmatrix}\, = \, - \, \begin{pmatrix}
s^{+}_{1}\\ \\
s^{+}_{2}\end{pmatrix} \]
and
\[\jj \ii \,\begin{pmatrix}
s^{+}_{1}\\ \\
s^{+}_{2}\end{pmatrix}\, = \, \jj \,\begin{pmatrix}
\ii \,s^{+}_{1}\\ \\
\ii \,s^{+}_{2}\end{pmatrix}\, = \, \begin{pmatrix}
\ii \,\overline{s}^{+}_{2}\\ \\
-\ii \,\overline{s}^{+}_{1}\end{pmatrix}\, = \,- \kk \,\begin{pmatrix}
s^{+}_{1}\\ \\
s^{+}_{2}\end{pmatrix} \]
Hence we have a quaternionic structure induced by $\ii$ and $\jj$.
\end{pf}

\end{document}